\documentclass[11pt]{article}
\usepackage[utf8]{inputenc}
\usepackage{amssymb,amsmath,amsthm,mathabx, graphicx,geometry,natbib,hyperref,xcolor,bm,mathtools,authblk,caption,subcaption, mathrsfs}

\newgeometry{left=1.1in,right=1.1in,top=1in,bottom=1.1in}

\hypersetup{
    colorlinks,
    linkcolor={blue},
    citecolor={blue},
    urlcolor={blue}
}
\newtheorem{theorem}{Theorem}

\newtheorem{corollary}{Corollary}
\newtheorem{lemma}{Lemma}

\newtheorem*{assumption*}{Assumption}

\theoremstyle{definition}
\newtheorem{example}{Example}
\newtheorem{definition}{Definition}

\newtheorem{innercustomthm}{Theorem}
\newenvironment{customthm}[1]
  {\renewcommand\theinnercustomthm{#1}\innercustomthm}
  {\endinnercustomthm}

\newtheorem{innercustomlem}{Lemma}
\newenvironment{customlem}[1]
  {\renewcommand\theinnercustomlem{#1}\innercustomlem}
  {\endinnercustomlem}

\newtheorem{innercustomcor}{Corollary}
\newenvironment{customcor}[1]
  {\renewcommand\theinnercustomcor{#1}\innercustomcor}
  {\endinnercustomcor}

  \newtheorem{innercustomexample}{Example}
\newenvironment{customexample}[1]
  {\renewcommand\theinnercustomexample{#1}\innercustomexample}
  {\endinnercustomexample}

\newcommand{\bl}[1]{{\mathbf #1}}
\newcommand{\bs}[1]{\boldsymbol #1}
\newcommand{\Pa}{\mathcal{P}_A}
\newcommand{\Po}{\mathcal{P}_0}

\newcommand{\mb}[1]{\mathbb #1}
\newcommand{\mc}[1]{\mathcal #1}

\begin{document}

\title{On the Existence of Unbiased Hypothesis Tests: An Algebraic Approach}
\author{Andrew McCormack}
\affil{Department of Mathematical and Statistical Sciences, University of Alberta}
\date{\today}

\maketitle

\begin{abstract}
In hypothesis testing problems the property of strict unbiasedness describes whether a test is able to discriminate, in the sense of a difference in power, between any distribution in the null hypothesis space and any distribution in the alternative hypothesis space. In this work we examine conditions under which unbiased tests exist for discrete statistical models. It is shown that the existence of an unbiased test can be reduced to an algebraic criterion; an unbiased test exists if and only if there exists a polynomial that separates the null and alternative hypothesis sets. This places a strong, semialgebraic restriction on the classes of null hypotheses that have unbiased tests. The minimum degree of a separating polynomial coincides with the minimum sample size that is needed for an unbiased test to exist, termed the unbiasedness threshold. It is demonstrated that Gr\"obner basis techniques can be used to provide upper bounds for, and in many cases exactly find, the unbiasedness threshold. Existence questions for uniformly most powerful unbiased tests are also addressed, where it is shown that whether such a test exists can depend subtly on the specified level of the test and the sample size. Numerous examples, concerning tests in contingency tables, linear, log-linear, and mixture models are provided. All of the machinery developed in this work is constructive in the sense that when a test with a certain property is shown to exist it is possible to explicitly construct this test.

\end{abstract}

\section{Introduction}

An enduring problem in statistics is the the determination of an appropriate sample size that is required in order for the performance of an inferential procedure, as measured by the power of a test \cite{SampleSizePoissonLogistic,signorini1991samplePoisson}, the expected width of a confidence interval \cite{GraybillSampleSizeWidth}, or the mean-squared error of an estimator \cite{adcock1997sample}, to be adequate for a given application. Questions of sample size determination commonly occur in the fields of experimental design and survey sampling where an investigator determines when a question of interest can feasibly be answered using the sample size available to the investigator \cite{ClinicalTrialSampleSize}. An even more fundamental question is whether an inferential procedure with specified, desirable properties exists at all \citep{ExistenceofUMPUTests,ExistenceofUMPSTests,NonExistenceofUnbiasedEstimators}.

In this article, we focus on the setting of categorical data and provide precise conditions on the classes of null hypotheses that possess unbiased (UB) tests. Here we are only interested in unbiased tests that have non-trivial power. Specifically, a null hypothesis has a non-trivial unbiased test if and only if it can be represented as the sub-level set of a \textit{single} polynomial. Simple models will be exhibited that violate this condition. Conversely, a large portion of null hypothesis sets are defined as the vanishing set of a collection of polynomials \citep[Ch 2]{LecturesonAlgebraicDrtonSullivant}. The results in this work show that for this important class of hypothesis there always exists a non-trivial unbiased test.

When an unbiased test does exist, the requisite minimum sample size for there to exist an unbiased test, termed the \textit{unbiasedness threshold}, is studied. The techniques developed for finding bounds of the unbiasedness threshold can be used to construct unbiased tests for arbitrarily complex null hypotheses that are described by polynomial constraints. Furthermore, in small-sample settings it is shown how uniformly most powerful unbiased (UMPU) tests can be constructed when they exist.  

Applied algebraic geometry \cite{cox2013ideals}, which can approximately be thought of as the manipulation of polynomial systems of equations, is the primary technical machinery that it utilized in this article. The burgeoning field of algebraic statistics uses such gadgetry to address statistical questions \cite{AlgStatPistoneWynn,SullivantAlgStats}. The seminal paper \citep{DiaconisSturmfelsMarkov} introduced the concept of a Markov basis as a way to sample from discrete distributions conditional on the value of a sufficient statistic, partially as a means of sampling from null distributions for hypothesis tests, such as Fisher's exact test. Subsequent work on Markov bases has studied how algebraic methods can be used to construct hypothesis tests for increasingly complex models \citep{MarkovBasesBook, DobraMarkovBasesGraphModel, GrossPetrovikGoodnessFitMarkovBasis, MarkovBasesHierarcSullivantHosten}. Other work on hypothesis testing in algebraic statistics has focused on constructing tests for statistical models that have singularities \citep{DrtonLRTSingularities,DrtonWaldTestSingularities,NilsUStatisticSingularity}. For example, Watanabe's singular learning theory \citep{WatanabeBook} provides methods for performing model selection for models with singularities.

The strand of research that is most closely linked to the present work is on finding the \textit{maximum likelihood (ML) threshold}, namely the number of samples needed for the maximum likelihood estimator to exist almost surely. In \cite{DerksenKron, DrtonKurikiKron,GrossSullivantMLETbhreshold, UhlerMLEGaussianGraph} the maximum likelihood threshold is described for different exponential families, such as structured covariance models and contingency tables. For exponential families with polynomial parameterizations finding the maximum likelihood estimator can be framed algebraically, as the likelihood equations usually involve solving a rational system of equations. The unbiasedness threshold introduced here can be seen as a testing analogue of the maximum likelihood threshold. Both thresholds provide a measure of the baseline complexity of the respective point estimation and hypothesis testing tasks. Unless prior information is present that discounts a portion of the alternative hypothesis space, if no non-trivial unbiased test exists, a statistician may need to either revise their null hypothesis, or obtain a larger sample size. As is the case for ML threshold, the UB threshold is typically most relevant in sparse settings where there are many parameters to be tested and the multinomial counts within categories is low. 

The scope of this paper is constrained to hypothesis tests for multinomial, categorical data, with a focus on tests for contingency tables. Notable related work examines tests for independence \citep{amiri2017comparisonofindeptests}, conditional independence \citep{BirchConditionalIndependenceTests, ExactConditionalTestsKreiner,TomonariHierarchicalSubspace}, symmetry \citep{bhapkar1979Symmetrytests}, and general linear hypotheses in the log-odds. A significant portion of the literature is focused on asymptotic tests like the chi-squared and likelihood ratio tests \citep{christensenLogLinearBook}. Exact tests that appear in the literature, like Fisher's exact test for independence, are usually based off the general form for the uniformly most powerful unbiased (UMPU) test of a linear hypothesis in a regular exponential family \citep[Thm 4.4.1] {TSPLehmannRomano}. Interestingly, in certain cases it will be shown that this UMPU test has no power for all sample sizes. We also highlight specific curved exponential family submodels \citep{BrownExponentialFamilies} that possess non-trivial UMPU tests, a setting that to the best of the author's knowledge has not been studied in generality.

A final contribution of this work is to introduce the method of power polynomials for constructing exact tests for complex hypotheses. The perspective taken here is to manipulate the power function, which in the discrete setting is a polynomial, of a test directly rather than search for a test statistic that induces a power function with desirable properties. While a general likelihood ratio test is asymptotically valid and efficient when the model is smooth, the size of the likelihood ratio test can be distorted for medium-to-small sample sizes and for models with singularities \citep{DrtonLRTSingularities}. In contrast, the proposed methodology results in exact tests that can be can be tuned to distribute power in regions of the alternative hypothesis that are anticipated to contain the true distribution \citep{PeterHoffPvalueIndirectInformation}. Working directly with power polynomials is best suited for small sample size regimes as the degree of the polynomial that has to be optimized increases with the sample size.  
% Symmetry and invariance considerations are shown to be helpful in lessening the computational overhead of searching for a power polynomial; invariance of a test under a subgroup of the permutation group is equivalent to constraining the power polynomial be a symmetric polynomial.
For exact testing methods \citep{agresti2001exact} these computational hurdles are to be expected, since there is a greater number of configurations of counts that could potentially occur for large sample sizes. For example, when using MCMC to find a null distribution \citep{smith1996montecarloexact} the Markov chain has to traverse a larger number of states when the sample size is large.  

An overview of this article is as follows: Section \ref{sec:BackgroundonTesting} includes background and central definitions regarding hypothesis testing. The fundamental correspondence between tests and power polynomials is detailed in Section \ref{sec:PowerPolynomialCorrespondence}. The subsequent section provides an exact characterization of which null hypotheses admit non-trivial unbiased tests in multinomial models. Subverting expectations, certain mixture models with fixed mixture components and certain log-linear models are shown to not admit non-trivial unbiased tests. A review of core concepts in applied algebraic geometry needed for the subsequent sections, such as ideals and varieties, is provided in Section \ref{sec:AlgebraBackground}. Algebraic methods for computing the unbiasedness threshold are presented in Section \ref{sec:UBThreshold}. For most null hypotheses of interest, we provide a general form for power polynomials associated with unbiased tests. The final section considers existence questions for UMPU tests geometrically in terms of a polytope. Examples are provided of non-exponential families that posses UMPU tests and UMPU tests that exist for small, but not large, sample sizes.
% Also explicated in this section is how permutation invariance can be paired with unbiasedness to simplify the space of power polynomials under consideration. 
Proofs of all of the results in this article are provided in the supplementary material.

\subsection{Notation}
The notation introduced in this article is summarized here for reference. Statistical quantities appearing in this work are: $\mc{X}$ for the sample space, $X_i$ for an observation in $\mc{X}$, $\theta$ for a general parameter in the parameter space $\Theta$, $\Po$ and $\Pa$ for null and alternative hypothesis sets, $\phi$ for a hypothesis test, $\alpha$ for the level or size of the test, $\beta_\phi$ for the power function of $\phi$, $\Tilde{\beta}$ for a separating polynomial, $\bs{\pi}$ for a probability vector in the simplex $\Delta_{k-1}$, and $\Tilde{\Delta}_{k-1}$ for the simplex projected onto its first $k-1$ coordinates. Sample sizes are always denoted by $n$ and $k$ always denotes the dimension of the multinomial model under consideration. The bracket notation $[n]$ represents the set of $n$ elements $\{1,\ldots,n\}$. Vectors and matrices appear in boldfont. The set of polynomials over a field $\mb{K}$, where $\mb{K}$ is either $\mb{R}$ or $\mb{C}$, in the free variables $z_1,\ldots,z_k$ is $\mb{K}[z_1,\ldots,z_k]$, or identically $\mb{K}[\bl{z}]$. A degree $n$ polynomial can be represented in the shorthand $f(\bl{z}) = \sum_{I} f_I \bl{z}^I$ where $\bl{z}^I \coloneqq \prod_{j = 1}^k z_j^{i_j}$ and the $f_I$ are scalars. Multiindices are represented by $I,J,L$ and consist of $k$ natural numbers that have a sum equal to $n$, or depending on context, bounded above by $n$. The variety associated with a collection of polynomials $W$ is $V_{\mb{R}^k}(W)$ and the ideal in $\mb{R}[\pi_1,\ldots,\pi_k]$ generated by a set $V \subseteq \mb{R}^k$ is $I_{\mb{R}^k}(V)$. If the subscript $\mb{R}^{k}$ is replaced with $\mb{R}^{k-1}$ the corresponding varieties and ideals are viewed as subsets of $\text{aff}(\Delta_{k-1}) \cong \mb{R}^{k-1}$ and $\mb{R}[\pi_1,\ldots,\pi_{k-1}]$ respectively. The ideal generated by the polynomials $f_{(1)},\ldots,f_{(m)}$ is denoted by $\langle f_{(1)},\ldots,f_{(m)} \rangle$. Subscripts in parentheses are used to denote the index of a polynomial within a collection, while subscripts not in parentheses represent the scalar coefficients of a polynomial. For example, $f_{2,(3)}$ represents the second coefficient in the polynomial $f_{(3)}$. All topological concepts, such as interiors, closures, and boundaries, refer to the standard Euclidean topology, not the Zariski topology.

\section{Background On Unbiased Tests}
\label{sec:BackgroundonTesting}
The present framework for hypothesis testing involves a statistical model $\mc{P} = \{P_\theta: \theta \in \Theta\}$ of probability distributions over a sample space $\mc{X}$, where it is desired to test the validity of the null hypothesis $\Po = \{P_\theta:\theta \in \Theta_0\} \subset \mc{P}$ against an alternative hypothesis $\Pa = \{P_\theta:\theta \in \Theta_A\} \subset \mc{P}$, with  $\mc{P} = \Po \bigsqcup \Pa$. That is, some data $X \sim P_\theta$ is observed and the goal is to make a decision as to whether $P_\theta$ is in $\Po$ or $\Pa$. This decision comes in the form of a randomized hypothesis test $\phi:\mc{X} \rightarrow [0,1]$ where a decision is made to reject the null hypothesis that $P_\theta \in \Po$ by:
\begin{enumerate}
    \item First computing $\phi(X) \in [0,1]$.
    \item Flipping a weighted coin that has a probability of $\phi(X)$ of landing heads, or equivalently independently drawing a $Z \sim \text{Bernoulli}(\phi(X))$ random variable.
    \item Rejecting the null hypothesis if the coin lands heads, or equivalently if $Z = 1$. 
\end{enumerate}
Non-randomized hypothesis tests are the subset of randomized tests of the form $\phi: \mc{X} \rightarrow \{0,1\}$, where no auxiliary randomization step is required. All hypothesis tests in this work will be assumed to be randomized unless otherwise stated. 

The performance of a randomized test is encapsulated by its power function.
\begin{definition}[Power function]
    The power function $\beta_\phi: \Theta \rightarrow [0,1]$ of the randomized test $\phi$ is equal to 
    \begin{align}
        \beta_\phi(\theta) = E_\theta\big(\phi(X)\big) =  \int_{\mc{X}} \phi(x) P_\theta(dx).
    \end{align}
\end{definition}
The quantity $\beta_\phi(\theta)$ is equal to the marginal probability of rejecting $\Po$ under the supposition that $X \sim P_\theta$. Committing a type-I error by rejecting the null hypothesis when it is true is often seen as a more serious error than failing to reject the null when it is false. As such, the maximum type-I error rate, specifically the size of the test $\sup_{\theta \in \Theta_0} \beta_\phi(\theta)$, is typically bounded above by some level $\alpha \in (0,1)$. Tests that have size less than or equal to $\alpha$ are said to be level-$\alpha$ tests. It is also desirable to have a test that rejects $\Po$ with high probability when the null hypothesis does not hold. The property of unbiasedness requires that the probability of rejecting $\Po$ when the alternative hypothesis holds is at least as large as the probability of rejecting $\Po$ when the null hypothesis holds. 

\begin{definition}[Non-trivial and strictly unbiased tests]
    A test $\phi$ is unbiased if its power function has the property that
    \begin{align}
    \label{eqn:UBDefinition}
        \beta_\phi(\theta_0) \leq \beta_\phi(\theta_A) \;\; \text{ for all } \theta_0 \in \Theta_0  \text{ and } \theta_A \in \Theta_A.
    \end{align}
    A test is said to be a non-trivial unbiased \textbf{(NTUB)} test if $\beta_\phi$ is not a constant function. A test is said to be strictly unbiased \textbf{(SUB)} if \eqref{eqn:UBDefinition} holds with strict inequality for all $\theta_0$ and $\theta_A$.
\end{definition}
The question of whether there exists an unbiased test is not an interesting one because the constant test function $\phi(x) = c$ is unbiased as it has the trivial power function $\beta_\phi(\theta) = c$. However, constant test functions are not of practical use as they completely disregard the observed data and amount to indiscriminate coin flipping. The focus of this work will therefore be on the existence of non-trivial unbiased and strictly unbiased tests, which together will be referred to as the existence of unbiased tests. Despite the constant test not being of use on its own, it is useful when combined with other tests as it can be used to modulate the size of a test, as is done in the proof in the following lemma.

\begin{lemma}
\label{lem:ExistenceLevel}
    The four different sets of level-$\alpha$, unbiased, non-trivial unbiased, and strictly unbiased tests are all convex sets of test functions. If there exists a size-$\alpha_0 \in [0,1)$ NTUB or SUB test, then for any $\alpha \in (0,1)$ there exists a size-$\alpha$ NTUB or SUB test respectively.   
\end{lemma}

Unbiased tests need not have favorable power attributes against $\Pa$. A much stronger condition is that a test be uniformly most powerful unbiased.

\begin{definition}[Uniformly most powerful unbiased test]
    A level-$\alpha$ test $\phi$ is a uniformly most powerful unbiased test (\textbf{UMPU}) if it is unbiased and if for every other unbiased, level-$\alpha$ test $\Tilde{\phi}$
    \begin{align*}
         \beta_{\Tilde{\phi}}(\theta_A) \leq \beta_\phi(\theta_A) ,\;\; \text{ for all } \theta_A \in \Theta_A.
    \end{align*}
\end{definition}

Unlike the existence question for unbiased tests, the existence of a UMPU test may depend on the level of the test that is required as the following example illustrates.

\begin{example}
    Consider the finite sample space $\mc{X} = \{a,b,c\}$, the null hypothesis $\Po = \{\theta_1\}$ and the alternative hypothesis $\Pa = \{\theta_1,\theta_2\}$.  The probability mass function of each $P_{\theta_i}$ is defined as
\begin{center}
\begin{tabular}{||c| c c c||} 
 \hline
 $\mc{X}$ & $a$ & $b$ & $c$ \\ 
 \hline\hline
 $P_{\theta_1}$ & 1/3 & 1/3 & 1/3 \\ 
 \hline
 $P_{\theta_2}$ & 3/4 & 1/4 & 0 \\
 \hline
 $P_{\theta_3}$ & 3/4 & 0 &  1/4
 \\
 \hline
\end{tabular}
\end{center}
\end{example}
The test that rejects when $a$ is observed is a UMPU test with size equal to $1/3$. There is no UMPU test for sizes larger than $1/3$.

However, we do have the following result that shows that UMPU tests are downward-closed with respect to the test size.
\begin{lemma}
\label{lem:UMPUvsSize}
    If there exists a size-$\alpha_0$ UMPU test of $\Po$ against $\Pa$ for some $\alpha_0 \in [0,1)$ then there also exists a size-$\alpha$ UMPU test for all $\alpha \in [0,\alpha_0]$.  
\end{lemma}

% \begin{proof}
%     Convexity for level-$\alpha$, UB and SUB tests follows easily by the linearity of the power function $\beta_{\lambda \phi_1 + (1-\lambda)\phi_2} = \lambda \beta_{\phi_1} + (1-\lambda)\beta_{\phi_2}$. Suppose $\phi_1,\phi_2$ are NTUB tests. There necessarily exists points $\theta_0 \in \Theta_0$ and $\theta_A$ where $\beta_{\phi_1}(\theta_0) < \beta_{\phi_1}(\theta_A)$ as otherwise $\beta_{\phi_1}$ would be constant. It follows that $ \lambda \beta_{\phi_1}(\theta_0) + (1-\lambda)\beta_{\phi_2}(\theta_0) <  \lambda \beta_{\phi_1}(\theta_A) + (1-\lambda)\beta_{\phi_2}(\theta_A)$, implying that $\beta_{\lambda \phi_1 + (1-\lambda)\phi_2}$ is not constant for $\lambda \in (0,1]$. By symmetry this holds for $\lambda \in [0,1)$.
    
%     To prove the last statement, we will take a convex combination of a size $\alpha_0$ NTUB or SUB test $\phi$ with either the constant test $1$ or the constant test $0$. If $\alpha > \alpha_0$ consider the test $\lambda \phi + (1-\lambda)$ that is NTUB or SUB when $\phi$ is, by the proof above. The size of this test is $\text{Size}(\lambda \phi + (1-\lambda)) = \lambda\text{Size}(\phi) + (1-\lambda) = \lambda \alpha_0 +(1-\lambda)$. Taking $\lambda = \tfrac{1-\alpha}{1-\alpha_0} \in (0,1)$ gives the desired test. The case where $\alpha < \alpha_0$ is handled similarly by considering $\tfrac{\alpha}{\alpha_0}\phi$.  
% \end{proof}

\section{Power Polynomials in Discrete Models}
\label{sec:PowerPolynomialCorrespondence}
It will be assumed throughout this paper that $\mc{P} = \{ \text{Multinomial}_k(\bs{\pi},n): \bs{\pi} \in \Delta_{k-1}\}$ consists of all multinomial distributions of counts on $k$ different categories that sum to $n$. Null and alternative hypothesis sets can easily be specified by a partition of the $(k-1)$-dimensional probability simplex $\Delta_{k-1} = \{\bs{\pi} \in \mb{R}^k: \bs{\pi} \geq \bl{0}, \bl{1}_k^\intercal \bs{\pi} = 1\}$. We let $\Po \subseteq \Delta_{k-1}$ and $\Pa = \Delta_{k-1} \backslash \Po$ represent subsets of the probability simplex as well as collections of multinomial distributions, where the intended meaning will be clear from context. Null hypothesis sets will always assumed to be non-trivial with $\emptyset \subsetneq \Po \subsetneq \Delta_{k-1}$. 

At a high-level, the motivation for considering models for discrete data is that the space of all test functions is finite dimensional as there are only finitely many possible outcomes of counts. Specifically, if $\bl{X} \sim \text{Multinomial}_k(\bs{\pi},n)$ then $\bl{X}$ can take on ${ n+k - 1 \choose k - 1}$ different values, and so by labeling these values by the numbers $1,\ldots,{n+k-1 \choose k - 1}$ a test function can be viewed as a function $\phi:\{1,\ldots,{n + k-1 \choose k - 1}\} \rightarrow [0,1]$, which is a point in the cube $[0,1]^{{n + k-1 \choose k - 1}}$. The dimension of this cube grows rapidly for increasing sample sizes at a rate of $O(n^{k-1})$, necessitating the imposition of further restrictions on the test function when searching for a viable test. In contrast, for continuous data with a sample space of $\mc{X} = \mb{R}$ and a univariate Gaussian model $\mc{P} = \{\mc{N}(\mu,\sigma^2): \mu \in \mb{R}, \; \sigma^2 > 0\}$, the space of test functions is a convex set with non-empty interior in the $L^\infty$ space of functions that are bounded almost everywhere with respect to the Lebesgue measure.

An even more important motivation for considering discrete models is that the power function is a multivariate polynomial in $\bs{\pi}$:
\begin{align}
\label{eqn:PowerPolynomial}
    \beta_\phi(\bs{\pi}) = E_{\bs{\pi}}\big(\phi(\bl{x})\big)  = \sum_{\bl{x}: \bl{x} \in \mb{N}^k, \bl{1}_k^\intercal \bl{x} = n} \phi(\bl{x}) { n \choose \bl{x}} \prod_{i = 1}^k \pi_i^{x_i} = \sum_{\bl{x}: \bl{x} \in \mb{N}^k, \bl{1}_k^\intercal \bl{x} = n} \phi(\bl{x})  { n \choose \bl{x}}\bs{\pi}^{\bl{x}}.
\end{align}
The set of polynomials with real or complex  coefficients in the variables $\pi_1,\ldots,\pi_k$ is denoted by $\mb{R}[\pi_1,\ldots,\pi_k]$ or $\mb{C}[\pi_1,\ldots,\pi_n]$ respectively. The above power polynomial is a homogeneous, degree $n$ polynomial, meaning that the total degree of each of the monomial terms $\bs{\pi}^{\bl{x}} \coloneqq  \prod_{i = 1}^k \pi_i^{x_i}$ is equal to $n$ because $\sum_{i = 1}^k x_i = n$. Rather than dealing with test functions $\phi$, we primarily work directly in the space of polynomial power functions. The following consequence of \eqref{eqn:PowerPolynomial} shows that nothing is lost in doing so.

\begin{lemma}
\label{lem:PowerFunctionCharacterization}
    The set of all possible power functions $\{\beta_\phi: \phi \text{ is a test function}\}$ is equal to the set of homogeneous, degree $n$ polynomials of the form $\{ \sum c_{\bl{x}} \bs{\pi}^\bl{x}: c_{\bl{x}} \in [0, {n \choose \bl{x}}]\}$. The power functions of non-randomized tests are likewise given by 
    \begin{align*}
        \big\{\beta_\phi: \phi \text{ is a non-randomized test function} \big\} = \bigg\{ \sum c_{\bl{x}} \bs{\pi}^\bl{x}: c_{\bl{x}} \in \big\{0, {n \choose \bl{x}}\big\} \bigg\}.
    \end{align*}
    The linear map $\phi \mapsto \beta_\phi$ is injective, so that there exists a unique test function corresponding to each power polynomial. 
\end{lemma}
Any polynomial $\beta$ that satisfies the constraints in Lemma \ref{lem:PowerFunctionCharacterization} will be referred to as a \textit{power polynomial}. One simple method for recovering a test $\phi$ given its associated power polynomial $\beta_\phi$ in \eqref{eqn:PowerPolynomial} is to differentiate the power polynomial: 
\begin{align*}
    \phi(\bl{x}) = \frac{1}{n!} \frac{\partial^n}{\partial_{\pi_1}^{x_1}\cdots \partial^{x_n}_{\pi_n}} \beta_\phi(\bs{\pi}).
\end{align*}

\subsection{Power Polynomials for Other Discrete Models}
\label{sec:PowerPolynomialOtherSampling}
Polynomials occur frequently when modeling discrete data with a notable class of examples coming from power series distributions \citep[Sec 1.4]{HirjiExactAnalysis} that include the multinomial, negative binomial, and Poisson distributions. In principle, the multinomial model is a completely general model for count data with the number of counts fixed. However, in some settings, such as sampling a contingency table with a fixed margin, the multinomial model is not the most natural model to use. A $p \times  q$ contingency table with  fixed row counts $n_1,\ldots,n_p$ is most naturally viewed as a product of independent multinomial models $\times_{i = 1}^p \text{Multinomial}_{q}(\bs{\pi}^{(i)},n_i)$ with each $\bs{\pi}^{(i)} \in \Delta_{q-1}$. To view this as a submodel of a multinomial model, a separate multinomial category is introduced for every unique arrangement of the entire table that has the appropriate row sums. For instance, in a $2 \times 2$ table with row counts constrained to $n_1 = n_2 = 2$ the multinomial distribution would have $9$ different categories with a category corresponding to each of the tables $\bl{T}$
\begin{align*}
    \begin{bmatrix}
        2 & 0
        \\
        2 & 0
    \end{bmatrix},  \begin{bmatrix}
        1 & 1
        \\
        2 & 0
    \end{bmatrix},     \begin{bmatrix}
        0 & 2
        \\
        2 & 0
    \end{bmatrix},    \begin{bmatrix}
        2 & 0
        \\
        1 & 1
    \end{bmatrix},  \begin{bmatrix}
        1 & 1
        \\
        1 & 1
    \end{bmatrix},     \begin{bmatrix}
        0 & 2
        \\
        1 & 1
    \end{bmatrix},    \begin{bmatrix}
        2 & 0
        \\
        0 & 2
    \end{bmatrix},  \begin{bmatrix}
        1 & 1
        \\
        0 & 2
    \end{bmatrix},     \begin{bmatrix}
        0 & 2
        \\
        0 & 2
    \end{bmatrix}.
\end{align*}
The probabilities for each of the above tables are further constrained to have the form
\begin{align*}
    \pi_{\bl{T}} = {2 \choose T_{11}, T_{12}} {2 \choose T_{21}, T_{22}} \prod_{i = 1}^2 \prod_{j = 1}^2 \tau_{ij}^{T_{ij}},
\end{align*}
for some $\tau_{ij}$ parameters.
In this way, the above fixed row count model $\mc{R} \subsetneq \Delta_{8}$ can be viewed as a submodel of a full multinomial model.  The theory presented in the subsequent sections therefore largely carries over to such restricted sampling settings. One notable distinction is that the relevant alternative hypothesis for a null hypothesis with $\Po \subset \mc{R} \subsetneq \Delta_{k'-1}$ for some $k'$, is usually not the complement $\Delta_{k'-1}\backslash\Po$ under the full multinomial model, but rather is $\mc{R} \backslash \Po$. It follows the property of a test being unbiased is a weaker condition when the alternative hypothesis set is restricted to $\mc{R} \backslash \Po$, as the constraint \eqref{eqn:UBDefinition} must hold only over the set $\mc{R}\backslash \Po$ that is smaller than $\Delta_{k'-1} \backslash \Po$. We do not examine such alternative sampling regimes further in this work.

\section{The Existence of Unbiased Tests}
\label{sec:ExistenceResults}
In this section we provide necessary and sufficient conditions for the existence of a level-$\alpha$ unbiased test in a multinomial model. 
\begin{definition}[Unbiasedness threshold]
    The NTUB and SUB unbiasedness thresholds are the respective minimum sample sizes $n$ that are needed for an NTUB or SUB test to exist. The NTUB or SUB threshold is infinite if there respectively does not exist a NTUB or SUB test for any sample size. 
\end{definition}
There are three observations that lead to our main existence result:
\begin{enumerate}
    \item The unbiasedness inequality \eqref{eqn:UBDefinition} is preserved by translation and multiplication by positive scalars in that $ \beta_\phi(\bs{\pi}_0) \leq \beta_\phi(\bs{\pi}_A)$ holds if and only if $ a\beta_\phi(\bs{\pi}_0) + b \leq a\beta_\phi(\bs{\pi}_A) + b$ holds for $a > 0$, $b \in \mb{R}$.
    \item A polynomial $f $ of degree less than or equal to $n$ has a unique representation $f = \sum_{i = 0}^n f_{(i)}$ where $f_{(i)}$ is homogeneous of degree $i$. Any $f$ can be homogenized into the homogeneous polynomial $f_{H}$ of degree $n$ defined by $
    f_H(\bs{\pi}) = \sum_{i = 0}^n \big(\sum_{i = 1}^k \pi_i\big)^{n-i} f_{(i)}(\bs{\pi})$. For any point $\bs{\pi}$ on the probability simplex $f_H(\bs{\pi}) = f(\bs{\pi})$ because $\sum_{i = 1}^n \pi_i = 1$. 
    \item If a power function $\beta$ corresponds to a size-$\alpha$ test the translated polynomial $\Tilde{\beta}(\bs{\pi}_0) \coloneqq \beta_\phi(\bs{\pi}_0) - \alpha$ is non-positive on $\Po$ since $\sup_{\bs{\pi} \in \Po} \Tilde{\beta}(\bs{\pi}) = 0$. Note that this supremum is attained when $\Po$ is closed and thus compact.
\end{enumerate}
A consequence of the first observation is that the hypercube constraints described in Lemma \ref{lem:PowerFunctionCharacterization} are not relevant when examining existence questions, as $\beta_\phi$ can be scaled and translated so that it equals any polynomial with degree bounded by $n$. The second observation implies that any polynomial of degree less than or equal to $n$  can always be homogenized into an equivalent power polynomial of degree equal to $n$. The final observation amounts to the main non-positivity and non-negativity constraints on candidate power polynomials. 

\begin{theorem}
    \label{thm:MainExistenceUniqueness}
Assume that the null hypothesis set is closed.  There exists an NTUB test of $\Po$ against $\Pa = \Delta_{k-1}\backslash \Po$ when the sample size is $n$ if and only if there exists a polynomial $\Tilde{\beta}$ of degree less than or equal to $n$ that is non-constant on $\Delta_{k-1}$ and has sub-level sets with
    \begin{align}
    \label{eqn:NTUBExistence}
        \Po \subseteq \{ \bs{\pi} \in \Delta_{k-1}: \Tilde{\beta}(\bs{\pi}) \leq 0\} \; \text{ and } \; \Pa \subseteq \{ \bs{\pi} \in \Delta_{k-1}: \Tilde{\beta}(\bs{\pi}) \geq 0\}.
    \end{align}
    There exists an SUB test if and only if there exists a similarly defined polynomial with sub-level sets of the form
    \begin{align}
    \label{eqn:SUBExistence}
             \Po = \{ \bs{\pi} \in \Delta_{k-1}: \Tilde{\beta}(\bs{\pi}) \leq 0\} \; \text{ and } \; \Pa = \{ \bs{\pi} \in \Delta_{k-1}: \Tilde{\beta}(\bs{\pi}) > 0\}.
    \end{align}
\end{theorem}
Any polynomial satisfying either the constraint \eqref{eqn:NTUBExistence} or \eqref{eqn:SUBExistence} will be referred to as a (NTUB or SUB) \textit{separating polynomial} for the null hypothesis $\Po$. Note that the assumption that $\Po$ is closed is necessary for a SUB test to exist. 

% \textcolor{red}{Provide a formulation in terms of prime cones/preorders maybe.}

A basic, closed, semialgebraic set is a set that has the form $\cap_{i = 1}^m \{\bl{y} \in \mb{R}^k: f_{(i)}(\bl{y}) \leq 0\}$ for a collection of polynomials $f_{(1)},\ldots,f_{(m)}$. Theorem \ref{thm:MainExistenceUniqueness} shows that for a null hypothesis to have a SUB test it must be a basic, closed, semialgebraic set. Furthermore, $\Po$ must equal the intersection of the simplex with the sub-level set of a \textit{single} separating polynomial $\Tilde{\beta}$. The unbiasedness threshold is precisely the minimum degree polynomial where $\Po$ admits a representation as in \eqref{eqn:NTUBExistence} or \eqref{eqn:SUBExistence}.  

The SUB constraint on $\Po$ imposed by Theorem \ref{thm:MainExistenceUniqueness} is stringent as it entails that $\Po$ must be described in terms of polynomials. The constraint for NTUB tests is less elegant, but is closely related to the SUB constraint in that if there exists a NTUB test for $\Po  \subseteq \{ \bs{\pi} \in \Delta_{k-1}: \Tilde{\beta}(\bs{\pi}) \leq 0\}$ then the enlarged null hypothesis set $\Po' = \{ \bs{\pi} \in \Delta_{k-1}: \Tilde{\beta}(\bs{\pi}) \leq 0\}$ has a SUB test. Example \ref{exm:ExponentialParameter} below will show that, while related, there can be situations where a NTUB test exists but a SUB test does not. 

Another simple yet powerful observation is that unbiased tests have the property of similarity on the boundary \citep{LehmannScheffeCompleteness1}, a consequence of polynomials being continuous.

\begin{lemma}[Similarity on the boundary]
\label{lem:SimilarityonBd}
 Assume that a NTUB or SUB test for $\Po$ exists and $\Tilde{\beta}$ is given as in \eqref{eqn:NTUBExistence} or \eqref{eqn:SUBExistence}. Then $\Tilde{\beta}(\bs{\pi}) = 0$ for every point $\bs{\pi}$ in the boundary of $\Po$ relative to $\Delta_{k-1}$. 
\end{lemma}
% \begin{proof}
%     By the definition of $\Tilde{\beta}$, it is known that $\Tilde{\beta}(\bs{\pi}) \leq 0$ for $\bs{\pi} \in \partial \Po$. As $\bs{\pi} \in \partial \Po$ there exists a sequence $\{\bs{\pi}_i \}_{i = 1}^\infty \subseteq \Pa$ with $\bs{\pi}_i \rightarrow \bs{\pi}$. Note that it is important that $\bs{\pi}$ is in the boundary of $\Po$ relative to $\Delta_{k-1}$ for this to hold. By continuity $0 \leq \Tilde{\beta}(\bs{\pi}_i) \rightarrow \Tilde{\beta}(\bs{\pi}) \leq 0$, as needed. 
% \end{proof}

\begin{figure}
    \centering
    \includegraphics[width=0.6\linewidth]{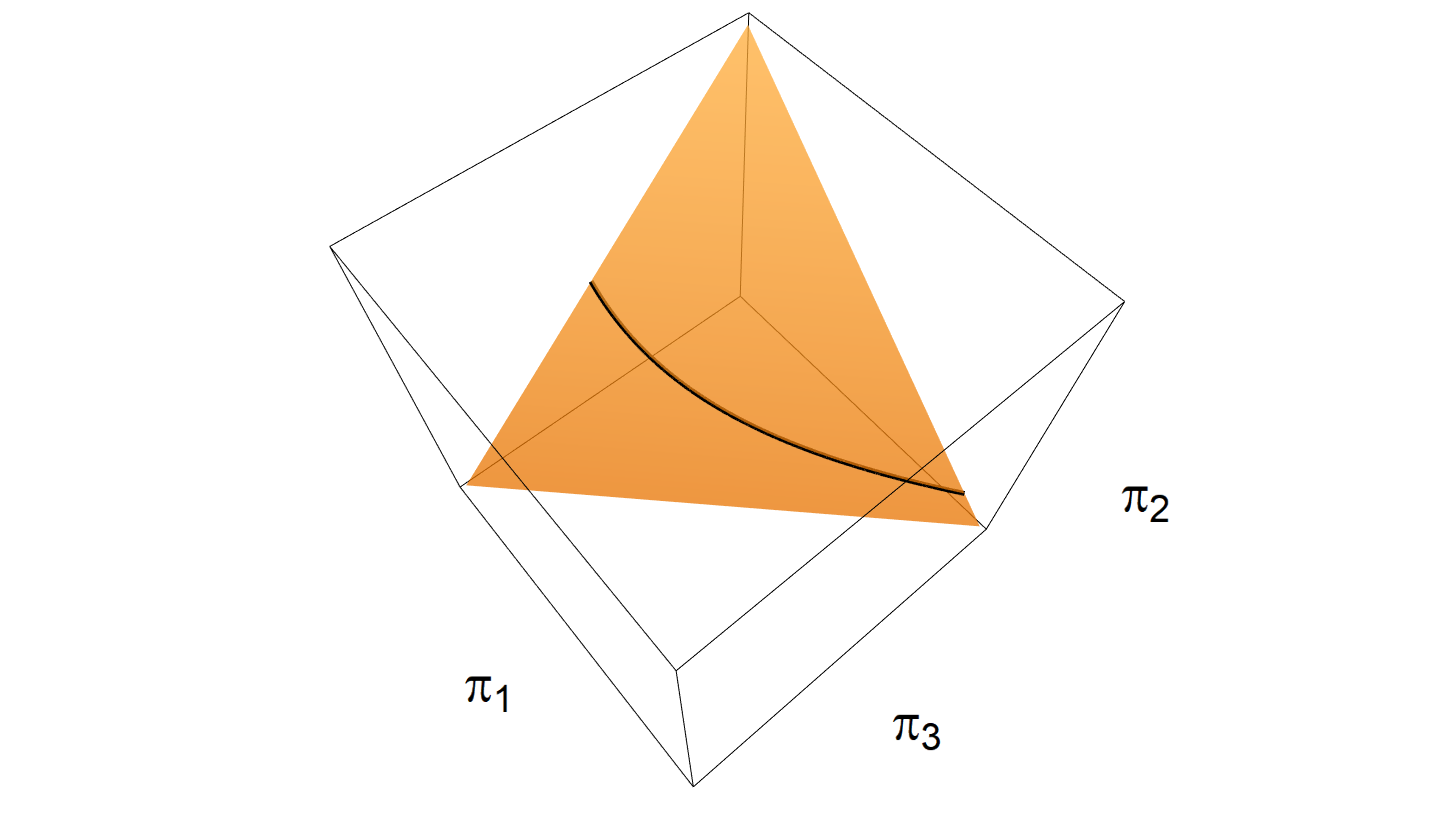}
    \caption{The model $\Po\subset \Delta_2$ (black curve) from Example \ref{exm:ExponentialParameter}.}
    \label{fig:ExponentModel}
\end{figure}

\begin{example}
\label{exm:ExponentialParameter}
    Consider the null hypothesis set $\Po = \{ \bs{\pi} \in \Delta_{2}: \bs{\pi} = \big(\tfrac{1}{2}\theta,\tfrac{1}{2}\exp(-\theta),1 - \tfrac{1}{2}\theta - \tfrac{1}{2}\exp(-\theta)\big), \theta \in \mb{R}\}$ shown in Figure \ref{fig:ExponentModel}. There does not exist a NTUB test of $\Po$ for any sample size. Suppose there existed a $\Tilde{\beta}$ as in \eqref{eqn:NTUBExistence}. We define $f(a,b) = \Tilde{\beta}(a,b,1-a-b)$ and $g(\theta) = f(\tfrac{1}{2}\theta,\tfrac{1}{2}\exp(-\theta))$.  As $\partial \Po = \Po$, similarity on the boundary implies that $g(\theta) = \Tilde{\beta}( \tfrac{1}{2}\theta,\tfrac{1}{2}\exp(-\theta),1 - \tfrac{1}{2}\theta - \tfrac{1}{2}\exp(-\theta)) = 0$ for all $\theta$ in a closed interval $I$. Since $g(\theta)$ is a holomorphic function that vanishes on an interval it must be the zero function \citep[Thm 3.7 IV]{ConwayComplexVar}. Writing $f(a,b) = \sum_{i,j = 0}^n  f_{ij} a^i b^j$, it has been shown that $\sum_{i,j = 0}^n  2^{-(i+j)} f_{ij} \theta^i \exp(-j\theta) = 0$ for all $\theta$. This implies $f_{ij} = 0$ for all $i,j$, and thus the contradiction that $\Tilde{\beta} = 0$, because the functions $\{\theta^i \exp(-j\theta)\}_{0 \leq i,j \leq n}$ are a linearly independent set.
    
    Consider modifying $\Po$ by embedding this set within a facet of $\Delta_3$. Specifically, let $\Po = \{\bs{\pi} \in \Delta_3: \bs{\pi} = (\theta,\exp(-\theta),1-\theta - \exp(-\theta),0)\}$ be contained in the facet of the simplex with $\pi_4 = 0$. The polynomial $\Tilde{\beta}(\bs{\pi}) = \pi_1$ satisfies the condition \eqref{eqn:NTUBExistence} for an NTUB test. However, for any candidate polynomial $\Tilde{\beta}$ the same argument provided in the previous paragraph shows that $\Tilde{\beta}(\pi_1,\pi_2,1-\pi_1-\pi_2,0)$ must be the zero polynomial, and so a SUB test for $\Po$ does not exist.     
\end{example}

The essence of Example \ref{exm:ExponentialParameter} is that  $\Po$ is not semialgebraic due to the non-polynomial, exponential function appearing in its parameterization. The next example demonstrates that even basic, closed, semialgebraic null hypothesis sets may not have NTUB tests if they cannot be described by a single polynomial.

\begin{figure}
     \centering
     \begin{subfigure}[b]{0.45\textwidth}
         \centering       \includegraphics[width=.95\textwidth, height = \textwidth]{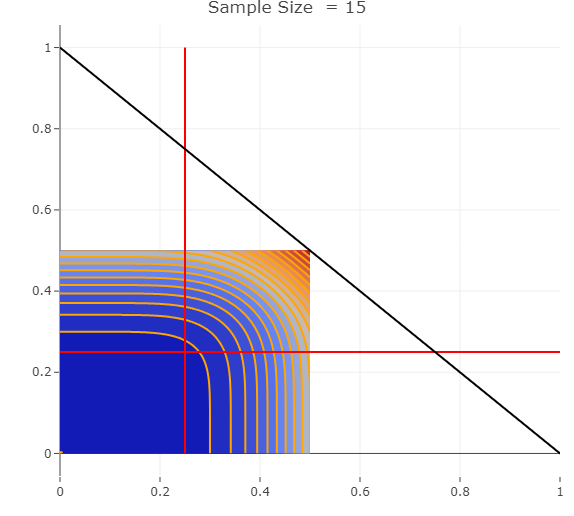}
     \end{subfigure}
     \hfill
     \begin{subfigure}[b]{0.5\textwidth}
         \centering     \includegraphics[width=\textwidth, height= .9\textwidth]{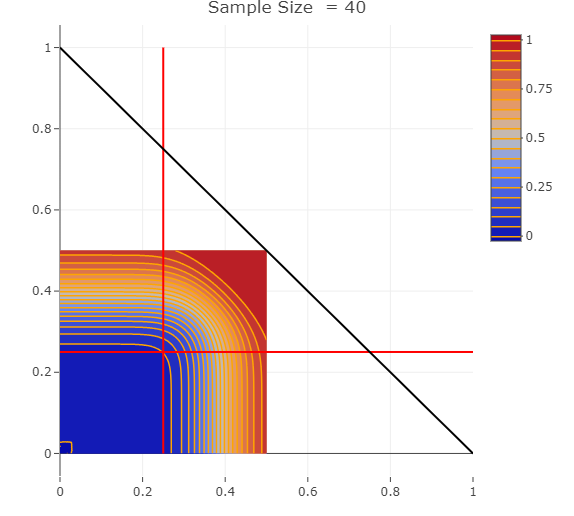}
     \end{subfigure}
     \caption{Contours of the power function $\beta(\bs{\pi})$ of the test statistic $\max(x_1,x_2)$ plotted as a function of $(\pi_1,\pi_2) \in [0,\frac{1}{2}]^2$ for the sample sizes $n = 15$ (left) and $n = 40$ (right) at increments of $0.05$. 
     }
     \label{fig:Square contour plot}
\end{figure}

\begin{example}
\label{ex:Polytope}
Consider the family of square null hypotheses $\mathcal{P}_{0,t} = \{\bs{\pi} \in \Delta_2: (\pi_1,\pi_2) \in [0,t]^2\}$. When $t \in [\frac{1}{2},1]$ there exists a SUB test of $\mathcal{P}_{0,t}$ since the polynomial $\Tilde{\beta}(\bs{\pi}) = -(\pi_1 - t)(\pi_2 - t)$ has sign 
\begin{align*}
    \text{sign}(\Tilde{\beta}(\bs{\pi})) = \begin{cases}
        -1 \;\;\;\; &\text{ if } (\pi_1,\pi_2) \in (t,\infty)^2 \text{ or } (\pi_1,\pi_2) \in (-\infty,t)^2.
        \\
        1 \;\;\;\; &\text{ if } (\pi_1,\pi_2) \in (t,\infty)\times (-\infty,t) \text{ or } (\pi_1,\pi_2) \in (-\infty,t) \times (t,\infty).
        \\
        0 \;\;\;\; &\text{ otherwise }.
    \end{cases}
\end{align*}
When $t < \frac{1}{2}$ this polynomial no longer suffices to provide an unbiased test because the region $(t,\infty)^2$ where $\Tilde{\beta}$ has a negative sign has a non-empty intersection with the alternative hypothesis set $\mathcal{P}_{A,t}$. This is illustrated in Figure \ref{fig:Square contour plot} where the red lines corresponding to $\mc{P}_{0,1/4}$ divide the plane into four regions, and the upper-right region has a non-empty intersection with the simplex $\Delta_2$. A natural test statistic to use for testing $\mc{P}_{0,1/4}$ is $\max(x_1,x_2)$. Figure \ref{fig:Square contour plot} shows the contours of the power function $\text{Pr}_{\bs{\pi}}( \max(x_1,x_2) > \tfrac{n}{4} + \sqrt{n}c)$. The constant $c$ is chosen so that the null hypothesis $\mc{P}_{0,1/4}$, the lower-left square enclosed in the by the red lines, has asymptotic level $0.95$. As the sample size increases it is seen that the test is consistent as the contours have increasing curvature at the point $\bs{\pi} = (\tfrac{1}{4},\tfrac{1}{4},\tfrac{1}{2})$ and more closely approximate the red boundary of the square $[0,\tfrac{1}{4}]^2$. However, for any finite sample size there will always exist points in the alternative hypothesis set that have power less than $0.05$. It is shown in the supplement that at the point $(\tfrac{1}{4},0,\tfrac{3}{4}) \in \mc{P}_{0,3/4}$ on the $\pi_2$-axis the asymptotic power of this test is strictly less than $0.05$.  
\end{example}

\begin{figure}
     \centering
     \begin{subfigure}[b]{0.49\textwidth}
         \centering       \includegraphics[width=.9\textwidth, height = .8\textwidth]{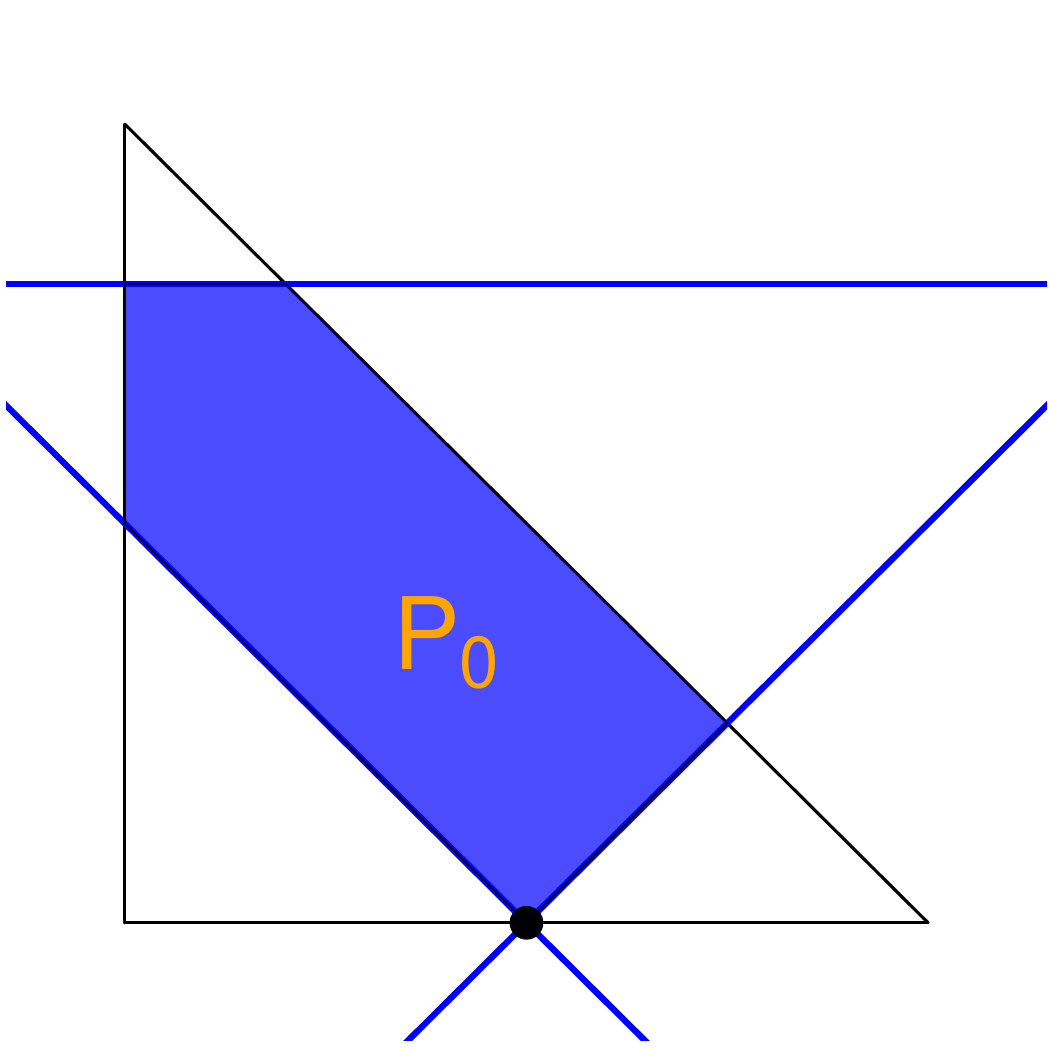}
     \end{subfigure}
     \hfill
     \begin{subfigure}[b]{0.49\textwidth}
         \centering     \includegraphics[width=.9\textwidth, height= .8\textwidth]{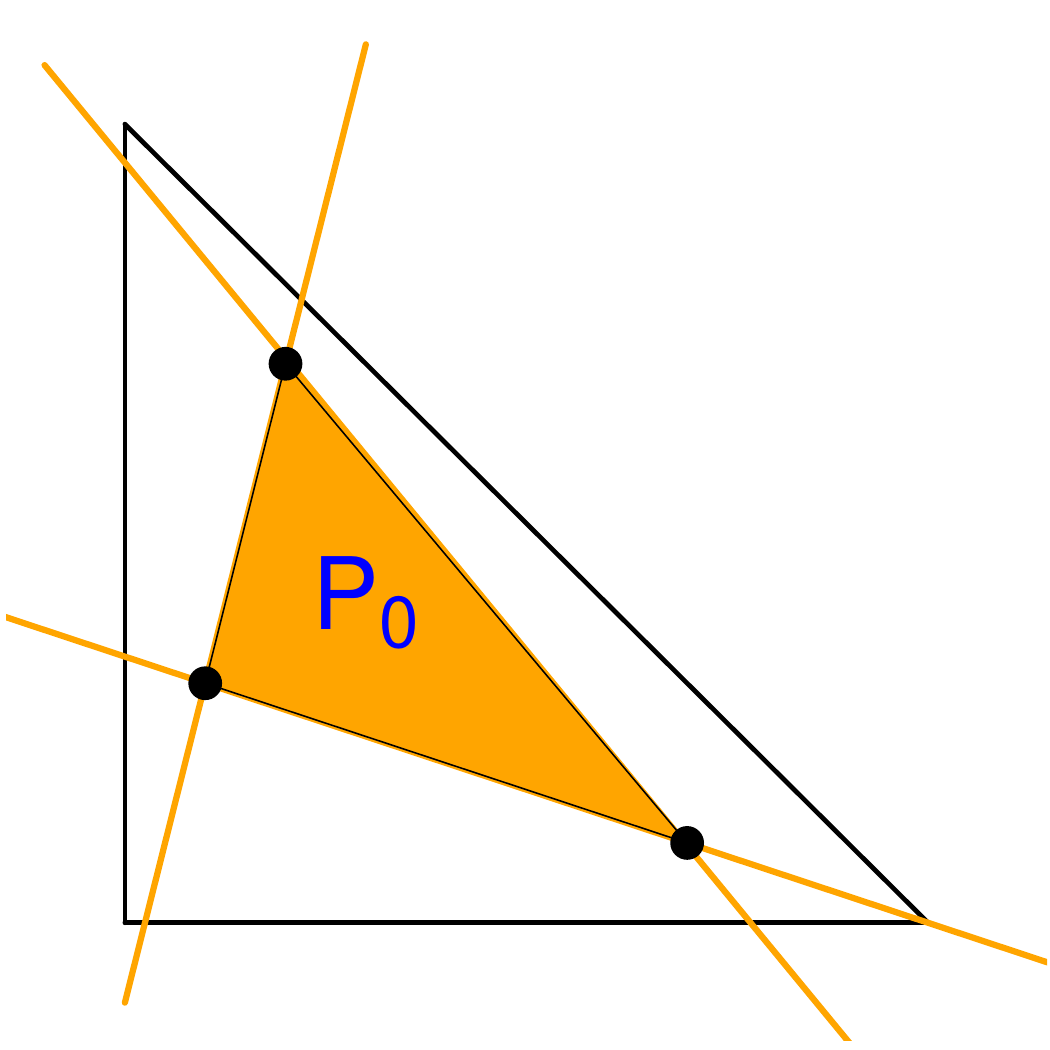}
     \end{subfigure}
     \caption{The polytope null hypothesis on the left has a SUB test, while the hypothesis on the right does not have an NTUB test, due to the three vertices in $\text{int}(\Tilde{\Delta}_2)$.}
     \label{fig:PolytopeNull}
\end{figure}

The reason why there is no unbiased test for $\mc{P}_{0,t}$ with $t < 1/2$ is that the sub-level set of a single polynomial is not able to represent the upper-right corner of the null hypothesis square. The same argument generalizes for other null hypotheses described by affine constraints that have vertices, or lower-dimensional faces, that are contained in the relative interior of the simplex. See the null hypothesis set on the right in Figure \ref{fig:PolytopeNull} for another example. Note that this null hypothesis set depicts a three-component mixture model with the mixture components determined by multinomial distributions located at the three vertices of the polytope. In general, there will not exist unbiased tests for mixture models with known mixture components unless the vertices of the polytope that the mixture defines are contained in the boundary of $\Delta_{k-1}$.

\begin{theorem}[Existence of UB tests for polytope null hypotheses]
\label{thm:PolytopeExistence}
Parameterize the probability simplex $(\pi_{1},\ldots ,\pi_{k}) \in \Delta_{k-1}$ by the projection onto its first $k-1$ coordinates $\Tilde{\Delta}_{k-1} = \{\bs{\pi}: \sum_{i = 1}^{k-1}\pi_i \leq 1, \; \pi_i \geq 0, \; i = 1,\ldots,k-1\} \subseteq \mb R^k$. Let a null hypothesis set be given by the $(k-1)$-dimensional polytope $\mc P_0 = \{\bs{\pi} \in \Tilde{\Delta}_{k-1}: \bl{a}_i^\intercal \bs{\pi} - b_i \geq 0, \; i = 1,\ldots,m\}$, where it is assumed that the facet of the polytope corresponding to the affine hyperplane $H_i = \{\bs{\pi}: \bl{a}_i^\intercal \bs{\pi} = b_i \}$ is $(k-2)$-dimensional and $H_i$ intersects $\mathrm{int}(\Tilde{\Delta}_{k-1})$.
There exists an NTUB or SUB test of $\Po$ if and only if for every $i \neq j$ the set $H_i \cap H_j \cap \Po$ is contained in the boundary of $\Tilde{\Delta}_{k-1}$.     
\end{theorem}
A minor generalization to the above theorem is the case where $\Po$ is not a full-dimensional polytope and is contained in the affine space $\text{aff}(\Po) = \{\bs{\pi} \in \Tilde{\Delta}_{k-1}:\bl{c}_i^\intercal \bs{\pi}  = d_i,\; i = 1,\ldots,r\}$. The polynomial $\Tilde{\beta}(\bs{\pi}) = (\bl{c}_1^\intercal \bs{\pi} - d_1)^2$ yields an NTUB separating polynomial. A SUB separating polynomial will not exist if there are vertices, edges, or more generally, lower-dimensional faces, that intersect the relative interior of $\text{aff}(\Po) \cap \Tilde{\Delta}_{k-1}$. If a face of $\Po$ that is not full-dimensional intersects only the boundary of $\Tilde{\Delta}_{k-1}$, a SUB test exists, such as in the case of the model on left in Figure \ref{fig:PolytopeNull} or $\mc{P}_{0,1/2}$ in  Example \ref{ex:Polytope}.

\begin{example}
The null hypothesis $\Po = \{\bs{\pi} \in \Delta_{k-1}: \pi_1 \leq \pi_2 \leq  \cdots \leq \pi_k\}$ that postulates an ordering for the probabilities of each category does not possess an NTUB test. The point $\bs{\pi} = \tfrac{1}{k} \bl{1}_k$ is a vertex of the polytope $\Po$ that is contained in $\text{int}(\tilde{\Delta}_{k-1})$. 
\end{example}

% where the point $(\tfrac{1}{2},\tfrac{1}{2},0)$ is in  $\{\bs{\pi}: (1,0,0) \bs{\pi} = \pi_1 = \tfrac{1}{2}\} \cap \{\bs{\pi}: (0,1,0) \bs{\pi} = \pi_2 = \tfrac{1}{2}\}$ and $\partial \Tilde{\Delta}_{k-1}$.

As a final example of an important class of distributions that do not possess NTUB tests, we investigate linear hypotheses, not on the multinomial probabilities, but on the natural parameters of the multinomial exponential family \citep{BrownExponentialFamilies}. An important theorem \citep[Thm 4.4.1]{TSPLehmannRomano} concerning the existence of uniformly most powerful unbiased (UMPU) tests in an exponential family that has a density function of the form
\begin{align*}
    p(x|\theta) = \exp\bigg(\sum_{i = 1}^{k-1} \theta_i T_i(x) - \kappa(\theta) \bigg),
\end{align*}
and has natural parameter space $\Theta$, 
states that there exists a UMPU test of the null hypothesis $\Po = \{\bs{\theta} \in \Theta : \theta_1 = c\}$ or $\Po = \{\bs{\theta} \in \Theta : \theta_1 \geq c\}$ against the alternative $\Pa = \Theta \backslash \Po$. After a possible reparameterization, this theorem implies that there exists a UMPU test for linear hypotheses on the natural parameters: $\Po = \{\bs{\theta}: \sum_{i = 1}^{k-1} a_i \theta_i = b\}$. In the multinomial exponential family any linear hypothesis on the natural parameters has the form
\begin{align}
\label{eqn:LinearHypothesisMultinomial}
    \Po = \bigg\{ \bs{\pi} \in  \text{int}(\Delta_{k-1}): \sum_{i = 1}^{k-1} a_i \log\big(\tfrac{\pi_i}{\pi_k}\big) = b \bigg\}. 
\end{align}
Surprisingly, the UMPU test constructed by this theorem can often be vacuous, with the UMPU test being the trivial test $\phi(\bl{x}) = \alpha$.

\begin{theorem}
\label{thm:LogLinearExistence}
When $k \geq 3$ there is an NTUB test of the log-linear null hypothesis \eqref{eqn:LinearHypothesisMultinomial} if and only if there exists a representation of $\Po$ where every $a_i$, $i = 1,\ldots,k-1$ is a rational number.
\end{theorem}

% \begin{theorem}
% \label{thm:LogLinearExistence}
% There is an SUB test for the 
% the Zariski closure of the log-linear null hypothesis
% \begin{align}
% \label{eqn:loglinearnullhypothesis}
%     \Po = \bigg\{ \bs{\pi} \in \text{int}(\Delta_{k-1}): \big[ \log\big(\tfrac{\pi_1}{\pi_k}\big),\ldots, \log\big(\tfrac{\pi_{k-1}}{\pi_k}\big),1\big] \bl{a}_i =   0, \;i = 1,\ldots,m\bigg\},
% \end{align}
%  if and only if there exists a representation of $\Po$ of the form \eqref{eqn:loglinearnullhypothesis} where every $\bl{a}_i$ has rational entries. There exists an NTUB test of $\Po$ if and only if there exists a representation of $\Po$ of the form \eqref{eqn:loglinearnullhypothesis} where at least one of the $\bl{a}_i$ vectors is in $\mb{Q}^k$.  
% \end{theorem}
By the above theorem if $\bl{a}^\intercal = (a_1,\ldots,a_{k-1})$ in \eqref{eqn:LinearHypothesisMultinomial} has one rational and one irrational entry then the UMPU test for linear hypotheses in \cite{TSPLehmannRomano} reduces to the trivial test. For example, there is no NTUB test of the linear hypothesis $\Po = \{\text{int}(\Delta_2) : \sqrt{2}\log(\tfrac{\pi_1}{\pi_3}) = \log(\tfrac{\pi_2}{\pi_3})\}$. The idea behind Theorem \ref{thm:LogLinearExistence} is that if $\bl{a} \in \mb{Q}^{k-1}$ then upon clearing denominators and exponentiating, a linear equation on the log-odds becomes a rational equation in the probabilities. The resulting polynomial is a binomial; further study of such systems of equations in log-linear models can be found in \cite{SturmfelsGeigerMeekToricAlgebraGraphical}. Conversely, when the log-linear equations cannot be represented rationally, upon exponentiating the resulting equations cannot be represented semialgebraically, as some of the exponents of the $\pi_i$s will necessarily be irrational numbers. 

\section{Algebra Background and Algebraic Null Hypotheses}
\label{sec:AlgebraBackground}
A focus of the subsequent sections will be on null hypotheses that have the form
\begin{align}
\label{eqn:AlgebraicNull}
    \Po = \{\bs{\pi} \in \Delta_{k-1}:  f_{(i)}(\bs{\pi}) = 0, \; i = 1,\ldots,m\}.
\end{align}
A collection of real solutions to a system of polynomial equations is referred to as a real algebraic variety, or for brevity, a variety. Specifically, the variety associated with a collection of polynomials $I \subseteq \mb{R}[\bs{\pi}]$ is equal to
\begin{align}
\label{eqn:VarietyDef}
    V_{\mb{R}^k}(I) \coloneqq \{ \bs{\pi} \in \mb{R}^k: f(\bs{\pi}) = 0 \text{ for all } f \in I\}.
\end{align}
A null hypothesis of the form \eqref{eqn:AlgebraicNull} is the intersection of the probability simplex with the variety $V_{\mb {R}^k}(\{f_{(1)},\ldots,f_{(m)}\})$ and will be referred to as an \textit{algebraic null hypothesis}. Strictly unbiased tests exist for all algebraic null hypotheses as \eqref{eqn:AlgebraicNull} has the form of the semialgebraic set \eqref{eqn:SUBExistence} by taking a sum of squares:
\begin{align}
    \Po = \bigg\{  \bs{\pi} \in \Delta_{k-1}: \sum_{i = 1}^m f_{(i)}^2(\bs{\pi}) \leq 0 \bigg\}.
\end{align}
A closely related algebraic counterpart to a variety is an ideal, which given any set of points $S \subset \mb{R}^k$ is the set of all polynomials that vanish on $S$:
\begin{align}
\label{eqn:IdealDef}
    I_{\mb{R}^k}(S) \coloneqq \{ f \in \mb{R}[\pi_1,\ldots,\pi_k]: f(\bs{\pi}) = 0  \text{ for all } \bs{\pi} \in S\}.
\end{align}
 As the complex numbers are algebraically closed, the interaction between ideals and varieties is especially civilized over $\mb{C}^k$ and $\mb{C}[\pi_1,\ldots,\pi_k]$. The functions $V_{\mb{C}^k}(\cdot)$ and $I_{\mb{C}^k}(\cdot)$ denote \eqref{eqn:VarietyDef} and \eqref{eqn:IdealDef} with the enlarged ranges of $\mb{C}^k$ and $\mb{C}[\pi_1,\ldots,\pi_k]$. Complex ideals and varieties will be needed as technical tools that provide insight into the real counterparts.

Any subset of $\Delta_{k-1}$ can be viewed either as a subset of the ambient space $\mb{R}^k$ or as a subset of the $(k-1)$-dimensional affine hyperplane $\text{aff}(\Delta_{k-1})$. The latter set can be parameterized using the coordinates $(\pi_1,\ldots,\pi_{k-1})$ by mapping $\pi_k \mapsto 1 - \sum_{i = 1}^{k-1} \pi_i$. The notation $V_{\mb{R}^{k-1}}(I)$ and $I_{\mb{R}^{k-1}}(S)$ will respectively be used to describe subsets of points in $\text{aff}(\Delta_{k-1}) \cong \mb{R}^{k-1}$ and polynomials in $\mb{R}[\pi_1,\ldots,\pi_{k-1}]$ after the coordinate change $\pi_k \mapsto 1 - \sum_{i = 1}^{k-1} \pi_i$ has been applied. The advantage of making this coordinate change is that the trivial, sum-to-one polynomial constraint $1 - \sum_{i = 1}^k \pi_i$ can be ignored as it is already accounted for by the parameterization. However, polynomial constraints are often most easily formulated in $\mb{R}^k$ rather than $\mb{R}^{k-1}$ as the latter reparameterization is does not treat the probabilities $\pi_1,\ldots,\pi_k$ symmetrically. 
% The notation $\Delta_{k-1}$ will be used to denote the probability simplex in $\mb{R}^k$. 
The projection of the probability simplex $\Delta_{k-1}$ into $\mb{R}^{k-1}$ will be denoted by $\Tilde{\Delta}_{k-1}$ and is defined by the expression provided in Theorem \ref{thm:PolytopeExistence}.

Like vector subspaces, ideals are closed under certain algebraic operations; if $f,g \in  I_{\mb{R}^k}(S)$ then $f + g \in  I_{\mb{R}^k}(S)$, and if $f \in  I_{\mb{R}^k}(S)$ then $hf \in  I_{\mb{R}^k}(S)$ for any $h \in \mb{R}[\pi_1,\ldots,\pi_k]$. The notion of the span of a collection of vectors generalizes to ideals where the ideal generated by a collection of polynomials is defined as
\begin{align}
\label{eqn:IdealGeneration}
    \langle f_{(1)},\ldots,f_{(m)} \rangle \coloneqq \bigg\{ \sum_{i = 1}^m h_{(i)}f_{(i)} : h_{(i)} \in \mb{R}[\pi_1,\ldots,\pi_k]\bigg\}.
\end{align}
An important notion that is an ideal-theoretic analogue of a vector space basis is that of a reduced Gr\"obner basis $\{g_{(1)},\ldots,g_{(m)}\} \subset I$. A Gr\"obner basis generates $I$ as in \eqref{eqn:IdealGeneration} and also has additional properties that ensure that if $f \in I$ then polynomial long division of $f$ by $\{g_{(1)},\ldots,g_{(m)}\}$ will result in a remainder of zero. A reduced Gr\"obner basis of an ideal is uniquely determined by a monomial order: a total order on the monomials of $\mb{R}[\bs{\pi}]$ that has the property that $\bs{\pi}^{\bs{\alpha}} \geq \bs{\pi}^{\bs{\beta}}$ implies $\bs{\pi}^{\bs{\alpha} + \bs{\gamma}} \geq \bs{\pi}^{\bs{\beta} + \bs{\gamma}}$ for $\bs{\alpha}, \bs{\beta}, \bs{\gamma} \in \mb{N}^k$. Note that in this notation $\bs{\pi}^{\bs{\alpha}}$ is defined as $\prod_{i = 1}^k \pi_i^{\alpha_i}$. Of interest will be graded monomial orderings that have the property that if $\text{deg}(\bs{\pi}^{\bs{\alpha}}) > \text{deg}(\bs{\pi}^{\bs{\alpha}})$ then $\bs{\pi}^{\bs{\alpha}} > \bs{\pi}^{\bs{\beta}}$ in the ordering. One example of a graded monomial order is the graded lexicographic order where $\bs{\pi}^{\bs{\alpha}} > \bs{\pi}^{\bs{\beta}}$ if and only if $\text{deg}(\bs{\pi}^{\bs{\alpha}}) > \text{deg}(\bs{\pi}^{\bs{\alpha}})$ or $\text{deg}(\bs{\pi}^{\bs{\alpha}}) = \text{deg}(\bs{\pi}^{\bs{\alpha}})$ and there exists a $j$ with $\alpha_i = \beta_i$, $i < j$ and $\alpha_j > \beta_j$. The reader is referred to \cite{cox2013ideals} for a more in-depth introduction to ideals, varieties, and Gr\"obner bases.

Similarity on the boundary of an unbiased test can be rephrased as an ideal constraint on the separating polynomial. For algebraic null hypotheses every point in $\Po$ is a boundary point. Algebraic methods for characterizing ideals will be utilized in the subsequent section to search for low-degree polynomials within the ideal. 
\begin{corollary}
\label{cor:ExistenceforAlgebraic}
        There exists an unbiased test for the algebraic null hypothesis \eqref{eqn:AlgebraicNull} when the sample size is $n$ if and only if there exists a degree $n$, non-zero polynomial $\Tilde{\beta} \in I_{\mb{R}^k}(\Po$) where 
    \begin{align}
    \label{eqn:AlgNullNTUB}
        \Pa \subseteq \{ \bs{\pi} \in \Delta_{k-1}: \Tilde{\beta}(\bs{\pi}) \geq 0\}
    \end{align}
    if the test is to be a NTUB test, and
    \begin{align}
    \label{eqn:AlgNullSUB}
                \Pa = \{ \bs{\pi} \in \Delta_{k-1}: \Tilde{\beta}(\bs{\pi}) >  0\}
    \end{align}
    if the test is to be a SUB test. 
\end{corollary}
% \begin{proof}
%     This result follows immediately from Lemma \ref{lem:SimilarityonBd} if it is shown that  $\partial \Po = \Po$. Let $\Po = \{\bs{\pi} \in \Delta_{k-1}: f(\bs{\pi}) = 0\} \subsetneq \Delta_{k-1}$. Suppose by contradiction that $\bs{\pi}_0 \in \Po \backslash \partial \Po$, so that there exists an open ball $U$ with $\bs{\pi}_0 \in U \subseteq \Po$. Then $f\vert_U = 0$ and in particular there exists an open ball $V \subset \mb{R}^{k-1}$ with $f(\pi_1,\pi_2,\ldots,1 - \sum_{i = 1}^{k-1}\pi_i) = 0$ for all $(\pi_1,\ldots,\pi_{k-1}) \in V$. Any polynomial that vanishes on an open ball must be the zero polynomial so $f(\pi_1,\pi_2,\ldots,1 - \sum_{i = 1}^{k-1}\pi_i) = 0$ over $\mb{R}^{k-1}$, implying the contradiction that $\Po = \Delta_{k-1}$.
% \end{proof}

\begin{example}
\label{ex:2x2Indep}
A motivating example of an algebraic null hypothesis is the following independence variety defined on a $2 \times 2$ contingency table 
\begin{align}
    \label{eqn:2x2IndepModel}
    \Po = \bigg\{ \bs{\pi} = \begin{bmatrix}
        \pi_{11} & \pi_{12}
        \\
        \pi_{21} & \pi_{22}
    \end{bmatrix} \in \Delta_3 : \det(\bs{\pi}) = \pi_{11}\pi_{22} - \pi_{12}\pi_{21} = 0 \bigg\}
\end{align}
that postulates that the marginal row and column counts are independent.
% Corollary \ref{cor:ExistenceforAlgebraic} applies here since each element $\bs{\pi} \in \Po$ is a limit of points in $\Pa$ as a small perturbation of the entries in $\bs{\pi}$ can produce a full-rank matrix. 
% It is apparent from the illustration of this set in Figure \textcolor{red}{xxx} that every point in $\Po$ is a limit of points in $\Pa$ as is stated in Corollary \ref{cor:ExistenceforAlgebraic}.
By the definition of $\Po$,  $\det(\bs{\pi}) \in I_{\mb{R}^{k-1}}(\Po)$ where $\det(\bs{\pi})$ is viewed as a polynomial in $\mb{R}[\pi_{11},\pi_{12},\pi_{21}]$ after a change of coordinates. Lemma \ref{lem:ParameterizedIdealRadical} in the next section can be used to show that $I_{\mb{R}^{k-1}}(\Po) = \langle \det(\bs{\pi}) \rangle = \langle \pi_{11}(1 - \pi_{11} - \pi_{12} - \pi_{21}) - \pi_{12}\pi_{21} \rangle$.  
\end{example}

% \textcolor{red}{Picture of indep variety --- relate to corollary 1. The ideal of the indep variety related to the corr. }

\section{Computing the Unbiasedness Threshold}
\label{sec:UBThreshold}

Due to the results in Sections \ref{sec:ExistenceResults} and \ref{sec:AlgebraBackground},  unless otherwise stated, it will be assumed in following sections that the null hypotheses considered are algebraic.

\subsection{A Sum of Squares Upper Bound on the Unbiasedness Threshold}
Finding an NTUB or SUB test for an algebraic null hypothesis amounts to finding a non-negative or positive polynomial in $I_{\mb{R}^k}(\Po)$ by Corollary \ref{cor:ExistenceforAlgebraic}. A straightforward way to construct non-negative polynomials is to form a sum of squares of polynomials in $I_{\mb{R}^k}(\Po)$.

\begin{definition}[Sum of squares]
    A polynomial $\Tilde{\beta}$ has a sum of squares (SOS) representation if it can be written in the form $\Tilde{\beta} = \sum_{i = 1}^m f_{(i)}^2$.
\end{definition}

If it is assumed that $\Po = V_{\mb{R}^k}(f_{(1)},\ldots,f_{(m)}) \cap \Delta_{k-1}$, with all $f_{(i)}$ non-constant, then $\Tilde{\beta} = \sum_{i = 1}^m f_{(i)}^2$  satisfies \eqref{eqn:AlgNullSUB} because $\Tilde{\beta}(\bs{\pi}_0) = 0$ if and only if $\bs{\pi}_0 \in V_{\mb{R}^k}(f_{(1)},\ldots,f_{(m)})$. The degree of  $\sum_{i = 1}^m f_{(i)}^2$ is twice the maximum of the degrees of $f_{(1)},\ldots,f_{(m)}$. In a similar way, $\Tilde{\beta} = f_{(i)}^2$ with degree $2 \,\text{deg}(f_{(i)})$ satisfies \eqref{eqn:AlgNullNTUB}.
% Recall that any $\Tilde{\beta}$ satisfying either \eqref{eqn:AlgNullNTUB} or \eqref{eqn:AlgNullSUB} can be normalized to give an NTUB or SUB test when the sample size is the degree of $\Tilde{\beta}$.
Consequently, $2 \max_{i}\text{deg}(f_{(i)})$ and $2 \min_{i}\text{deg}(f_{(i)})$ provide upper bounds on the unbiasedness thresholds of SUB and NTUB tests. These upper bounds may not always equal the unbiasedness threshold as the extra constraint that $\Tilde{\beta}$ is SOS has been imposed. The goal of this subsection is to find the smallest degree $\Tilde{\beta}$ under the assumption that $\Tilde{\beta}$ is SOS. In the next subsection it will be shown that for a large collection of models these upper bounds equal the UB thresholds.  

It is possible that a null hypothesis $\Po$ can be represented as an algebraic variety in different ways. For instance, in Example \ref{ex:2x2Indep} the null hypothesis is equal to $ V_{\mb{R}^k}( \det(\bs{\pi})) \cap \Delta_3$ by definition, but $\Po$ can also be represented as $V_{\mb{R}^{k}}(\det(\bs{\pi})^t) \cap \Delta_3$ for $1 < t \in \mb{N}$. In finding an upper bound for the UB threshold by taking a sum of squares it is desired to find the most degree-economical representation of $\Po$ of the form $V_{\mb{R}^{k}}(f_{(1)},\ldots,f_{(m)}) \cap \Delta_{k-1}$. The following lemma shows that this can be done by computing a Gr\"obner basis with respect to a graded monomial order.

\begin{theorem}
\label{thm:UBThresholdSOSUpperBounds}
    Assume that $\Po$ is an algebraic null hypothesis and let $I_{\mb{R}^{k-1}}(\Po)$ have the reduced Gr\"obner basis $\mc{G} = \{g_{(1)},\ldots, g_{(m)}\}$ with respect to a graded monomial order.  The minimum degree of a non-zero polynomial in $I_{\mb{R}^{k-1}}(\Po)$ is equal to $\min_{i}\mathrm{deg}(g_{(i)})$. The lowest degree of a SOS polynomial $\Tilde{\beta}$ corresponding to an NTUB test is $2\min_{i}\mathrm{deg}(g_{(i)})$.  Let $\mc{G}_i$ consist of the subset of polynomials in $\mc{G}$ with degree less than or equal to $i$. Define $d = \min \{i:  V_{\mb{R}^{k-1}}(\mc{G}_i) \cap \Delta_{k-1} = \Po \}$ to be the smallest $i$ such that the polynomials in $\mc{G}_i$ cut out exactly the null hypothesis set $\Po$. The lowest degree SOS polynomial corresponding to a SUB test is $2d$.
\end{theorem}

The application of the above theorem requires that a Gr\"obner basis for $I_{\mb{R}^{k-1}}(\Po)$ be computed. If a generating set for $I_{\mb{R}^{k-1}}(\Po)$ is known, a Gr\"obner basis can be computed by Buchberger's algorithm \cite[Sec 2.7]{cox2013ideals}. Finding a set of generating polynomials for $I_{\mb{R}^{k-1}}(\Po)$ is not always straightforward as even if $\Po = V_{\mb{R}^{k-1}}(f_{(1)},\ldots,f_{(m)}) \cap \Delta_{k-1}$ the set of polynomials that vanish on $\Po$, namely $I_{\mb{R}^{k-1}}(\Po)$, can be larger than $\langle f_1, \ldots, f_m \rangle$. However, even if generators of $I_{\mb{R}^{k-1}}(\Po)$ cannot be found, it is still possible to obtain an upper bound on the unbiasedness thresholds if a set of polynomials $\{f_{(1)},\ldots,f_{(m)}\}$ with $\Po = V_{\mb{R}^{k-1}}(f_{(1)},\ldots,f_{(m)}) \cap \Delta_{k-1}$ is known. As an algebraic statistical model is by definition the zero set of a collection of polynomials, the SOS approach always provides an upper bound on the unbiasedness thresholds. 

\begin{corollary}
    Assume that $\Po = V_{\mb{R}^{k-1}}(f_{(1)},\ldots,f_{(m)}) \cap \Delta_{k-1}$ with every $f_{(i)}$ non-zero. The NTUB unbiasedness threshold is bounded above by $2 \min_{i} \mathrm{deg}(f_{(i)})$. The SUB unbiasedness threshold is bounded above by $2\max_{i} \mathrm{deg}(f_{(i)})$. 
\end{corollary}
% \begin{proof}
%     The proof of this Corollary follows from the proof of Theorem \textcolor{red}{REF}. In particular, it only needs to be shown that $\Tilde{\beta} = f_{(i)}^2$ where $f_{(i)}$ is the polynomial with the smallest degree in $\{f_{(1)},\ldots,f_{(m)}\}$ and $\Tilde{\beta} = \sum_{i = 1}^m f_{(i)}^2$ respectively are polynomials  associated with NTUB and SUB tests. 
% \end{proof}

It is usually the case that an algebraic statistical model is described by means of a polynomial parameterization $\varphi: S \rightarrow \Po$ onto the null hypothesis set for a subset $S \subset \mb{R}^p$. A polynomial parameterization is a function $\varphi = (\varphi_1,\ldots,\varphi_k)$ where each $\varphi_i$ is a multivariate polynomial. For instance, the independence variety introduced in Example \ref{ex:2x2Indep} is parameterized by the function $\varphi: [0,1]^2 \rightarrow \Po$ with 
\begin{align*}
    \varphi(a,b) = \begin{bmatrix}
        a
        \\
        1-a
    \end{bmatrix} \begin{bmatrix}
        b & 1-b
    \end{bmatrix} = \begin{bmatrix}
        ab & a(1-b)
        \\
        (1-a)b & (1-a)(1-b)
    \end{bmatrix}.
\end{align*}
This parameterization clarifies the independence model terminology, as a single observation $\bl{X} \in \mb{N}^{2 \times 2}$ with $n = 1$ from a distribution in $\Po$ is obtained by independently drawing $i \sim \text{Bern}(a)$ and $j \sim \text{Bern}(b)$ and recording a count of one in the $(i+1)$st row and $(j+1)$st column of $\bl{X}$. Under mild conditions it is possible to compute a set of generators of $I_{\mb{R}^{k-1}}(\Po)$ for parameterized, algebraic models.

% The idea of a polynomial parameterization can be extended to a rational parameterization. 
% \textcolor{red}{Deal with these rational parameterizations}
% A rational mapping  $\varphi: \mb{C}^p \rightarrow V$ of a variety $V \subseteq \mb{C}^k$ is a function $\varphi = (\varphi_1,\ldots,\varphi_k)$ where $\varphi(\bs{x}) \in V$ for all $\bs{x} \in \mb{C}^k$ and there exists at least one $\bs{x}$ where the denominators of all $\varphi_i$ do not vanish \citep[Def 5.5.4]{cox2013ideals}. 
% where $\varphi = (\varphi_1,\ldots,\varphi_k)$ and each $\varphi_i$ is a rational function. Rational parameterizations      

\begin{lemma}
\label{lem:ParameterizedIdealRadical}
Let $\Po = V_{\mb{R}^{k-1}}(f_{(1)},\ldots,f_{(m)}) \cap \Delta_{k-1}$ and 
    let $\varphi:\mb{C}^p \rightarrow V_{\mb{C}^{k-1}}(f_{(1)},\ldots,f_{(m)})$ be a polynomial parameterization onto a (complex) variety. If $\Po = \varphi(S)$ for some set $S \subset \mb{R}^p$ with a non-empty interior then $I_{\mb{R}^{k-1}}(\Po) =  \sqrt[\mb{C}]{\langle f_{(1)},\ldots,f_{(m)} \rangle} \cap \mb{R}[\bs{\pi}]$.
\end{lemma}
% \begin{proof}
%     First assume that $f \in I_{\mb{R}}(\Po)$. It follows that $f(\varphi(\bs{x})) = 0$ for all $\bs{x} \in S$. As $S$ has a non-empty interior $f \circ \varphi$ is the zero polynomial. It follows that $f \in I_{\mb{C}}( V_{\mb{C}}(f_{(1)},\ldots,f_{(m)}))$ because the image of $\varphi$ is $V_{\mb{C}}(f_{(1)},\ldots,f_{(m)})$ by assumption.  The nullstellensatz implies $ f \in \sqrt[\mb{C}]{\langle f_{(1)},\ldots,f_{(m)} \rangle} \cap \mb{R}[\bs{\pi}]$. Conversely, if  $ f \in \sqrt[\mb{C}]{\langle f_{(1)},\ldots,f_{(m)} \rangle} \cap \mb{R}[\bs{\pi}]$
%  then by the nullstellensatz $f \in I_{\mb{C}}( V_{\mb{C}}(f_{(1)},\ldots,f_{(m)}))$, implying that $f$ vanishes on $\Po \subset V_{\mb{C}}(f_{(1)},\ldots,f_{(m)})$. 
%  \end{proof}
Lemma \ref{custlem:RationalParamIdeal} in the supplement extends Lemma \ref{lem:ParameterizedIdealRadical} to null hypotheses that have a parameterization in terms of rational, but not necessarily polynomial functions. Featured in Lemma \ref{lem:ParameterizedIdealRadical} is the radical \citep[Sec 4.2]{cox2013ideals} of the ideal $\langle f_{(1)},\ldots, f_{(m)}\rangle$  defined as $\sqrt[\mb{C}]{\langle f_{(1)},\ldots,f_{(m)} \rangle} \coloneqq \{f: \exists k \; f^k \in \langle f_{(1)},\ldots,f_{(m)} \rangle\}$. The main importance of Lemma \ref{lem:ParameterizedIdealRadical} is that there exist procedures for computing a Gr\"obner basis of such radical ideals when $f_{(1)},\ldots,f_{(m)}$ are provided. In the independence model example
\begin{align*}
    \sqrt[\mb{C}]{ \langle \pi_{11}(1 - \pi_{11} - \pi_{12} - \pi_{21}) - \pi_{12}\pi_{21} \rangle} = \langle \pi_{11}(1 - \pi_{11} - \pi_{12} - \pi_{21})- \pi_{12}\pi_{21} \rangle
\end{align*}
% $\sqrt[\mb{C}]{ \langle \pi_{11}\pi_{22} - \pi_{12}\pi_{21} \rangle} = \langle \pi_{11}\pi_{22} - \pi_{12}\pi_{21} \rangle$
has $\{\pi_{11}(1 - \pi_{11} - \pi_{12} - \pi_{21}) - \pi_{12}\pi_{21}  \}$ as a Gr\"obner basis under any graded, monomial order. Theorem \ref{thm:UBThresholdSOSUpperBounds} implies that $n = 4$ is the optimal sum of squares upper bound for both the NTUB and SUB unbiasedness thresholds. In the next section it will be shown that $n = 4$ is equal to the unbiasedness threshold.

\begin{example}
\label{ex:IndepModelpq}
The independence model can be extended to $p \times q$ contingency tables by the parameterization $\varphi(\bs
{a},\bs{b}) = \bs{a}\bs{b}^\intercal$ with $\bs{a} \in \Delta_{p-1}$, $\bs{b} \in \Delta_{q-1}$. This is equivalent to requiring that the matrix of probabilities $\bs{\pi} \in \mb{R}^{p \times q}$ be a rank one matrix. Thus
\begin{align*}
    \Po = \{\bs{\pi} \in  \Delta_{pq - 1}: \det(\bs{\pi}_{I,J}) = \pi_{i_1j_1}\pi_{i_2j_2} - \pi_{i_1j_2}\pi_{i_2j_1} = 0,\;   \forall \;  (i_1,i_2) \in [p] \times [p], \; (j_1,j_2) \in [q] \times [q] \}
\end{align*}
expresses the independence model as an algebraic model. In this general case there is a polynomial constraint for every $2 \times 2$ minor in the matrix $\bs{\pi}$. Other constraints can be incorporated into the independence model; for instance it might be conjectured that the first row and the first column of the table have the marginal probabilities that satisfy: $2\sum_{i = 1}^p \pi_{i1} = \sum_{j = 1}^q \pi_{1j}$. Expressed in terms of the parameterization this adds the constraint that $2b_1 = a_1$. We examine the independence model for a $2 \times 3$ table with the added constraint that $2b_1 = a_1$. In particular, $\Po = V_{\mb{R}^k}(I) \cap \Delta_5$ where $I = \langle \pi_{11}\pi_{22} - \pi_{12}\pi_{21},  \pi_{11}\pi_{23} - \pi_{21}\pi_{13},  \pi_{12}\pi_{23} - \pi_{13}\pi_{22}, \pi_{11} + 2\pi_{21} - \pi_{12} - \pi_{13}\rangle$. The software \texttt{Macaulay 2} \citep{Macaulay2}, along with the expression in Lemma \ref{lem:ParameterizedIdealRadical}, can be used to compute the following reduced Gr\"obner basis of 
$I_{\mb{R}^{k-1}}(\Po)$ with respect to the graded reverse lexicographic ordering:
\begin{align*}
   \mc{G} = \big\{& \pi_{11}-\pi_{12}-\pi_{13}+2\,\pi_{21}, \;\;\;  \pi_{12}\pi_{21}-\pi_{12}\pi_{22}-\pi_{13}\pi_{22}+2\pi_{21}\pi_{22},
      \\
&2\pi_{12}^{2}+4\pi_{12}\pi_{13}+2\pi_{13}^{2}-4\pi_{13}\pi_{21}+2\pi_{21}^{2}-4\pi_{12}\pi_{22}-4\pi_{13}\pi_{22}+8\pi_{21}\pi_{22}-\pi_{12}-\pi_{13}+2\pi_{21},
\\
&2\pi_{13}^{2}\pi_{21}-4\pi_{13}\pi_{21}^{2}+2\,\pi_{21}^{3}+2\pi_{12}\pi_{13}\pi_{22}+2\,\pi_{13}^{2}\pi_{22}-8\pi_{13}\pi_{21}\pi_{22}+ \cdots
\\
&6\pi_{21}^{2}\pi_{22}-4\pi_{12}\pi_{22}^{2}-4\pi_{13}\pi_{22}^{2}+8\pi_{21}\pi_{22}^{2}-\pi_{13}\pi_{21}+2\pi_{21}^{2} \big\}.
\end{align*}
The Gr\"obner basis has one polynomial with degree one, two polynomials with degree two, and one polynomial with degree three. The degree three polynomial is redundant in that the degree two and one polynomials are enough to cut out the null hypothesis: $V_{\mb{R}}(\mc{G}_2) \cap \Delta_{5} = \Po$. However, the linear polynomial alone is not enough to cut out null hypothesis as $\Po \subsetneq V_{\mb{R}}(\mc{G}_1) \cap \Delta_5$. We conclude from Theorem \ref{thm:UBThresholdSOSUpperBounds} and Lemma \ref{lem:ParameterizedIdealRadical} that the SOS upper bound of the NTUB and SUB unbiasedness thresholds are $2$ and $4$ respectively, which are twice the minimum and maximum degrees of polynomials in $\mc{G}_2$. The SOS function $\Tilde{\beta}$ associated with the non-trivial unbiased test when $n = 2$ is the square of a linear polynomial. The associated non-trivial unbiased test is effectively only testing if the constraint  $2\pi_{11} + 2\pi_{21}=  \pi_{11} + \pi_{12} + \pi_{13}$ holds, while ignoring the other independence constraints. Although it requires a slightly larger sample size, the strictly unbiased test has non-trivial power against any alternative hypothesis.  
\end{example}

\begin{example}
\label{ex:MotzkinHypothesis}
There exist positive polynomials that cannot be represented as a sum of squares. One well-known example of such a polynomial is the Motzkin polynomial $f(x,y,z) = z^6 + x^2 y^4 + y^2 x^4 - 3x^2 y^2 z^2 \geq 0$ \citep{MotzkinPolynomial}. Define the Motzkin null hypothesis to be
        \begin{align*}
            \Po = \{(\pi_1,\pi_2,\pi_3,\pi_4) \in \Delta_3: f(\pi_1,\pi_2,\pi_3) = 0\}.
        \end{align*}
% As shown in Figure \textcolor{red}{FIGURE REF} 
The null hypothesis set is non-empty as it can be checked that $(\tfrac{1}{4},\tfrac{1}{4},\tfrac{1}{4},\tfrac{1}{4}) \in \Po$. The polynomial $\Tilde{\beta} = f$ yields a SUB test by Corollary \ref{cor:ExistenceforAlgebraic}. However, $f$ cannot be represented as a sum of squares, implying that the SOS approach cannot always construct the separating polynomial of every unbiased test.
\end{example}

\subsection{Unbiasedness Threshold Examples}
In the previous section, upper bounds for the unbiasedness thresholds were found by constructing SOS power polynomials. In this section the unbiasedness thresholds for algebraic models are explicitly computed for a collection of examples. We first discuss a simpler setting where the ideal generated by $\Po$ is principal and then turn to the general case. 

\subsubsection{Principal Vanishing Ideals}
\label{sec:PrincipalVanishing}
An ideal $I = \langle f \rangle$ is said to be principal when it is generated by a single element. Sums of squares formed from the principal ideal $I$ are especially simple as $f^2$ will divide any sum of squares. As such, the sum of squares polynomial in $I$ with the smallest degree will be a positive scalar multiple of $f^2$. The next theorem uses this idea to compute the unbiasedness threshold and characterizes every unbiased test for any model with a principal ideal that satisfies some mild regularity conditions.

\begin{theorem}
\label{thm:PrincipalIdealUBThresholdUMP}
   Suppose that $\Po$ is algebraic, $\overline{\mathrm{int}(\tilde{\Delta}_{k-1}) \cap \Po} = \Po$, and $I_{\mb{R}^{k-1}}(\Po) = \langle f \rangle$ is a principal ideal generated by a single polynomial. If $f$ has a non-zero gradient throughout $\Po$, the NTUB and SUB unbiasedness thresholds are both equal to $2 \,\mathrm{deg}(f)$. For every size $0<\alpha<1$ when the sample size is $n = 2 \,\mathrm{deg}(f)$ there is a unique, uniformly most powerful NTUB test with size $\alpha$ that is also the unique, uniformly most powerful SUB test. This test has a power function of the form $\beta = c_{\alpha} f^2 + \alpha$, where $c_{\alpha} \in (0,\infty)$ is the largest constant that makes $\beta$ a power polynomial. Every size-$\alpha$ NTUB test of $\Po$ has a power function of the form $cf^2 + \alpha$, $c \in (0,c_\alpha]$ when $n = 2 \,\mathrm{deg}(f)$. More generally, if $\beta$ is the power function of a level-$\alpha$, NTUB test for a sample size of $n$ there exists a polynomial $h$ that is non-negative on $\Delta_{k-1}$ with degree $n - 2\mathrm{deg}(f)$ such that $\beta = f^2h + \alpha$. For a SUB test $h$ must in addition be strictly positive on $\Pa$. 
\end{theorem}

As discussed in Example \ref{ex:MotzkinHypothesis}, the UB threshold of the Motzkin hypothesis is at most six, and the UB threshold for the hypothesis $\Po = \{\bs{\pi} \in \Delta_{k-1}: \pi_1 = 0\}$ is one since $\beta(\bs{\pi}) = \pi_1$ is a separating polynomial. This does not contradict the above theorem since both of these hypotheses fall outside of its purview: the Motzkin hypothesis because $\nabla f(\bs{\pi}_0) = 0$ on $\Po$, and the latter hypothesis because $\Po$ is entirely contained in the boundary of $\tilde{\Delta}_{k-1}$.

Theorem \ref{thm:PrincipalIdealUBThresholdUMP} can be applied to construct the UMPU tests for log-linear null hypotheses \eqref{eqn:LinearHypothesisMultinomial} as long as $\Po$ is algebraic. One limitation of Theorem \ref{thm:PrincipalIdealUBThresholdUMP} as compared to the classical test \cite[Thm 4.4.1]{TSPLehmannRomano} for linear hypotheses in the natural parameters of an exponential family is that Theorem \ref{thm:PrincipalIdealUBThresholdUMP} does not indicate how to construct UMPU tests, or even whether UMPU tests exist, for sample sizes that are larger than the UB threshold. However, Theorem \ref{thm:PrincipalIdealUBThresholdUMP} applies more generally beyond exponential subfamilies that are linear in the natural parameters to a large class of curved exponential subfamilies \cite{BrownExponentialFamilies}. Less is known about general conditions for the existence of UMPU tests when $\Po$ is a curved exponential family submodel. One relevant contribution is the paper \citep{SenGuptaExponentialUMPU} that exhibits UMPU tests for curved exponential families of multivariate exponential distributions. In next example a test for a hypothesis that is not linear in the log-odds is explicitly constructed.

% An important theorem \citep[Thm 4.4.1]{TSPLehmannRomano} concerning the existence of uniformly most powerful unbiased (UMPU) tests in an exponential family that has a density function of the form
% \begin{align*}
%     p(x|\theta) = \exp\bigg(\sum_{i = 1}^k \theta_i T_i(x) - \kappa(\theta) \bigg),
% \end{align*} 
% is that there exists a UMPU test of the null hypothesis $\Po = \{\bs{\theta} \in \Theta : \theta_1 = c\}$ or $\Po = \{\bs{\theta} \in \Theta : \theta_1 \geq c\}$ against the alternative $\Pa = \Theta \backslash \Po$, where $\Theta$ is the natural parameter space. These null hypothesis sets are relatively simple as they represent an exponential family submodel or a closed halfspace in the natural parameter space. Less is known about general conditions for the existence of UMPU tests when $\Po$ is a curved exponential family submodel \cite{BrownExponentialFamilies} of a full exponential family; see  \citep{SenGuptaExponentialUMPU} for examples of UMPU tests in curved exponential families of multivariate exponential distributions. Theorem \ref{thm:PrincipalIdealUBThresholdUMP} provides an existence result for UMPU tests at the unbiasedness threshold for a large class of curved multinomial families.

\begin{example}
For $k > 2$ the null hypothesis set $\mathcal{P}_{0,\delta} = \{\bs{\pi} \in \Delta_{k-1}: \sum_{i = 1}^k (\pi_i-\tfrac{1}{k})^2 = \delta^2\}$ consisting of  probability vectors $\bs{\pi}$ contained within a radius-$\delta$ sphere of the uniform probability vector $(\tfrac{1}{k},\ldots,\tfrac{1}{k})$ is a curved exponential subfamily of the exponential family of multinomial distributions. Assume that $\delta$ is less than $1 - \tfrac{1}{k}$, so that $\Po$ has a non-empty intersection with $\text{relint}(\Delta_{k-1})$. The $(k-1)$-sphere in $\mb{C}^{k}$ does not have a polynomial parameterization, but it does possess a rational parameterization defined by the stereographic projection:
\begin{align*}
    \varphi(z_1,\ldots,z_{k-1}) = \bigg(\frac{\sum_{i = 1}^{k-1}z_i^2 - 1}{\sum_{i = 1}^{k-1}z_i^2 + 1}, \frac{2z_1}{\sum_{i = 1}^{k-1} z_i^2 + 1},\ldots, \frac{2z_{k-1}}{\sum_{i = 1}^{k-1} z_i^2 + 1},  \bigg) \in \mb{S}^{k-1}.
\end{align*}
Lemma \ref{custlem:RationalParamIdeal} in the supplement shows that $I_{\mb{R^{k-1}}}(\Po) = \langle  \sum_{i = 1}^{k-1} (\pi_i-\tfrac{1}{k})^2 + (1 - \sum_{i = 1}^{k-1}\pi_i - \tfrac{1}{k})^2 - \delta^2 \rangle$. By Theorem \ref{thm:PrincipalIdealUBThresholdUMP}, there exists a unique UMPU test of $\Po$ when $n = 4$. The UMPU test with size-$\alpha$ has the power polynomial 
\begin{align*}
    \beta(\bs{\pi}) = c_\alpha \bigg( (\bs{\pi}^\intercal \bs{\pi})^2 - 2 \bs{\pi}^\intercal \bs{\pi} (\delta + \tfrac{1}{k}) + (\delta + \tfrac{1}{k})^2\bigg) + \alpha
\end{align*}
In the supplement we show how to find the maximum $c_\alpha$ and construct the UMPU hypothesis test that corresponds to this power polynomial. 
\end{example}

A similar result to Theorem \ref{thm:PrincipalIdealUBThresholdUMP} can be obtained for semialgebraic null hypotheses. The main distinction is that when $\Po =  \{\bs{\pi} \in \Delta_{k-1}: f(\bs{\pi}) \leq 0\}$ the polynomial $f$ is a separating polynomial, and unlike the algebraic case, this polynomial does not need to be squared as it is already positive on $\Pa$. Assuming some regularity conditions, the unbiasedness threshold of this semialgebraic null hypothesis is thus half the unbiasedness threshold of the corresponding algebraic null hypothesis  $\Po =  \{\bs{\pi} \in \Delta_{k-1}: f(\bs{\pi}) = 0\}$. The intuition here is that a test only needs power in the single region $\{ \bs{\pi} \in \Delta_{k-1}: f(\bs{\pi}) > 0\}$ to test the semialgebraic null hypothesis, while to test the algebraic null hypothesis a test needs non-trivial power in two regions:  $ \{\bs{\pi} \in \Delta_{k-1}: f(\bs{\pi}) < 0\}$ and $ \{\bs{\pi} \in \Delta_{k-1}: f(\bs{\pi}) > 0\}$. Twice as many observations are required to obtain non-trivial power in both of these regions.

% to test the algebraic null hypothesis the power function 

% there are two opposite directions of alternative hypotheses, $ \{\bs{\pi} \in \Delta_{k-1}: f(\bs{\pi}) < 0\}$ and $ \{\bs{\pi} \in \Delta_{k-1}: f(\bs{\pi}) > 0\}$; twice as many observations are needed to have non-trivial power in both directions as compared to the null hypothesis $\Po =  \{\bs{\pi} \in \Delta_{k-1}: f(\bs{\pi}) \leq 0\}$ which only has a single direction of alternative hypotheses. 

% The intuition here is that u the alternative hypothesis extends in two opposite directions: $ \{\bs{\pi} \in \Delta_{k-1}: f(\bs{\pi}) < 0\}$ and $ \{\bs{\pi} \in \Delta_{k-1}: f(\bs{\pi}) > 0\}$; twice as many observations are needed to have non-trivial power in both directions as compared to the null hypothesis $\Po =  \{\bs{\pi} \in \Delta_{k-1}: f(\bs{\pi}) \leq 0\}$ which only has a single direction of alternative hypotheses. 

% The intuition here is that to test the algebraic null hypothesis there are two opposite directions of alternative hypotheses, $ \{\bs{\pi} \in \Delta_{k-1}: f(\bs{\pi}) < 0\}$ and $ \{\bs{\pi} \in \Delta_{k-1}: f(\bs{\pi}) > 0\}$; twice as many observations are needed to have non-trivial power in both directions as compared to the null hypothesis $\Po =  \{\bs{\pi} \in \Delta_{k-1}: f(\bs{\pi}) \leq 0\}$ which only has a single direction of alternative hypotheses. 

\begin{theorem}
    Let $\Po = \{\bs{\pi} \in \Delta_{k-1}: f(\bs{\pi}) \leq 0\}$ and assume that $I_{\mb{R}^{k-1}}( \partial \Po) = \langle f \rangle$ where $\partial \Po$ is the relative boundary of $\Po$ in $\Delta_{k-1}$.  The NTUB and SUB unbiasedness thresholds are both equal to $\mathrm{deg}(f)$. For every size $0<\alpha<1$, when the sample size is $n = \,\mathrm{deg}(f)$ there is a unique, uniformly most powerful NTUB test that is also the unique, uniformly most powerful SUB test. This test has a power function of the form $\beta = c_{\alpha} f + \alpha$, where $c_{\alpha} \in (0,\infty)$ is the largest constant that makes $\beta$ a power polynomial. Every size-$\alpha$ NTUB test of $\Po$ has a power function of the form $cf + \alpha$, $c \in (0,c_\alpha]$ when $n = \mathrm{deg}(f)$.  More generally, if $\beta$ is the power function of a level-$\alpha$, NTUB test for a sample size of $n$ there exists a polynomial $h$ that is non-negative on $\Delta_{k-1}$ with degree $n - \mathrm{deg}(f)$ such that $\beta = fh + \alpha$. For a SUB test $h$ must in addition be strictly positive on $\Pa$.
\end{theorem}

To conclude this section, we use Theorem \ref{thm:PrincipalIdealUBThresholdUMP} to obtain the UB threshold for a union of hypotheses. A union of a collection of null hypotheses, $\Po = \cup_{i = 1}^l \mathcal{P}_i$, each with the principal ideal $I_{\mb{R}^{k-1}}(\mathcal{P}_i) = \langle f_{(i)} \rangle$, also has a principal ideal. As a result, if the $f_{(i)}$ polynomials do not share a common factor, then the UB threshold is the sum of the individual UB thresholds of each submodel $\mathcal{P}_i$. 

\begin{corollary}
\label{cor:UnionofPrincipal}
        Assume that $f_{(1)},\ldots,f_{(l)}$ are distinct, irreducible polynomials, $I_{\mb{R}^{k-1}}(\mathcal{P}_i) = \langle f_{(i)} \rangle$, each $\mathcal{P}_i$ is algebraic, $\overline{\mathrm{int}(\tilde{\Delta}_{k-1}) \cap \mathcal{P}_i} = \mathcal{P}_i$,  and $\nabla f_{(i)}$ is non-zero on $\mathcal{P}_i$. Then $I_{\mb{R}^{k-1}}(\cup_{i = 1}^l \mathcal{P}_i) = \langle \prod_{i = 1}^l f_{(i)} \rangle$ and the SUB and NTUB thresholds of $\Po = \cup_{i = 1}^l \mathcal{P}_i$ are $2 \sum_{i = 1}^l\mathrm{deg}(f_{(i)})$. 
\end{corollary}

One typical application of Corollary \ref{cor:UnionofPrincipal} is when each $\mathcal{P}_i = \{ \bs{\pi} \in \Delta_{k-1}: \sum_{j = 1}^k a_{ij} \pi_j = b_i\} \subsetneq \Delta_{k-1}$ is a distinct hyperplane that intersects the relative interior of the simplex. The resulting UB threshold of the union of hyperplanes $\cup_{i = 1}^l \mathcal{P}_i$ is $2l$. For example, the hypothesis $\Po = \{ \bs{\pi} \in \Delta_{k-1}: \exists i \neq j, \pi_i = \pi_j\}$ 
that at least two coordinates of the vector $\bs{\pi}$ are equal has an UB threshold of $2{ k \choose 2}$. This is a setting where there is a very real possibility that the sample size does not exceed the UB threshold; if $k = 100$ the UB threshold is approximately $10$,$000$.

More generally, it is possible to obtain lower bounds for the NTUB and SUB thresholds of a union of null hypotheses via the following lemma.
\begin{lemma}
    Let $\mathcal{P}_i$ have the NTUB and SUB thresholds of $d_{1i}, d_{2i}$ respectively with corresponding separating polynomials $f_{(1i)}$ and $f_{(2i)}$ for $i = 1,\ldots,l$. The polynomials $\prod_{i = 1}^l f_{(1i)}$ and $\prod_{i = 1}^l f_{(2i)}$ are NTUB and SUB separating polynomials with degrees $\sum_{i = 1}^l d_{1i}$,  $\sum_{i = 1}^l d_{2i}$ for the union null hypothesis $\Po = \cup_{i = 1}^l \mathcal{P}_i$.
\end{lemma}

\subsubsection{Non-principal Vanishing Ideals}
If $I_{\mb{R}^{k-1}}(\Po) = \langle f_{(1)},\ldots,f_{(m)} \rangle$ for an algebraic null hypothesis then $\sum_{i = 1}^m h_{(i)}^2f_{(i)}^2$ is an example of a SOS in the vanishing ideal. A wide selection of separating polynomials can be constructed by varying the coefficient polynomials $h_{(1)},\ldots, h_{(m)}$. The $h_{(i)}$ can be viewed as weights that accentuate the power in certain regions in the alternative hypothesis space. Specifically, the alternative hypothesis can be written as a union of sets $\Delta_{k-1}\backslash V_{\mb{R}^{k-1}}(f_{(i)})$ as follows:
\begin{align}
\label{eqn:AlgAlternativeAsUnion}
    \Pa = \Delta_{k-1} \backslash \Po = \Delta_{k-1} \big\backslash \bigg(\bigcap_{i = 1}^m V_{\mb{R}^{k-1}}(f_{(i)}) \bigg) = \bigcup_{i = 1}^m \bigg(\Delta_{k-1}  \backslash V_{\mb{R}^{k-1}}(f_{(i)}) \bigg).
\end{align}
The separating polynomial $f_{(i)}^2$ corresponds to a SUB test for testing the null hypothesis $\Delta_{k-1}  \cap V_{\mb{R}^{k-1}}(f_{(i)})$ against the alternative $\Delta_{k-1}  \backslash V_{\mb{R}^{k-1}}(f_{(i)})$. By taking an $h_{(i)}$-weighted sum of these $f_{(i)}$, each respective term in a separating power polynomial of the form $\sum_{i = 1}^m h_{(i)}^2f_{(i)}^2$ can be viewed as increasing the power of the corresponding test against the alternative $\Delta_{k-1}  \backslash V_{\mb{R}^{k-1}}(f_{(i)})$. In the simplest case the $h_{(i)}^2$ are non-negative scalars and $\sum_{i = 1}^m h_{(i)}^2f_{(i)}^2$ is a sum of squares with non-negative weights. There is a tradeoff here as once a separating polynomial is re-normalized to a power polynomial, improving the power in one region of the alternative hypothesis space may come at the expense of reducing the power of the test in other regions.

\begin{example}
    Let $\Po = \{ \bs{\pi} \in \Delta_{k-1}: \bl{a}_i^\intercal \bs{\pi} = b_i, i = 1,\ldots,l\}$ be an affine null hypothesis set that intersects the relative interior of $\Delta_{k-1}$. The unbiasedness threshold for $\Po$ is two as no degree one polynomial is  separating, for if $f$ has degree one then $\{\bs{\pi}: f(\bs{\pi}) < 0\}$ is a halfspace that cannot be contained in $\Po$. When $n = 2$ a collection of SUB separating polynomials of the form $\sum_{i = 1}^l c_i (\bl{a}_i^\intercal \bs{\pi} - b_i)^2$ can be obtained by modifying the various weights $c_i > 0$; the larger the weight $c_i$ is, the higher the power of the corresponding test against alternatives where $(\bl{a}_i^\intercal \bs{\pi} - b_i)^2$ is large. A specific example of an affine null hypothesis is the null hypothesis of symmetry \citep{bowker1948symmetrytest} for square contingency tables where $\bs{\pi} \in \mb{R}^{p \times p}$ must satisfy $\pi_{ij} - \pi_{ji} = 0$ for $i,j \in [p]$. 
\end{example}

Another example of a non-principal vanishing ideal is the $p \times q$ independence model in Example \ref{ex:IndepModelpq} that has a vanishing ideal generated by the polynomials $f_{(i)} = \det(\bs{\pi}_{I,J})$. A separating power polynomial of the form $\sum_{I,J} c_{I,J}^2  \det(\bs{\pi}_{I,J})^2$, $c_{I,J} > 0$ has a weight of $c_{I,J}^2$ that imparts power against the alternative that the $2 \times 2$ submatrix $\bs{\pi}_{I,J}$ has rank larger than one. Using the SOS upper bound from Theorem \ref{thm:UBThresholdSOSUpperBounds}, the aim is to show that this upper bound is also a lower bound, meaning that there do not exist polynomials with degree less than the upper bound that satisfy the conditions of Corollary \ref{cor:ExistenceforAlgebraic}. The following result uses a symmetrization argument to establish the UB thresholds for contingency tables with bounded rank; these tables generalize tables with independent rows and column variables.

\begin{theorem}
\label{thm:BoundedRankUBThreshold}
    The NTUB and SUB unbiasedness thresholds for the null hypothesis $\Po = \{ \bs{\pi} \in \Delta_{pq-1}: \bs{\pi} \; \mathrm{has \; rank \; less \; than } \; r \}$ on the space of $p \times q$ dimensional contingency tables is equal to $2r$. When $n = 2r$ one such SUB power polynomial has the form $\beta(\bs{\pi}) = c\sum_{I,J} \det(\bs{\pi}_{I,J})^2 + \alpha$ for some $c > 0$, where $I,J$ are index sets of size $r$. 
    % The NTUB and SUB unbiasedness threshold for the null hypothesis $\Po = \{ \bs{\pi} \in \Delta_{pq-1}: \bs{\pi} \; \mathrm{has \; rank \; less \; than } \; r \}$ on the space of $p \times q$ dimensional contingency tables is equal to $2r$. A separating polynomial corresponding to such a SUB test has the form $\Tilde{\beta} = \sum_{I,J} c_{I,J} \det(\bs{\pi}_{I,J})^2$ where the index sets $I,J$ have size $r$, $\bs{\pi}_{I,J}$ is any $r \times r$ submatrix of $\bs{\pi}$, and $c_{I,J} >  0$. 
\end{theorem}
The bounded rank hypothesis is closely related to finite mixtures of the independence model. An $(r-1)$-component mixture of the independence model has a contingency table of the form 
\begin{align}
\label{eqn:IndepMixtureModel}
  \bs{\pi} =  \sum_{i = 1}^{r-1} c_i \bl{a}^{(i)} (\bl{b}^{(i)})^\intercal, \;\; \bl{a}^{(i)} \in \Delta_{p-1},\;\; \bl{b}^{(i)} \in \Delta_{q-1}, \; \bl{c}^\intercal 
 = (c_1,\ldots,c_r) \in \Delta_{r-2}. 
\end{align}
Any mixture model with less than $r$ independent mixture components is contained in the null hypothesis set in Theorem \ref{thm:BoundedRankUBThreshold}. However, it is not true that every point in the null hypothesis set can be represented in the form \eqref{eqn:IndepMixtureModel} due to the non-negativity constraints on $\bl{a}^{(i)}$, $\bl{b}^{(i)}$, and $\bl{c}$  \citep{NonnegativeRankCohen}.

The next result is the counterpart to Theorem \ref{thm:PrincipalIdealUBThresholdUMP} that provides a representation of UB power polynomials for non-principal ideals.

\begin{theorem}
\label{thm:NonPrincipalPowerPolyRepresentation}
  Suppose that $\Po$ is algebraic, $\overline{\text{int}(\tilde{\Delta}_{k-1}) \cap \Po} = \Po$, $I_{\mb{R}^{k-1}}(\Po) = \langle f_{(1)},\ldots,f_{(m)} \rangle$, and the matrix $[\nabla f_{(1)} \cdots \nabla f_{(m)}] \in \mb{R}^{k-1 \times m}$ has rank $m$ on $\Po$, then every NTUB test of $\Po$ for a sample size of $n$ has a power polynomial of the form $\beta(\bs{\pi}) = \bl{f}^\intercal(\bs{\pi}) \bl{H}(\bs{\pi}) \bl{f}(\bs{\pi}) + \alpha$. Here $\bl{f}^\intercal = (f_{(1)},\ldots,f_{(m)})$, the entries of $\bl{H}$ are polynomials $h_{(ij)}$ with $\bl{f}^\intercal \bl{H} \bl{f} \geq 0$ on $\Delta_{k-1}$. If the test is SUB we require in addition that $\bl{f}^\intercal \bl{H} \bl{f} > 0$ on $\Pa$. If the $f_{(1)},\ldots,f_{(m)}$ are a Gr\"obner basis with respect to a graded monomial ordering then the representation of $\tilde{\beta}$ can be constructed so that $\mathrm{deg}(\tilde{\beta}) = \max_{i,j} \mathrm{deg}(h_{(ij)} f_{(i)}f_{(j)})$. The NTUB threshold is equal to the SOS threshold described in Theorem \ref{thm:UBThresholdSOSUpperBounds}. When the degrees of $f_{(1)},\ldots,f_{(m)}$ are all the same the NTUB threshold is also equal to the SUB threshold. 
\end{theorem}

The non-negativity condition in the above theorem is more subtle than in Theorem \ref{thm:PrincipalIdealUBThresholdUMP}; it is sufficient, but not necessary, for $\bl{H}(\bs{\pi)}$ to be a positive semidefinite matrix for all $\bs{\pi} \in \Delta_{k-1}$. When $\bl{H}$ is a positive semidefinite matrix of scalars the separating polynomial $\tilde{\beta} = \bl{f}^\intercal \bl{H} \bl{f}$ is a sum of squares and the representation is known as the Gram matrix method \citep[Sec 3.1]{PowersCertificatesofPositivity} for finding a sum of squares.

\section{The Existence of UMPU Tests}
\label{sec:UMPUTesting}
 In Theorem \ref{thm:PrincipalIdealUBThresholdUMP} it was shown that any algebraic hypothesis with a principal ideal has an UMPU test when the sample size equals the UB threshold. The reason that there is a UMPU test is simply because there is only one UB test. This raises the question of how to construct UB tests with adequate power for sample sizes larger than the UB threshold? Machinery for constructing UMPU tests for models with principal ideals is detailed in this section, providing proofs that UMPU tests do not always exist for algebraic null hypotheses, and that there can exist UMPU tests for null hypotheses that are not exponential families.

As an illustrative example, consider the null hypothesis $\Po = \{(\pi_1,\pi_2) \in \Delta_1: \pi_1 = \pi_2\}$ in the model $\text{Binomial}(\pi_1,n)$. Using Theorem \ref{thm:PrincipalIdealUBThresholdUMP}, any NTUB power polynomial for $\Po$ must be of the form
\begin{align}
\label{eqn:BinomialHalfPowerFnc}
    \beta(\pi_1,\pi_2) = f^2 h + \alpha = (\pi_1 - \pi_2)^2 \bigg(\sum_{i = 0}^{n'} h_i \pi_1^i \pi_2^{n'-i}\bigg) + \alpha,
\end{align}
where $n' \coloneqq n-2$. To maximize the power, the coefficients $h_i \in \mb{R}$ of the polynomial $h$ should be taken to be as large as possible since $h$ is pointwise non-decreasing on $\Delta_{k-1}$ with respect to the $h_i$. The coefficients $h_i$ must also be constrained so that $\beta$ satisfies the box-constraints in Lemma \ref{lem:PowerFunctionCharacterization}.
\begin{definition}[Coefficient polytope]
If $I_{\mb{R}^{k-1}}(\Po) = \langle f \rangle$ and $\Tilde{f}$ is the homogenization of $f$, the coefficient polytope $\mc{C}_{n,\alpha}(\Po)$ is the subset of monomial coefficients of degree $n' \coloneqq n - 2\text{deg}(f)$, $(h_J) \in \mb{R}^{ {n' + k - 1 \choose k-1}}$ that satisfy the linear inequality constraints
\begin{align}
\label{eqn:ConstraintPolytopeEquations}
  - {n \choose L} \alpha \leq  \sum_{(I,J): I + J = L} (\Tilde{f}^2)_I h_J \leq {n \choose L} (1-\alpha),
\end{align}
for all multiindices $L$ of size $n$. The scalar $(\Tilde{f}^2)_I$ is the coefficient of $\bs{\pi}^I$ in $\tilde{f}^2$. 
\end{definition}

The next observation is that if $\beta$ in \eqref{eqn:BinomialHalfPowerFnc} is UMPU then $\beta(1,0) = h_{n'} + \alpha$ and $\beta(0,1) = h_0 + \alpha$ must be as large as possible while remaining in the coefficient polytope. Specifically, if $\text{Proj}_{(0,n')}( \mc{C}_{n,\alpha}(\Po))$ is the projection of the coefficient polytope onto $(h_0,h_{n'})$, then in order for a UMPU test to exist this projected set must have a maximum point with respect to the componentwise ordering of vectors in $\mb{R}^2$. The constraints involving $h_0$ are
\begin{align*}
    -\alpha \leq &h_0 \leq (1-\alpha)
    \\
    -n \alpha  \leq &h_1-2h_0 \leq  n(1-\alpha)
    \\
    -{n \choose 2} \alpha \leq & h_2 -2 h_1 + h_0 \leq {n \choose 2}(1-\alpha). 
\end{align*}
% \textcolor{red}{It is more complicated than this --- I think the largest $\alpha$ is not }

% From the first constraint $h_0$ is at most $(1-\alpha)$, and it can be shown that $h_0 = h_n = 2^{n-1}\alpha$

% $(h_0, h_{n'}) = (1-\alpha,1-\alpha)$ is in $\text{Proj}_{(0,n')}( \mc{C}_{n,\alpha}(\Po))$ \textcolor{red}{Have to be a bit careful about this --- this requires that the level is larger than $(1/2)^n$ I believe. Basically, when this constraint is satisfied with equality it means we reject when $x = 0,n$ --- the size of the test must be so that this possibility can occur. This perspective might give us a nice sequential method for finding a UMPU test though when it exists via peeling.}. 
To find a candidate UMPU polynomial we set $(h_0,h_{n'}) = (h_0^*, h_{n'}^*) \in \text{Proj}_{(0,n')}( \mc{C}_{n,\alpha}(\Po))$ to the componentwise maximum in \eqref{eqn:BinomialHalfPowerFnc} and maximize $\beta^* = \beta - \tilde{f}^2(h_0^*\pi_1^{n'} + h_n^*\pi_2^{n'}) = (\pi_1 - \pi_2)^2 \bigg(\sum_{i = 1}^{n'-1} h_i \pi_1^i \pi_2^{n'-i}\bigg) + \alpha$ with respect to $h_1,\ldots,h_{n'-1}$. As $t \rightarrow 0,1$ we have $\beta^*(1-t,t) \rightarrow h_{1},h_{n'-1}$, and any UMPU test must maximize the pair $(h_1,h_{n'-1})$ over $\text{Proj}_{(1,n'-1)}(\mc{C}_{n,\alpha}(\Po) \cap \{h:(h_0,h_{n'}) = (h_0^*, h_{n'}^*) \})$ with respect to the componentwise order. Once $h_1,h_{n'-1}$ are found, the process repeats with $h_2,h_{n'-2}$, sequentially peeling off the highest degree coefficients remaining in $h$ until every coefficient of $h$ is determined. The ordering that monomials should be peeled off does not depend on the choice of $\Po$ and is related to following convex peeling \citep{BarnettPeelingDepth} of the simplex, motivated in part by Polya's theorem \citep[Sec 9.1]{PowersCertificatesofPositivity}.

\begin{definition}[Convex peeling]
\label{def:ConvexPeeling}
    Let $\mc{S}^{(0)} = \{(i_1,\ldots,i_k): i_k \in \mb{N}, \sum_{j = 1}^{k} i_j = n'\}$ and define $\Delta^{(0)}_{k-1} = \text{conv}(\mc{S}^{(0)}) = \Delta_{k-1}$. If $\mc{V}^{(i-1)} \subset \mc{S}^{(i-1)}$ are the vertices of $\Delta^{(i-1)}_{k-1}$, define for $i \geq 1$
    \begin{align*}
        \mc{S}^{(i)} & = \mc{S}^{(i-1)} \backslash V^{(i-1)}
        \\
        \Delta^{(i)}_{k-1} & = \text{conv}(\mc{S}^{(i)}).
    \end{align*} 
 The set $\Delta_{k-1}^{(i)}$ is the $i$th convex peeling of $\Delta_{k-1}$ with respect to the integer lattice points in $\Delta_{k-1}$.   
\end{definition}
Monomials with exponents in the vertex sets $V^{(i)}$ for $i = 0,1,2,3$ and $k = 3$ respectively have the form
\begin{align*}
    \{\mc{\pi}_i^{n'}: i \in [k]\}, \;  \{\mc{\pi}_i^{n'-1}\pi_j: i \neq j\}, \; \{\mc{\pi}_i^{n'-2}\pi_j^2: i \neq j\},  \{ \pi_i^{n'-2}\pi_j \pi_l: i \neq j \neq l\} \cup   \{ \pi_i^{n'-3}\pi_j^3: i \neq j\}.
\end{align*}
For a $1$-dimensional simplex, the $i$th vertex set is $\mc{V}^{(i)} = \{i+1,n'-i\}$, corresponding to the order in which the $h_i$ coefficients were determined in \eqref{eqn:BinomialHalfPowerFnc}.
\begin{theorem}
\label{thm:ConstraintNeccPolytopeConditions}
A necessary condition for a level-$\alpha$  UMPU test for $\Po$ to exist is that
    \begin{enumerate}
        \item For $V^{(0)} = \{I_{01},\ldots,I_{0l_0}\}$ the set $\text{Proj}_{(I_{01},\ldots,I_{0l_0})}\big( \mc{C}_{n,\alpha}(\Po)\big)$ has a componentwise maximum element $(h_{I_{01}}^*,\ldots,h_{I_{0l_0}}^*)$. The vertex set $V^{(0)}$ is described in Definition \ref{def:ConvexPeeling}.
        \item Inductively, for every $j$ with $V^{(j)} = \{I_{j1},\ldots,I_{jl_j}\}$ the set
        \begin{align*}
\text{Proj}_{(I_{j1},\ldots,I_{jl_j})}\big( \mc{C}_{n,\alpha}(\Po) \cap \{ (h_I) : h_{I_{j'a}} = h_{I_{j'a}}^*, \forall j' < j, a \leq l_{j'}\}\big)
        \end{align*}set 
has a componentwise maximum $(h^*_{I_{j1}},\ldots,h^*_{I_{jl_j}})$. 
    \end{enumerate}
    If a UMPU test exists it has the power polynomial $f^2h^* + \alpha$, where $h^*$ is the polynomial  corresponding to the above polynomial coefficients $(h^*_I)$.
\end{theorem}

The coefficient polytope is defined by the hyperplane equations \eqref{eqn:ConstraintPolytopeEquations}. For determining the projections of the coefficient polytope a vertex representation \citep[Sec 1.1]{ZieglerLectureonPolytopes} is more useful, since  $\text{Proj}_{(I_1,\ldots,I_l)}\big(\text{conv}(\bl{v}_1,\ldots,\bl{v}_m)\big) = \text{conv}\big(\text{Proj}_{(I_1,\ldots,I_l)}(\bl{v}_1,\ldots,\bl{v}_m)\big)$. The conditions in Theorem \ref{thm:ConstraintNeccPolytopeConditions} can be checked by finding the vertex representation of $\mc{C}_{n,\alpha}(\Po)$ and examining the structure of the vertices.

\begin{lemma}
\label{lem:UMPUSufficientConditions}
    If $\bl{v}_1,\ldots,\bl{v}_m$ are the vertices of $\mc{C}_{n,\alpha}(\Po)$ and there exists a vertex that is the maximum in the componentwise ordering, then there exists a UMPU test with $(h_I)$ equal to the maximizing vector. When $n' = 1$ this is a also a necessary condition for the existence of a UMPU test.
\end{lemma}

\begin{example}
The hyperplane representation of $\mc{C}_{n,\alpha}(\Po)$ for $\Po = \{\bs{\pi} \in \Delta_{2}: \pi_1 + \pi_2 = \pi_3\}$ can be converted to its vertex representation by solving a linear program. The   vertices of $\mc{C}_{n = 3,\alpha = 0.05}(\Po)$ are found to be
\begin{align*}
\begin{bmatrix}
    -0.05
    \\
    0.05
    \\
    -0.05
\end{bmatrix},
\begin{bmatrix}
    0.05
    \\
    0.1
    \\
    -0.05
\end{bmatrix},
\begin{bmatrix}
    -0.05
    \\
    -0.05
    \\
    -0.05
\end{bmatrix},
\begin{bmatrix}
    0.05
    \\
    -0.05
    \\
    -0.05
\end{bmatrix},
\begin{bmatrix}
    -0.05
    \\
    -0.05
    \\
    0.05
\end{bmatrix},'
\begin{bmatrix}
    0.05
    \\
    -0.05
    \\
    0.1
\end{bmatrix},
\begin{bmatrix}
   -0.05
    \\
    0.05
    \\
    0.05
\end{bmatrix},
\begin{bmatrix}
    0.15
    \\
    0.15
    \\
    0.15
\end{bmatrix}.
  %   -0.05   0.05  -0.05
 % 0.05   0.1   -0.05
% -0.05  -0.05  -0.05
%  0.05  -0.05  -0.05
% -0.05  -0.05   0.05
%  0.05  -0.05   0.1
% -0.05   0.05   0.05
 % 0.15   0.15   0.15
\end{align*}
The last vertex is the componentwise maximum and so $h(\bs{\pi}) = 0.15(\pi_1 + \pi_2 - \pi_3)^2(\pi_1 + \pi_2 + \pi_3) + 0.05 = 0.15(\pi_1 + \pi_2 - \pi_3)^2 + 0.05$ is the power function of the UMPU test of $\Po$.  As the model $\Po$ is not an exponential family, due to Lemma 2.6 in \citep{InfortmationGeometryAy}, the existence of a UMPU test is not a consequence of the existence of UMPU tests for log-linear null hypotheses that was mentioned in Section \ref{sec:ExistenceResults}. 

Suppose that $\Po$ is modified to $\{\bs{\pi} \in \Delta_{2}: 2\pi_1 + \pi_2 = \pi_3\}$. The vertices of $\mc{C}_{n = 3,\alpha = 0.05}(\Po)$ are
\begin{align*}
\begin{bmatrix}
    -0.05
    \\
    -0.025
    \\
    -0.0125
\end{bmatrix},
\begin{bmatrix}
    0.0625
    \\
    -0.025
    \\
    0.1
\end{bmatrix},
\begin{bmatrix}
    -0.0375
    \\
    -0.0375
    \\
  0
\end{bmatrix},
\begin{bmatrix}
    0.0125
    \\
    -0.05
    \\
    0.05
\end{bmatrix},
\begin{bmatrix}
    0.05
    \\
    -0.5
    \\
    0.0875
\end{bmatrix},
\begin{bmatrix}
    0.0375
    \\
    -0.0375
    \\
    0
\end{bmatrix},
\begin{bmatrix}
    0.03438
    \\
    -0.025
    \\
    -0.0125
\end{bmatrix},
\\
\begin{bmatrix}
    0.05
    \\
    -0.05
    \\
    0.05
\end{bmatrix},
\begin{bmatrix}
    -0.03438
    \\
    0.09219
    \\
    -0.0125
\end{bmatrix},
\begin{bmatrix}
    -0.0125
    \\
    0.06875
    \\
    -0.0125
\end{bmatrix},
\begin{bmatrix}
    0.05
    \\
    0.1
    \\
    0.05
\end{bmatrix},
\begin{bmatrix}
    0.0625
    \\
    0.0875
    \\
    0.1
\end{bmatrix},
\begin{bmatrix}
    -0.05
    \\
    0.03125
    \\
    -0.0125
\end{bmatrix}.
 %    -0.05     -0.025    -0.0125
%  0.0625   -0.025     0.1
% -0.0375   -0.0375   -1.692e-18
 % 0.0125   -0.05      0.05
%  0.05     -0.05      0.0875
 % 0.0375   -0.0375   -2.961e-18
%  0.03438  -0.025    -0.0125
%  0.05     -0.05      0.05
%  0.03438   0.09219  -0.0125
% -0.0125    0.06875  -0.0125
%  0.05      0.1       0.05
%  0.0625    0.0875    0.1
% -0.05      0.03125  -0.0125
\end{align*}
There is no componentwise maximum vertex since the third-last vertex maximizes the second coordinate, but does not maximize the first and second coordinates. By Lemma \ref{lem:UMPUSufficientConditions} we conclude that $\Po$ does not have a UMPU test at this size and sample size. 
\end{example}

\begin{customexample}{8} \textbf{(continued)}
\label{ex:SphericalExampleContinued}
 The coefficient polytope vertices for the spherical hypothesis are provided in the supplement for a radius of $\delta^2 = 1/6$, $\alpha = 0.05$ and $n = 5$, where it is seen that this hypothesis has a UMPU test. For $n = 6$ the vertices satisfy the necessary conditions in Theorem \ref{thm:ConstraintNeccPolytopeConditions}, but not the sufficient conditions in Lemma \ref{lem:UMPUSufficientConditions}. It is possible to compute the vertices of $\mc{C}_{n,\alpha}(\Po)$ for larger $n$, although the number of vertices grows large very quickly. When $n = 8$ the coefficient polytope is defined by $72$ halfspaces and has $17{,}267$ vertices. However, in many problems the coefficient polytope has symmetries that arise from symmetries in $\Po$. These symmetries can be utilized to simplify the problem of finding a UMPU test. 
\end{customexample}

An overview of how the ideas above can be extended to non-principal ideals is now provided. Theorem \ref{thm:NonPrincipalPowerPolyRepresentation} shows that unbiased power functions have a representation of the form $\beta = \sum_{i,j} \tilde{f}_{(i)} \tilde{f}_{(j)} h_{(ij)} + \alpha$. Optimizing this expression over every $h_{(ij)}$ is challenging; however, in many cases where ideal is non-principal there is symmetry present in $\Po$. Precisely, there is some subgroup $G$ of the symmetric group $\mc{S}_k$ that fixes $\Po$ in the sense that for all  $\tau \in G$ and $\bs{\pi}_0 \in \Po$ we have $\tau(\bs{\pi}_0) \in \Po$.  This symmetry can be used to place additional constraints on $\beta$ by requiring that $\beta$ corresponds to an invariant test. A test $\phi$ is invariant under $G$ if it is a symmetric function with $\phi(X_1,\ldots,X_k) = \phi(X_{\tau(1)},\ldots,X_{\tau(k)})$ for all $\tau \in G$ \citep[Ch 6]{TSPLehmannRomano}. Invariant tests give rise to power polynomials that are symmetric functions and vice versa.

As an example, consider the symmetry model \citep{bowker1948symmetrytest} for a $p\times p$ contingency table: $\Po = \{\bs{\pi} \in \Delta_{p^2-1}: \bs{\pi}^\intercal = \bs{\pi}, \forall \, i,j \in [p]\}$. This null hypothesis set is invariant under every transposition that sends $\pi_{ij} \mapsto \pi_{ji}$ and fixes every other coordinate, and also under every index relabeling of the form $\pi_{ij} \mapsto \pi_{\sigma(i)\sigma(j)}$ for any $\sigma$ in the symmetric group $\mc{S}_p$. If we restrict to invariant tests the class of unbiased power polynomials has the following form.

\begin{lemma}
    The power polynomial $\beta$ of every unbiased, invariant test of the symmetry hypothesis $\Po = \{\bs{\pi} \in \Delta_{p^2-1}: \pi_{ij} = \pi_{ji}\}$ can be written in the form
    \begin{align*}
        \beta(\bs{\pi}) = \sum_{i,j = 1}^p (\pi_{ij} - \pi_{ji})^2 h_{(ij)} + \alpha,
    \end{align*}
    where the polynomial $h_{(ij)}$ is invariant under the transposition $\pi_{kl} \mapsto \pi_{lk}$ for all $l,k \in [p]$.
\end{lemma}
Furthermore, it can also be shown that all of the $h_{(ij)}$ polynomials are equal to each other up to a relabeling of the indices. Hence, invariance reduces the problem of finding a uniformly most powerful unbiased and invariant test to that of finding a single $h_{(ij)}$ polynomial. As in the principal ideal case, the coefficients of $h_{(ij)}$ should be taken to be as large as possible subject to the constraint that $\beta$ must be a valid power polynomial.

\section{Conclusion}
In this paper we have addressed existence questions for unbiased and UMPU tests, and have clarified how the existence of such tests can depend on both the size and sample size. The approach taken emphasizes studying the power polynomials directly, instead of indirectly through a test function. Many avenues for future work are raised by the present article.

Power functions for sampling distributions beyond multinomial distributions may be amenable to analyses of the power function that are similar to what is done here. This is especially the case for distributions on finite sample spaces, such as for product multinomial sampling mentioned in Section \ref{sec:PowerPolynomialOtherSampling}. Semialgebraically constrained Gaussian distributions play a central role in algebraic statistics and it would be of interest to develop criteria for whether unbiased tests exist for such distributions. The space of power functions is much larger when the sample space is not discrete, suggesting that NTUB tests will exist for large classes of null hypotheses. 

A second direction for future work is to find UMPU, or at least admissible, tests by optimizing over the space of power polynomials. As illustrated in Section \ref{sec:UMPUTesting} directly finding the coefficients of $h$ quickly becomes cumbersome as the sample size grows. By analyzing the coefficient polytope, exploiting symmetries that are present
% , and potentially using fast Fourier transform methods \citep{hirji1997reviewFFTExact}, 
it is hoped that algorithms can be developed for finding UMPU tests in a runtime that is polynomial in the sample size. Such algorithms would certainly depend on the algebraic structure of the null hypothesis, as this determines the coefficient polytope.   

To close, we mention a few open questions regarding the present work. When the ideal associated with a null hypothesis is non-principal the form for power functions in Theorem \ref{thm:NonPrincipalPowerPolyRepresentation} illustrates that each of the polynomials $f_{(i)}$ can be weighted differently in an NTUB test. Does this imply that no UMPU test exists whenever a null hypothesis has a non-principal ideal? Even when the null hypothesis ideal is principal it is not obvious which generating polynomials lead to hypotheses that possess UMPU tests. In Example \ref{ex:SphericalExampleContinued} it was demonstrated that a UMPU test may exist for some sample size but not for a larger sample size. Is it true in general that UMPU tests are downward-closed with respect to the sample size in the sense that if a UMPU test exists for some sample size then it also exists for all smaller sample sizes that are larger than the UB threshold?

\section*{Acknowledgments}
The author acknowledges the support of the Natural Sciences and Engineering Research Council of Canada (NSERC), [RGPIN-2025-03968, DGECR-2025-00237].

Cette recherche a été financée par le Conseil de recherches en sciences naturelles et en génie du Canada (CRSNG), [RGPIN-2025-03968, DGECR-2025-00237].

\bibliographystyle{abbrv}
\bibliography{biblio}

\begin{thebibliography}{10}

\bibitem{adcock1997sample}
C.~Adcock.
\newblock Sample size determination: a review.
\newblock {\em Journal of the Royal Statistical Society: Series D (The Statistician)}, 46(2):261--283, 1997.

\bibitem{agresti2001exact}
A.~Agresti.
\newblock Exact inference for categorical data: recent advances and continuing controversies.
\newblock {\em Statistics in medicine}, 20(17-18):2709--2722, 2001.

\bibitem{amiri2017comparisonofindeptests}
S.~Amiri and R.~Modarres.
\newblock Comparison of tests of contingency tables.
\newblock {\em Journal of biopharmaceutical statistics}, 27(5):784--796, 2017.

\bibitem{MarkovBasesBook}
S.~Aoki, H.~Hara, and A.~Takemura.
\newblock {\em Markov bases in algebraic statistics}.
\newblock Springer Series in Statistics. Springer, New York, 2012.

\bibitem{InfortmationGeometryAy}
N.~Ay, J.~Jost, H.~V. L\^e, and L.~Schwachh\"ofer.
\newblock {\em Information geometry}, volume~64 of {\em Ergebnisse der Mathematik und ihrer Grenzgebiete. 3. Folge. A Series of Modern Surveys in Mathematics [Results in Mathematics and Related Areas. 3rd Series. A Series of Modern Surveys in Mathematics]}.
\newblock Springer, Cham, 2017.

\bibitem{ExistenceofUMPUTests}
S.~K. Bar-Lev and B.~Reiser.
\newblock An exponential subfamily which admits {UMPU} tests based on a single test statistic.
\newblock {\em Ann. Statist.}, 10(3):979--989, 1982.

\bibitem{BarnettPeelingDepth}
V.~Barnett.
\newblock The ordering of multivariate data.
\newblock {\em J. Roy. Statist. Soc. Ser. A}, 139(3):318--355, 1976.

\bibitem{bhapkar1979Symmetrytests}
V.~P. Bhapkar.
\newblock On tests of marginal symmetry and quasi-symmetry in two and three-dimensional contingency tables.
\newblock {\em Biometrics}, pages 417--426, 1979.

\bibitem{BirchConditionalIndependenceTests}
M.~W. Birch.
\newblock The detection of partial association. {II}. {T}he general case.
\newblock {\em J. Roy. Statist. Soc. Ser. B}, 27:111--124, 1965.

\bibitem{bowker1948symmetrytest}
A.~H. Bowker.
\newblock A test for symmetry in contingency tables.
\newblock {\em Journal of the american statistical association}, 43(244):572--574, 1948.

\bibitem{BrownExponentialFamilies}
L.~D. Brown.
\newblock {\em Fundamentals of statistical exponential families with applications in statistical decision theory}, volume~9 of {\em Institute of Mathematical Statistics Lecture Notes---Monograph Series}.
\newblock Institute of Mathematical Statistics, Hayward, CA, 1986.

\bibitem{ClinicalTrialSampleSize}
Y.~Cheng, F.~Su, and D.~A. Berry.
\newblock Choosing sample size for a clinical trial using decision analysis.
\newblock {\em Biometrika}, 90(4):923--936, 2003.

\bibitem{christensenLogLinearBook}
R.~Christensen.
\newblock {\em Log-linear models and logistic regression}.
\newblock Springer Texts in Statistics. Springer-Verlag, New York, second edition, 1997.

\bibitem{NonnegativeRankCohen}
J.~E. Cohen and U.~G. Rothblum.
\newblock Nonnegative ranks, decompositions, and factorizations of nonnegative matrices.
\newblock {\em Linear Algebra Appl.}, 190:149--168, 1993.

\bibitem{ConwayComplexVar}
J.~B. Conway.
\newblock {\em Functions of one complex variable}, volume~11 of {\em Graduate Texts in Mathematics}.
\newblock Springer-Verlag, New York-Berlin, second edition, 1978.

\bibitem{cox2013ideals}
D.~Cox, J.~Little, and D.~OShea.
\newblock {\em Ideals, varieties, and algorithms: an introduction to computational algebraic geometry and commutative algebra}.
\newblock Springer Science \& Business Media, 2013.

\bibitem{DerksenKron}
H.~Derksen and V.~Makam.
\newblock Maximum likelihood estimation for matrix normal models via quiver representations.
\newblock {\em SIAM J. Appl. Algebra Geom.}, 5(2):338--365, 2021.

\bibitem{DiaconisSturmfelsMarkov}
P.~Diaconis and B.~Sturmfels.
\newblock Algebraic algorithms for sampling from conditional distributions.
\newblock {\em Ann. Statist.}, 26(1):363--397, 1998.

\bibitem{DobraMarkovBasesGraphModel}
A.~Dobra.
\newblock Markov bases for decomposable graphical models.
\newblock {\em Bernoulli}, 9(6):1093--1108, 2003.

\bibitem{DrtonLRTSingularities}
M.~Drton.
\newblock Likelihood ratio tests and singularities.
\newblock {\em Ann. Statist.}, 37(2):979--1012, 2009.

\bibitem{DrtonKurikiKron}
M.~Drton, S.~Kuriki, and P.~Hoff.
\newblock Existence and uniqueness of the {K}ronecker covariance {MLE}.
\newblock {\em Ann. Statist.}, 49(5):2721--2754, 2021.

\bibitem{LecturesonAlgebraicDrtonSullivant}
M.~Drton, B.~Sturmfels, and S.~Sullivant.
\newblock {\em Lectures on algebraic statistics}, volume~39 of {\em Oberwolfach Seminars}.
\newblock Birkh\"auser Verlag, Basel, 2009.

\bibitem{DrtonWaldTestSingularities}
M.~Drton and H.~Xiao.
\newblock Wald tests of singular hypotheses.
\newblock {\em Bernoulli}, 22(1):38--59, 2016.

\bibitem{SturmfelsGeigerMeekToricAlgebraGraphical}
D.~Geiger, C.~Meek, and B.~Sturmfels.
\newblock On the toric algebra of graphical models.
\newblock {\em Ann. Statist.}, 34(3):1463--1492, 2006.

\bibitem{AlgStatPistoneWynn}
P.~Gibilisco, E.~Riccomagno, M.~P. Rogantin, and H.~P. Wynn, editors.
\newblock {\em Algebraic and geometric methods in statistics}.
\newblock Cambridge University Press, Cambridge, 2010.

\bibitem{GraybillSampleSizeWidth}
F.~A. Graybill and R.~D. Morrison.
\newblock Sample size for a specified width confidence interval on the variance of a normal distribution.
\newblock {\em Biometrics}, 16:636--641, 1960.

\bibitem{Macaulay2}
D.~R. Grayson and M.~E. Stillman.
\newblock Macaulay2, a software system for research in algebraic geometry.
\newblock Available at \url{http://www2.macaulay2.com}.

\bibitem{GrossPetrovikGoodnessFitMarkovBasis}
E.~Gross, S.~Petrovi\'c, and D.~Stasi.
\newblock Goodness of fit for log-linear network models: dynamic {M}arkov bases using hypergraphs.
\newblock {\em Ann. Inst. Statist. Math.}, 69(3):673--704, 2017.

\bibitem{GrossSullivantMLETbhreshold}
E.~Gross and S.~Sullivant.
\newblock The maximum likelihood threshold of a graph.
\newblock {\em Bernoulli}, 24(1):386--407, 2018.

\bibitem{TomonariHierarchicalSubspace}
H.~Hara, T.~Sei, and A.~Takemura.
\newblock Hierarchical subspace models for contingency tables.
\newblock {\em J. Multivariate Anal.}, 103:19--34, 2012.

\bibitem{HirjiExactAnalysis}
K.~F. Hirji.
\newblock {\em Exact analysis of discrete data}.
\newblock Chapman \& Hall/CRC, Boca Raton, FL, 2006.

\bibitem{PeterHoffPvalueIndirectInformation}
P.~Hoff.
\newblock Smaller {$p$}-values via indirect information.
\newblock {\em J. Amer. Statist. Assoc.}, 117(539):1254--1269, 2022.

\bibitem{MarkovBasesHierarcSullivantHosten}
S.~Hosten and S.~Sullivant.
\newblock A finiteness theorem for {M}arkov bases of hierarchical models.
\newblock {\em J. Combin. Theory Ser. A}, 114(2):311--321, 2007.

\bibitem{ExactConditionalTestsKreiner}
S.~Kreiner.
\newblock Analysis of multidimensional contingency tables by exact conditional tests: techniques and strategies.
\newblock {\em Scand. J. Statist.}, 14(2):97--112, 1987.

\bibitem{TSPLehmannRomano}
E.~L. Lehmann and J.~P. Romano.
\newblock {\em Testing statistical hypotheses}.
\newblock Springer Texts in Statistics. Springer, Cham, fourth edition, 2021.

\bibitem{LehmannScheffeCompleteness1}
E.~L. Lehmann and H.~Scheff\'{e}.
\newblock Completeness, similar regions, and unbiased estimation. {I}.
\newblock {\em Sankhy\={a}}, 10:305--340, 1950.

\bibitem{MillerSturmfelsCombinatorial}
E.~Miller and B.~Sturmfels.
\newblock {\em Combinatorial commutative algebra}, volume 227 of {\em Graduate Texts in Mathematics}.
\newblock Springer-Verlag, New York, 2005.

\bibitem{MotzkinPolynomial}
T.~S. Motzkin.
\newblock The arithmetic-geometric inequality.
\newblock In {\em Inequalities ({P}roc. {S}ympos. {W}right-{P}atterson {A}ir {F}orce {B}ase, {O}hio, 1965)}, pages 205--224. Academic Press, New York-London, 1967.

\bibitem{ExistenceofUMPSTests}
D.~Paindaveine.
\newblock On {UMPS} hypothesis testing.
\newblock {\em Ann. Inst. Statist. Math.}, 76(2):289--312, 2024.

\bibitem{PowersCertificatesofPositivity}
V.~Powers.
\newblock {\em Certificates of positivity for real polynomials---theory, practice, and applications}, volume~69 of {\em Developments in Mathematics}.
\newblock Springer, Cham, [2021] \copyright 2021.

\bibitem{SenGuptaExponentialUMPU}
A.~SenGupta.
\newblock Optimal tests in multivariate exponential distributions.
\newblock In {\em The exponential distribution}, pages 351--376. Gordon and Breach, Amsterdam, 1995.

\bibitem{SampleSizePoissonLogistic}
G.~Shieh.
\newblock Sample size calculations for logistic and {P}oisson regression models.
\newblock {\em Biometrika}, 88(4):1193--1199, 2001.

\bibitem{signorini1991samplePoisson}
D.~F. Signorini.
\newblock Sample size for poisson regression.
\newblock {\em Biometrika}, 78(2):446--450, 1991.

\bibitem{smith1996montecarloexact}
P.~W. Smith, J.~J. Forster, and J.~W. McDonald.
\newblock Monte carlo exact tests for square contingency tables.
\newblock {\em Journal of the Royal Statistical Society Series A: Statistics in Society}, 159(2):309--321, 1996.

\bibitem{NonExistenceofUnbiasedEstimators}
A.~Somekh-Baruch, A.~Leshem, and V.~Saligrama.
\newblock On the non-existence of unbiased estimators in constrained estimation problems.
\newblock {\em IEEE Trans. Inform. Theory}, 64(8):5549--5554, 2018.

\bibitem{NilsUStatisticSingularity}
N.~Sturma, M.~Drton, and D.~Leung.
\newblock Testing many constraints in possibly irregular models using incomplete {$U$}-statistics.
\newblock {\em J. R. Stat. Soc. Ser. B. Stat. Methodol.}, 86(4):987--1012, 2024.

\bibitem{SullivantAlgStats}
S.~Sullivant.
\newblock {\em Algebraic statistics}, volume 194 of {\em Graduate Studies in Mathematics}.
\newblock American Mathematical Society, Providence, RI, 2018.

\bibitem{UhlerMLEGaussianGraph}
C.~Uhler.
\newblock Geometry of maximum likelihood estimation in {G}aussian graphical models.
\newblock {\em Ann. Statist.}, 40(1):238--261, 2012.

\bibitem{WatanabeBook}
S.~Watanabe.
\newblock {\em Algebraic geometry and statistical learning theory}, volume~25 of {\em Cambridge Monographs on Applied and Computational Mathematics}.
\newblock Cambridge University Press, Cambridge, 2009.

\bibitem{ZieglerLectureonPolytopes}
G.~M. Ziegler.
\newblock {\em Lectures on polytopes}, volume 152 of {\em Graduate Texts in Mathematics}.
\newblock Springer-Verlag, New York, 1995.

\end{thebibliography}

\newpage

\appendix 
\section*{\huge{Supplementary Material: Proofs}}

%%%%%%%%%%%%%%%%%%%%%%%%%%%%%%%%%%%%%%%%%%%%%
%%%%%%%%%%%%%%%%%%%%%%%%%%%%%%%%%%%%%%%%%%%%%

\begin{customlem}{1}
\label{custlem:ExistenceLevel}
    The four different sets of level-$\alpha$, unbiased, non-trivial unbiased, and strictly unbiased tests are all convex sets of test functions. If there exists a size-$\alpha_0 \in [0,1)$ NTUB or SUB test, then for any $\alpha \in (0,1)$ there exists a size-$\alpha$ NTUB or SUB test respectively.
\end{customlem}
\begin{proof}
    Convexity for level-$\alpha$, UB, and SUB tests follows easily by the linearity of the power function $\beta_{\lambda \phi_1 + (1-\lambda)\phi_2} = \lambda \beta_{\phi_1} + (1-\lambda)\beta_{\phi_2}$. Suppose $\phi_1,\phi_2$ are NTUB tests. There necessarily exist points $\theta_0 \in \Theta_0$ and $\theta_A$ where $\beta_{\phi_1}(\theta_0) < \beta_{\phi_1}(\theta_A)$ as otherwise $\beta_{\phi_1}$ would be constant. It follows that $ \lambda \beta_{\phi_1}(\theta_0) + (1-\lambda)\beta_{\phi_2}(\theta_0) <  \lambda \beta_{\phi_1}(\theta_A) + (1-\lambda)\beta_{\phi_2}(\theta_A)$, implying that $\beta_{\lambda \phi_1 + (1-\lambda)\phi_2}$ is not constant for $\lambda \in [0,1]$.
    
    To prove the last statement, we will take a convex combination of a size-$\alpha_0$ NTUB or SUB test $\phi$ with either the constant test $1$ or the constant test $0$. If $\alpha > \alpha_0$ consider the test $\lambda \phi + (1-\lambda)$ that is NTUB or SUB when $\phi$ is, by the proof above. The size of this test is $\text{Size}(\lambda \phi + (1-\lambda)) = \lambda\text{Size}(\phi) + (1-\lambda) = \lambda \alpha_0 +(1-\lambda)$. Taking $\lambda = \tfrac{1-\alpha}{1-\alpha_0} \in (0,1)$ gives the desired test. The case where $\alpha < \alpha_0$ is handled similarly by considering $\tfrac{\alpha}{\alpha_0}\phi$.  
\end{proof}

%%%%%%%%%%%%%%%%%%%%%%%%%%%%%%%%%%%%%%%%%%%%%%%%%
%%%%%%%%%%%%%%%%%%%%%%%%%%%%%%%%%%%%%%%%%%%%%%%%%

\begin{customlem}{2}
    If there exists a level-$\alpha_0$ UMPU test of $\Po$ against $\Pa$ for some $\alpha_0 \in [0,1)$ then there also exists a level-$\alpha$ UMPU test for all $\alpha \in [0,\alpha_0]$.  
\end{customlem}
\begin{proof}
    Let $\beta_\phi$ be the power function of a size-$\alpha_0$ UMPU test. For $\alpha \leq \alpha_0$ and $c = \alpha/\alpha_0 \leq 1$ the test $c\phi$ is a size-$\alpha$ test with power function $c\beta_\phi$. Assume by contradiction that this test was not UMPU and so $\beta_{\tilde{\phi}}(\theta_A) > c\beta_\phi(\theta_A)$ for some NTUB, size-$\alpha$ test $\tilde{\phi}$. Then $(1- c) \beta_\phi(\theta_A) +  \beta_{\tilde{\phi}}(\theta_A) > \beta_\phi(\theta_A)$ and $(1-c)\phi + \tilde{\phi}$ is a level-$\alpha$ NTUB test, contradicting the assumption that $\phi$ is UMPU.
\end{proof}

%%%%%%%%%%%%%%%%%%%%%%%%%%%%%%%%%%%%%%%%%%%%%%%%
%%%%%%%%%%%%%%%%%%%%%%%%%%%%%%%%%%%%%%%%%%%%%%%%

\begin{customthm}{1}
    \label{custthm:MainExistenceUniqueness}
Assume that the null hypothesis set is closed.  There exists an NTUB test of $\Po$ against $\Pa = \Delta_{k-1}\backslash \Po$ when the sample size is $n$ if and only if there exists a polynomial $\Tilde{\beta}$ of degree less than or equal to $n$ that is non-constant on $\Delta_{k-1}$ and has sub-level sets with
    \begin{align}
    \label{custeqn:NTUBExistence}
        \Po \subseteq \{ \bs{\pi} \in \Delta_{k-1}: \Tilde{\beta}(\bs{\pi}) \leq 0\} \; \text{ and } \; \Pa \subseteq \{ \bs{\pi} \in \Delta_{k-1}: \Tilde{\beta}(\bs{\pi}) \geq 0\}.
    \end{align}
    There exists an SUB test if and only if there exists a similarly defined polynomial with sub-level sets of the form
    \begin{align}
    \label{custeqn:SUBExistence}
             \Po = \{ \bs{\pi} \in \Delta_{k-1}: \Tilde{\beta}(\bs{\pi}) \leq 0\} \; \text{ and } \; \Pa = \{ \bs{\pi} \in \Delta_{k-1}: \Tilde{\beta}(\bs{\pi}) > 0\}.
    \end{align}
\end{customthm}
\begin{proof}
Suppose that there exists a NTUB or SUB unbiased test $\phi$ with size-$\alpha$. If $\phi$ is NTUB the translated power polynomial $\Tilde{\beta} = \beta_\phi - \alpha$ is non-constant on $\Delta_{k-1}$ and satisfies \eqref{custeqn:NTUBExistence} because $\sup_{\bs{\pi} \in \Po} \Tilde{\beta}(\bs{\pi}) = 0$ by the size assumption, and $0 = \sup_{\bs{\pi} \in \Po} \Tilde{\beta}(\bs{\pi}) \leq \Tilde{\beta}(\bs{\pi}_A)$ for all $\bs{\pi}_A \in \Pa$ by the unbiasedness assumption. If $\phi$ is SUB then as $\Po$ is assumed to be closed there exists a $\bs{\pi}_0 \in \Po$ with that $\beta_\phi(\bs{\pi}_0) - \alpha = \Tilde{\beta}(\bs{\pi}_0) = 0$.  By the SUB assumption $0 = \Tilde{\beta}(\bs{\pi}_0) < \Tilde{\beta}(\bs{\pi}_A)$ for all $\bs{\pi}_A \in \Pa$, proving condition \eqref{custeqn:SUBExistence}.

The converse statements are shown by taking a polynomial $\Tilde{\beta}$ that satisfies the sub-level set assumptions and normalizing it so that it is equal to some power polynomial corresponding to either a NTUB or SUB test. Without loss of generality $\Tilde{\beta}$ can be assumed to be homogeneous as homogenization does not alter any of the assumptions on $\Tilde{\beta}$ outlined in the theorem. There exists a large enough constant $b \geq 0$ such that the coefficient of the monomial $\bs{\pi}^{\bl{x}}$ in the polynomial $\Tilde{\beta} + b(\sum_{i = 1}^k \pi_i)^n$ is non-negative because the coefficient of $\bs{\pi}^{\bl{x}}$ in $b(\sum_{i = 1}^k \pi_i)^n$ is $b{n \choose \bl{x}} \geq 0$. There also exists a small enough constant $a > 0$ such that the coefficient of the monomial $\bs{\pi}^{\bl{x}}$ in the polynomial  $\beta^* \coloneqq a\big(\Tilde{\beta} + b(\sum_{i = 1}^k \pi_i)^n\big)$ is bounded above by ${n \choose \bl{x}}$. By Lemma \ref{lem:PowerFunctionCharacterization} the polynomial $\beta^*$ is a power polynomial. It is non-constant on $\Delta_{k-1}$ if $\Tilde{\beta}$ is. The assumptions \eqref{custeqn:NTUBExistence} and \eqref{custeqn:SUBExistence} respectively imply the following two conditions on $\beta^*$:
\begin{align*}
      \Po & \subseteq \{ \bs{\pi} \in \Delta_{k-1}: \beta^*(\bs{\pi}) \leq ab\} \; \text{ and } \; \Pa \subseteq \{ \bs{\pi} \in \Delta_{k-1}: \beta^*(\bs{\pi}) \geq ab\}, 
      \\
                  \Po & = \{ \bs{\pi} \in \Delta_{k-1}: \beta^*(\bs{\pi}) \leq ab\} \; \text{ and } \; \Pa = \{ \bs{\pi} \in \Delta_{k-1}: \beta^*(\bs{\pi}) > ab\}.
\end{align*}
These two conditions respectively imply \eqref{eqn:UBDefinition} and strict inequality version of \eqref{eqn:UBDefinition} for the power polynomial $\beta^*$.
\end{proof}

%%%%%%%%%%%%%%%%%%%%%%%%%%%%%%%%%%%%%%%%%%%%%%%%%
%%%%%%%%%%%%%%%%%%%%%%%%%%%%%%%%%%%%%%%%%%%%%%%%%

\begin{customlem}{3}[Similarity on the boundary]
\label{custlem:SimilarityonBd}
 Assume that a NTUB or SUB test for $\Po$ exists and $\Tilde{\beta}$ is given as in \eqref{eqn:NTUBExistence} or \eqref{eqn:SUBExistence}. Then $\Tilde{\beta}(\bs{\pi}) = 0$ for every point $\bs{\pi}$ in the boundary of $\Po$ relative to $\Delta_{k-1}$. 
\end{customlem}
\begin{proof}
    By the definition of $\Tilde{\beta}$, it is known that $\Tilde{\beta}(\bs{\pi}) \leq 0$ for $\bs{\pi} \in \partial \Po$. As $\bs{\pi} \in \partial \Po$ there exists a sequence $\{\bs{\pi}_i \}_{i = 1}^\infty \subseteq \Pa$ with $\bs{\pi}_i \rightarrow \bs{\pi}$. Note that it is important that $\bs{\pi}$ is in the boundary of $\Po$ relative to $\Delta_{k-1}$ for this to hold. By continuity $0 \leq \Tilde{\beta}(\bs{\pi}_i) \rightarrow \Tilde{\beta}(\bs{\pi}) \leq 0$, as needed. 
\end{proof}

%%%%%%%%%%%%%%%%%%%%%%%%%%%%%%%%%%%%
%%%%%%%%%%%%%%%%%%%%%%%%%%%%%%%%%%%%

\begin{customthm}{2}[Existence of UB Tests for Polytope Null Hypotheses]
Parameterize the probability simplex $(\pi_{1},\ldots ,\pi_{k}) \in \Delta_{k-1}$ by the projection onto its first $k-1$ coordinates $\Tilde{\Delta}_{k-1} = \{\bs{\pi}: \sum_{i = 1}^{k-1}\pi_i \leq 1, \; \pi_i \geq 0, \; i = 1,\ldots,k-1\} \subseteq \mb R^k$. Let a null hypothesis set be given by the $(k-1)$-dimensional polytope $\mc P_0 = \{\bs{\pi} \in \Tilde{\Delta}_{k-1}: \bl{a}_i^\intercal \bs{\pi} - b_i \geq 0, \; i = 1,\ldots,m\}$, where it is assumed that the facet of the polytope corresponding to the affine hyperplane $H_i = \{\bs{\pi}: \bl{a}_i^\intercal \bs{\pi} = b_i \}$ is $(k-2)$-dimensional and $H_i$ intersects $\mathrm{int}(\Tilde{\Delta}_{k-1})$.
There exists an NTUB or SUB test of $\Po$ if and only if for every $i \neq j$ the set $H_i \cap H_j \cap \Po$ is contained in the boundary of $\Tilde{\Delta}_{k-1}$.       
\end{customthm}
\begin{proof}
    Assume that $H_i \cap H_j \cap \Po$ is contained in the boundary of $\Tilde{\Delta}_{k-1}$ for every $i \neq j$. The polynomial $\Tilde{\beta}(\bs{\pi}) = -\prod_{i = 1}^m (\bl{a}_i^\intercal \bs{\pi} - b_i)$ is non-zero and non-positive on $\mc{P}_0$. We wish to show that $\Tilde{\beta}$ is positive on $\Pa$, and thus after normalization gives a power polynomial for a SUB test. For $\bs{\pi} \in \mc P_A$ at least one of the constraints $\bl{a}_i^\intercal \bs{\pi} - b_i \geq 0$ does not hold, say $\bl{a}_1^\intercal \bs{\pi} - b_1 < 0$. If $\bl{a}_i^\intercal \bs{\pi} - b_i > 0$ for $i = 2,\ldots, m$ then $\Tilde{\beta}(\bs{\pi}) > 0$, so by contradiction we assume that $\bl{a}_i^\intercal \bs{\pi} - b_i \leq 0$ for some $i \in \{2,\ldots,l\}$. Let $\bl{p}_1 \in H_1 \cap \Po \cap \text{int}(\Tilde{\Delta}_{k-1})$, implying that $\bl{a}_1^\intercal \bl{p}_1 - b_1 = 0$ and $\bl{a}_i^\intercal \bl{p}_1 - b_i > 0$ for $i > 1$. The latter inequality follows since $H_1 \cap H_i \cap \Po$ is contained in the boundary of $\Tilde{\Delta}_{k-1}$. We consider the convex combination $\bs{\pi}_\lambda = (1-\lambda)\bs{\pi} + \lambda \bl{p}_1$ which is in $\text{int}(\Tilde{\Delta}_{k-1})$ for all $\lambda \in (0,1)$. Choose the smallest $\lambda_0 \in (0,1)$ such that $\bl{a}_i^\intercal \bs{\pi}_{\lambda_0} - b_i \geq 0$ for all $i = 2,\ldots,m$ and without loss of generality assume that $\bl{a}_2^\intercal \bs{\pi}_{\lambda_0} - b_2 = 0$. Note that it remains the case that $\bl{a}_1^\intercal \bs{\pi}_{\lambda_0} - b_1 < 0$. The next step is to similarly choose a point $\bl{p}_2 \in H_2 \cap \Po \cap \text{int}(\Tilde{\Delta}_{k-1})$ and consider the convex combination $\bs{\pi}_\delta' = (1-\delta)\bs{\pi}_{\lambda_0} + \delta \bl{p}_2$ which also lies in $\text{int}(\Tilde{\Delta}_{k-1})$. For all $\delta$ we have
    \begin{align*}
        \bl{a}_2^\intercal \bs{\pi}'_\delta  - b_2 = (1-\delta)(\bl{a}_2^\intercal \bs{\pi}_{\lambda_0} - b_2) + \delta(\bl{a}_2^\intercal \bl{p}_2 - b_2) = 0.
    \end{align*}
    As $\bl{a}_1^\intercal \bl{p}_2 - b_2 > 0$, we choose the smallest $\delta_0 \in (0,1)$ such that $\bl{a}_1^\intercal \bs{\pi}'_{\delta_0} - b_1 = 0$. The point $\bs{\pi}'_{\delta_0}$ is in $H_1 \cap H_2 \cap \Po \cap \text{int}(\Tilde{\Delta}_{k-1})$, contradicting the assumed structure of the hyperplanes. It follows that $\Tilde{\beta}$ is a polynomial associated with a strictly unbiased test of $\Po$. 

    To prove the converse we assume that there is a non-zero polynomial $\Tilde{\beta}$ defined on $\Tilde{\Delta}_{k-1}$ that is non-positive on $\Po$ and non-negative on $\Pa$. It will eventually be shown that if there exists a $\bs{\pi}_0 \in H_i \cap H_j \cap \Po \cap \text{int}(\tilde{\Delta}_{k-1})$ this leads to a contradiction. As $\Tilde{\beta}$ must vanish on $\partial \Po$, which has full dimensional faces, it must be the case that $\bl{a}_i^\intercal \bs{\pi} - b_i$ divides $\Tilde{\beta}$. Thus $\Tilde{\beta}(\bs{\pi}) = g(\bs{\pi})\prod_{i = m}^l (\bl{a}_i^\intercal \bs{\pi} - b_i)^{c_i}$ for some integers $c_i$, where $\bl{a}_i^\intercal \bs{\pi} - b_i$, $i = 1,\ldots,m$ does not divide $g(\bs{\pi})$. 
    
    We begin by showing that every $c_i$ must be odd. Suppose by contradiction that $c_1$ is even and pick a $\bs{\pi}_1 \in \text{relint}(\Po \cap H_1) \cap \text{int} (\Tilde{\Delta}_{k-1})$. Such a point exists by the assumption that every facet of $\Po$ is full-dimensional and intersects the interior of $\Tilde{\Delta}_{k-1}$. In particular, $\bs{\pi}_1 \notin \bigcup_{i = 2}^m H_i$ and thus $\bl{a}_1^\intercal \bs{\pi}_1 - b_1 = 0$ and $\bl{a}_i^\intercal \bs{\pi}_1 - b_i > 0$, $i > 1$. The goal is to show that $g(\bs{\pi}_1) = 0$ as this will show that $g$ vanishes on the relatively open set $\text{relint}(\Po \cap H_1) \cap \text{int} (\Tilde{\Delta}_{k-1})$ contained in $H_1$. As a consequence, $g$ vanishes on $H_1$, leading to the contradiction that $(\bl{a}_1^\intercal \bs{\pi} - b_1)$ divides $g$. 
    
    Choose a $\bs{\tau} \in \text{int}(\Po)$, which exists because $\Po$ is assumed to be a $k$-dimensional polytope.
Define $\bs{\pi}_\lambda = \lambda\bs{\pi}_1 + (1-\lambda)\bs{\tau}$. There exists a small enough $\epsilon$ such that  $\bs{\pi}_{\lambda} \in \Pa$ and $\bl{a}_i^\intercal \bs{\pi}_\lambda - b_i > 0$ for $i > 1$ when $\lambda \in (1, 1 + \epsilon]$. As $c_1$ is even, the sign of $\prod_{i = 1}^l (\bl{a}_i^\intercal \bs{\pi} - b_i)^{c_i}$ does not change as $\lambda$ passes from the interval $[0,1]$ to $(1,1+\epsilon]$. However, as $\lambda$ increases past one, $\bs{\pi}_{\lambda}$ moves from $\Po$ into $\Pa$ so the sign of $\Tilde{\beta}(\bs{\pi}_\lambda)$ must either change at $\lambda = 1$ or alternatively $\Tilde{\beta}(\bs{\pi}_\lambda) = 0$ in a neighbourhood of $\lambda = 1$. In the former case $g(\bs{\pi}_\lambda)$ changes sign at $\lambda = 1$ so $g(\bs{\pi}_1) = 0$. In the latter case $\Tilde{\beta}(\bs
{\pi}_\lambda)$ equals zero on an open interval of $\lambda$ and thus $\Tilde{\beta}(\bs{\pi}_\lambda) = 0$ for all $\lambda \in \mb{R}$. This further implies that $g(\bs{\pi}_\lambda) = 0$ for all $\lambda \in \mb{R}$ and in particular, $g(\bs{\pi}_1) = 0$. 

% The above argument remains valid if $\bs{\tau}$ is replaced with $\Tilde{\bs{\tau}} \in U$ where $U$ is an open ball contained in $\text{int}(\Po)$. That is, define $\bs{\pi}_{\lambda,\Tilde{\bs{\tau}}} = \lambda \bs{\pi}_1 + (1-\lambda)\Tilde{\bs{\tau}}$ for $\Tilde{\bs{\tau}} \in U$. If for any $\Tilde{\bs{\tau}} \in U$, $\Tilde{\beta}(\bs{\pi}_{\lambda,\Tilde{\bs{\tau}}})$ changes sign at $\lambda = 1$, $g(\bs{\pi}_1) = 0$ and we are done. The alternative case is that for all  $\Tilde{\bs{\tau}} \in U$, $\Tilde{\beta}(\bs{\pi}_{\lambda,\Tilde{\bs{\tau}}})$ vanishes on an open interval of $\lambda = 1$.
% We conclude that for $\Tilde{\bs{\tau}} \in U$ and $\lambda \in \mb{R}$ $\Tilde{\beta}(\bs{\pi}_{\lambda,\Tilde{\bs{\tau}}}) = 0$. As the image of $\bs{\pi}_{\lambda,\Tilde{\bs{\tau}}}$ over $\Tilde{\bs{\tau}} \in U$ and $\lambda \in \mb{R}$ contains the open set $U$, we obtain that $\Tilde{\beta} = 0$, a contradiction. 

Next we assume that $\Tilde{\beta}(\bs{\pi}) = g(\bs{\pi})\prod_{i = 1}^m (\bl{a}_i^\intercal \bs{\pi} - b_i)^{c_i}$ has every $c_i$ odd. It at this point in the proof that we invoke the assumption that there exists a $\bs{\pi}_0 \in \Po \cap \big(\cap_{i = 1}^l H_i\big) \backslash \big(\cap_{i = l+1}^m H_i \big)$, $l \geq 2$ where $\bs{\pi}_0$ is also in $\text{int}(\Tilde{\Delta}_{k-1})$. 
% We show that a $\Tilde{\beta}$ associated with a NTUB test cannot exist by utilizing the assumption that there exists a $\bs{\pi}_0 \in \Po \cap \big(\cap_{i = 1}^l H_i\big) \backslash \big(\cap_{i = l+1}^m H_i \big)$, $l \geq 2$ where $\bs{\pi}_0$ is also in $\text{int}(\Tilde{\Delta}_{k-1})$. 
We work in a small enough open neighbourhood $U \subset \text{int}(\Tilde{\Delta}_{k-1})$ containing $\bs{\pi}_0$ that has an empty intersection with $\cap_{i = l+1}^m H_i$, and thus $\prod_{i = l+1}^l (\bl{a}_i^\intercal \bs{\pi} - b_i)^{c_i}$ does not change sign in $U$. The hyperplanes $H_1,\ldots,H_l$ divide $U$ into different open, convex sets where the signs of $(\bs{a}_1^\intercal \bs{\pi}-b_1,\ldots,\bl{a}_l^\intercal \bs{\pi} - b_l)$ are constant within each set. One of these regions, the region where the signs of these linear polynomials are all positive, is contained in $\Po$. As $l \geq 2$ there exists at least three different regions, as if there were only two, a single hyperplane would divide $U$ into two regions, implying that $l = 1$. Only a single region is contained in $\Po$, with the remaining regions, of which there are at least two, being contained in $\Pa$. There exist two such regions $R_1,R_2 \subset \Pa$ that are neighbouring in the sense that the sign of $(\bs{a}_1^\intercal \bs{\pi}-b_1,\ldots,\bl{a}_l^\intercal \bs{\pi} - b_l)$ differs only in a single coordinate in $R_1$ as compared to $R_2$. Without loss of generality assume that the sign of $\bs{a}_1^\intercal \bs{\pi} - a_1$ is positive in $R_1$ and negative in $R_2$, meaning that the hyperplane $H_1$ is the boundary between $R_1$ and $R_2$. The goal will be to show that $g(\bs{\pi})$ must vanish on $H_1$, which yields the contradiction that $\bl{a}_1^\intercal \bs{\pi} - b_1$ divides $g$. A similar proof to the proof used to show that each $c_i$ is even will be used. Let $V \subset R_1$ be open and take $\bs{p} \in \text{int}(R_2)$. For any $\bs{\tau} \in V$ define the line segment $\bl{p}_{\lambda,\bs{\tau}} = (1-\lambda)\bs{\tau} + \lambda\bl{p}$ connecting $\bs{\tau}$ to $\bl{p}$. Define $\bl{p}_{\lambda_{\bs{\tau}} ,\bs{\tau}}$ to be the unique point along this segment that lies in $H_1$. The product $\prod_{i = 1}^m (\bl{a}_i^\intercal \bs{\pi} - b_i)^{c_i}$ changes sign at $\bl{p}_{\lambda_{\bs{\tau}} ,\bs{\tau}}$ along the segment  $\bl{p}_{\lambda,\bs{\tau}}$ because $\bl{a}_1^\intercal \bs{\pi} - b_1$ changes sign in $H_1$ while every $\bl{a}_i^\intercal \bs{\pi} - b_i$, $i > 1$ has the same sign in $R_1$ and $R_2$ by construction. As $R_1,R_2 \subset \Pa$ the function $\Tilde{\beta}$ cannot have two different signs in the region $R_1 \cup R_2$. It follows that $g(\bl{p}_{\lambda,\bs{\tau}})$ must change sign at $\lambda = \lambda_{\bs{\tau}}$ so $g(\bl{p}_{\lambda_{\bs{\tau}},\bs{\tau}}) = 0$. As $\bs{\tau}$ varies over $V$ the set of points $\bl{p}_{\lambda_{\bs
\tau},\bs{\tau}}$ includes a relatively open set in $H_1$. The function $g$ vanishes on a relatively open set of the hyperplane $H_1$, which gives the desired contradiction that $\bl{a}_1^\intercal \bs{\pi} - b_1$ divides $g$. 
 \end{proof}

%%%%%%%%%%%%%%%%%%%%%%%%%%%%%%%%%%%%%%%%
%%%%%%%%%%%%%%%%%%%%%%%%%%%%%%%%%%%%%%%%
%%%%%%%%%%%%%%%%%%%%%%%%%%%%%%%%%%%%%%%%

\begin{customthm}{3}
  When $k \geq 3$ there is an NTUB test of the log-linear null hypothesis 
\begin{align}
\label{custeqn:LinearHypothesisMultinomial}
    \Po = \bigg\{ \bs{\pi} \in  \text{int}(\Delta_{k-1}): \sum_{i = 1}^{k-1} a_i \log\big(\tfrac{\pi_i}{\pi_k}\big) = b \bigg\}. 
\end{align}
if and only if there exists a representation of $\Po$ where every $a_i$, $i = 1,\ldots,k-1$ is a rational number. 
\end{customthm}
\begin{proof}
      Assume first that $(a_1\ldots,a_{k-1})^\intercal = \bl{a} \in \mb{Q}^{k-1}$. By clearing denominators, without loss of generality we assume that $\bl{a} \in \mb{Z}^{k-1}$. A $\bs{\pi} \in \text{int}(\Delta_{k-1})$ 
satisfies \eqref{custeqn:LinearHypothesisMultinomial} if and only if
\begin{align}
\label{custeqn:PolynomialLogLinEquations}
    \bigg(\prod_{j = 1}^{k-1} \pi_j^{a_{j}}\bigg) - e^{b}\pi_k^{- \sum_{j = 1}^{k-1}a_{j}} = 0.
\end{align}
 Clearing denominators in this rational expression we obtain a polynomial equation of the form $\bs{\pi}^I - e^b\bs{\pi}^J = 0$ with $I,J$ disjoint, where every point in $\Po$ satisfies this equation. The polynomial $(\bs{\pi}^I - e^b\bs{\pi}^J)^2$ is an NTUB separating polynomial as long as it does not identically vanish on $\Delta_{k-1}$. At least one of $\bs{\pi}^I$ or $\bs{\pi}^J$ is not equal to one, suppose it is $\bs{\pi}^I$. Evaluating the polynomial $(\bs{\pi}^I - e^b\bs{\pi}^J)^2$ at the average of the standard basis vectors in $J$, $\frac{1}{\vert J \vert} \sum_{j \in J} \bl{e}_j \in \Delta_{k-1}$ produces a non-zero scalar.

To prove the converse, we invoke Theorem 4.4.1 of \citep{TSPLehmannRomano} and show that this produces the trivial test. We assume without loss of generality that $a_1 = 1$ and that $a_2$ is irrational. This can always be assumed by relabeling the indices of the $\pi_i$s and by dividing the log-linear equation by $a_1$. Defining $\theta_i = \log\big(\frac{\pi_1}{\pi_k}\big)$ a multinomial distribution can be parameterized as 
\begin{align*}
  p(\bl{x}|\bs{\theta}) & =   {n \choose x_1,\ldots,x_k} \exp\bigg(\sum_{i = 1}^{k-1} x_i \theta_i - n\log(\pi_k) \bigg)
  \\
  & = {n \choose x_1,\ldots,x_k}e^{x_1 b(k-1)} \exp\bigg( x_1 \big(\sum_{i = 1}^{k-1}a_i\theta_i -b\big) + \sum_{i = 2}^{k-1}(x_i - a_ix_1) \theta_i - n\log(\pi_k) \bigg). 
\end{align*}
Define $\bl{y} = (x_2 - a_2 x_1,\ldots, x_{k-1} - a_{k-1}x_{k-1})$. 
The UMPU test of the linear hypothesis $\sum_{i = 1}^{k-1} a_i\theta_i - b = 0$ has the form by equation (4.16) of \citep{TSPLehmannRomano}: 
\begin{align*}
    \phi(x_1, \bl{y}) = \begin{cases}
        1 \;\;\; & x_1 < C_1(\bl{y}) \text{ or } x_1 > C_2(\bl{y})
        \\
        \gamma_i(\bl{y}) \;\;\; &x_1 = C_i(\bl{y}), \; i = 1,2
        \\
        0 \;\;\; & C_1(\bl{y}) < x_1 < C_2(\bl{y}),
    \end{cases}
\end{align*}
with the $C_i(\cdot)$ and $\gamma_i(\cdot)$ determined by the equations
\begin{align*}
    E(\phi(x_1,\bl{y})| \bl{y}) = \alpha, \;\; E(x_1\phi(x_1,\bl{y})| \bl{y}) = \alpha E(x_1| \bl{y}),
\end{align*}
where the expectations in these equations must hold for every distribution in $\Po$. The key observation is that if the value of $x_2 - a_2 x_1$ is known, then the value of $x_1$ is also known because $a_2$ is irrational and $x_1,x_2 \in \mb{N}$. The conditional distribution of $x_1$ given $\bl{y}$ is thus a degenerate point mass at $x_1$, implying that $C_1(\bl{y}) = C_2(\bl{y}) = x_1$ and $\gamma_i(\bl{y}) = \alpha$. The resulting test is the trivial test. As it is the UMPU test, there does not exist an NTUB test for this hypothesis.  
\end{proof}

%%%%%%%%%%%%%%%%%%%%%%%%%%%%%%%%%%%%%%%%
%%%%%%%%%%%%%%%%%%%%%%%%%%%%%%%%%%%%%%%%

\begin{customcor}{1}
\label{custcor:ExistenceforAlgebraic}
        There exists an unbiased test for the algebraic null hypothesis \eqref{eqn:AlgebraicNull} when the sample size is $n$ if and only if there exists a degree $n$, non-zero polynomial $\Tilde{\beta} \in I_{\mb{R}^k}(\Po$) where 
    \begin{align}
    \label{custeqn:AlgNullNTUB}
        \Pa \subseteq \{ \bs{\pi} \in \Delta_{k-1}: \Tilde{\beta}(\bs{\pi}) \geq 0\},
    \end{align}
    for the test to be a NTUB test, and
    \begin{align}
    \label{custeqn:AlgNullSUB}
                \Pa = \{ \bs{\pi} \in \Delta_{k-1}: \Tilde{\beta}(\bs{\pi}) >  0\},
    \end{align}
    for the test to be a SUB test. 
\end{customcor}
\begin{proof}
    This result follows immediately from Lemma \ref{lem:SimilarityonBd} if it is shown that  $\partial \Po = \Po$. Let $\Po = \{\bs{\pi} \in \Delta_{k-1}: f_{(i)}(\bs{\pi}) = 0, \;i = 1,\ldots,m\} \subsetneq \Delta_{k-1}$. Suppose by contradiction that $\bs{\pi}_0 \in \Po \backslash \partial \Po$, so that there exists a relatively open ball $U \subset \text{aff}(\Delta_{k-1})$ with $\bs{\pi}_0 \in U \subseteq \Po$. Then $f_{(i)}\vert_U = 0$ for every $i$ and in particular there exists an open ball $V \subset \mb{R}^{k-1}$ with $f_{(i)}(\pi_1,\pi_2,\ldots,1 - \sum_{i = 1}^{k-1}\pi_i) = 0$ for all $(\pi_1,\ldots,\pi_{k-1}) \in V$. Any polynomial that vanishes on an open ball must be the zero polynomial so $f_{(i)}(\pi_1,\pi_2,\ldots,1 - \sum_{i = 1}^{k-1}\pi_i) = 0$ over $\mb{R}^{k-1}$ for every $i$, implying the contradiction that $\Po = \Delta_{k-1}$.
\end{proof}

%%%%%%%%%%%%%%%%%%%%%%%%%%%%%%%%%%%%%%%%%%%%%%%%%%%
%%%%%%%%%%%%%%%%%%%%%%%%%%%%%%%%%%%%%%%%%%%%%%%%%%%

\begin{customthm}{4}
\label{custthm:SOSUBThresholdUpperBound}
    Assume that $\Po$ is an algebraic null hypothesis and let $I_{\mb{R}^{k-1}}(\Po)$ have the reduced Gr\"obner basis $\mc{G} = \{g_{(1)},\ldots, g_{(m)}\}$ with respect to a graded monomial order.  The minimum degree of a non-zero polynomial in $I_{\mb{R}^{k-1}}(\Po)$ is equal to $\min_{i}\text{deg}(g_{(i)})$. The lowest degree of a SOS polynomial $\Tilde{\beta}$ corresponding to an NTUB test is $2\min_{i}\text{deg}(g_{(i)})$.  Let $\mc{G}_i$ consist of the subset of polynomials in $\mc{G}$ with degree less than or equal to $i$. Define $d = \min \{i:  V_{\mb{R}^{k-1}}(\mc{G}_i) \cap \Delta_{k-1} = \Po \}$ to be the smallest $i$ such that the polynomials in $\mc{G}_i$ cut out exactly the null hypothesis set $\Po$. The lowest degree SOS polynomial corresponding to a SUB test is $2d$.
\end{customthm}
\begin{proof}
We first prove the case where we seek an SOS separating polynomial for an NTUB test. Assume that $\Tilde{\beta} = \sum_{i = 1}^l f_{(i)}^2$ is separating with each $f_{(i)}$ non-zero. As each $f_{(i)}^2$ is non-negative it is clear that each $f_{(i)}$ is in $I_{\mb{R}^{k-1}}(\Po)$ since $\Tilde{\beta}$ must vanish on $\Po$. Any $f_{(i)}$ that is non-zero must be non-zero on $\Pa$ as otherwise $f_{(i)}$ would vanish on $\text{aff}(\Delta_{k-1})$, implying that $f_{(i)}$ is the zero polynomial. Note that this would not be true if it was assumed that $f_{(i)} \in I_{\mb{R}^{k}}(\Po)$ since the polynomial $1 - \sum_{i = 1}^k \pi_i$ vanishes on $\Delta_{k-1}$ but is not the zero polynomial in $\mb{R}^k$. Reparameterizing, by using only $(k-1)$ the coordinates $(\pi_1,\ldots,\pi_{k-1})$, avoids this difficulty. There is no cancellation in the highest degree terms of $\Tilde{\beta}$ and so $2\text{deg}(f_{(i)}) \leq \text{deg}(\Tilde{\beta})$. If $f_1$ has the smallest degree out of $f_{(1)},\ldots,f_{(l)}$, the polynomial $f_{(1)}^2$ is separating with degree less than or equal to $\Tilde{\beta}$. Thus, it suffices to consider the square of a single, non-zero polynomial in $I_{\mb{R}^{k-1}}(\Po)$ with minimal total degree. Let $f \in I_{\mb{R}^{k-1}}(\Po)$ be a non-zero polynomial with the smallest possible degree. As $\{g_{(1)},\ldots,g_{(m)}\}$ is a Gr\"obner basis for $I_{\mb{R}^{k-1}}(\Po)$ the leading term of $f$, $\text{LT}(f)$, is divisible by the leading term of one of the $g_{(i)}$. The leading terms of $f$ and $g_{(i)}$ are the monomials with the largest total degree in $f$ and $g_{(i)}$ respectively, as the monomial order is graded. As $\text{LT}(g_{(i)})$ divides $\text{LT}(f)$ the result $\text{deg}(g_{(i)}) \leq \text{deg}(f)$ follows, and one such minimum degree polynomial in $I_{\mb{R}^{k-1}}(\Po)$ is a minimum degree, non-zero polynomial in $\mc{G}$.

To prove the SUB test case, we first argue that the degree $2\max_{f \in \mc{G}_d} \text{deg}(f)$ polynomial  $\Tilde{\beta} = \sum_{i = 1}^l f_{(i)}^2$, $f_{(i)} \in \mc{G}_d$ gives a polynomial associated with a SUB test. Note that because the hypothesis $\Po$ is assumed to be algebraic $\Po = V_{\mb{R}^{k-1}}(\mc{G}) \cap \Delta_{k-1}$ and there does exist a smallest degree $d$ with $\Po = V_{\mb{R}^{k-1}}(\mc{G}_d) \cap \Delta_{k-1}$. The polynomial $\tilde{\beta}$ vanishes on $\Po$ since $f_{(i)} \in \mc{G}_d \subset I_{\mb{R}^{k-1}}(\Po)$. For any $\bs{\pi}_A \in \Pa$ as $\bs{\pi}_A \notin \Po = V_{\mb{R}^{k-1}}(\mc{G}_d) \cap \Delta_{k-1}$, there exists an $f_{(i)} \in \mc{G}_d$ with $f_{(i)}(\bs{\pi}_A) \neq 0$. This shows that $\Tilde{\beta}(\bs{\pi}_A) > 0$ for every $\bs{\pi}_A \in \Pa$. 

Conversely, let's assume that $\sum_{i = 1}^l f_{(i)}^2$ is a polynomial associated with a SUB test and we wish to show that the degree of this polynomial at least as large as $2\max_{g \in \mc{G}_d} \text{deg}(f)$. Each $f_{(i)}$ must be in $I_{\mb{R}^{k-1}}(\Po)$ as $\Tilde{\beta}$ vanishes on $\Po$. Moreover, $V_{\mb{R}^{k-1}}(f_{(1)},\ldots,f_{(l)}) \cap \Delta_{k-1} = \Po$, as otherwise there would exist a $\bs{\pi}_A \in \Pa$ with $\Tilde{\beta}(\bs{\pi}_A) = 0$. By polynomial long division, every $f_{(i)}$ has a representation in terms of $\mc{G}$ of the form $f_{(i)} = \sum_{j = 1}^m h_{(ij)} g_{(j)}$ where $h_{(ij)} \in \mb{R}[\bs{\pi}]$ and $\text{deg}(g_{(j)}) \leq \text{deg}(f_{(i)})$ for every $j$ with $h_{(ij)} \neq 0$. Also note that as we have a graded ordering $h_{(ij)} = 0$ whenever $\text{deg}(g_{(j)}) > \text{deg}(f_{(i)})$. Let $b = \max_{i \in \{1,\ldots,l\}} \text{deg}(f_{(i)})$. We claim that $d \leq b$, which proves that $2d$ is the lowest attainable degree bound on an SOS representation of $\Tilde{\beta}$. If $b < d$ then $\{f_{(1)}\ldots,f_{(l)} \} \subset \langle \mc{G}_b \rangle$ since each $f_{(i)}$ is in $\langle \mc{G} \rangle$ and $h_{(ij)} = 0$ when $\text{deg}(f_{(j)}) > b$. This containment yields
\begin{align*}
    \Po = V_{\mb{R}^{k-1}}(f_{(1)},\ldots,f_{(l)}) \cap \Delta_{k-1} \supset V_{\mb{R}^{k-1}}(\mc{G}_b) \cap \Delta_{k-1} \supset V_{\mb{R}^{k-1}}(\mc{G}_d) \cap \Delta_{k-1} = \Po,
\end{align*}
contradicting the definition of $d$ as the minimum degree where enough polynomials in $\mc{G}_i$ are included to cut out all of $\Po$.
\end{proof}

%%%%%%%%%%%%%%%%%%%%%%%%%%%%%%%%%%%%%%%%%%
%%%%%%%%%%%%%%%%%%%%%%%%%%%%%%%%%%%%%%%%%%

\begin{customcor}{2}
    Assume that $\Po = V_{\mb{R}^{k-1}}(f_{(1)},\ldots,f_{(m)}) \cap \Delta_{k-1}$ with every $f_{(i)}$ non-zero. The NTUB unbiasedness threshold is bounded above by $2 \min_{i} \text{deg}(f_{(i)})$. The SUB unbiasedness threshold is bounded above by $2\max_{i} \text{deg}(f_{(i)})$. 
\end{customcor}
\begin{proof}
    The proof of this corollary follows from the proof of Theorem \ref{custthm:SOSUBThresholdUpperBound}. In particular, it only needs to be shown that $\Tilde{\beta} = f_{(i)}^2$, where $f_{(i)}$ is the polynomial with the smallest degree in $\{f_{(1)},\ldots,f_{(m)}\}$, and $\Tilde{\beta} = \sum_{i = 1}^m f_{(i)}^2$ respectively are polynomials  associated with NTUB and SUB tests. 
\end{proof}

%%%%%%%%%%%%%%%%%%%%%%%%%%%%%%%%%%%%%%%%%%%%%
%%%%%%%%%%%%%%%%%%%%%%%%%%%%%%%%%%%%%%%%%%%%%%

\begin{customlem}{5}
\label{custlem:RealRadicalParam}
Let $\Po = V_{\mb{R}^{k-1}}(f_{(1)},\ldots,f_{(m)}) \cap \Delta_{k-1}$ and 
    let $\varphi:\mb{C}^p \rightarrow V_{\mb{C}^{k-1}}(f_{(1)},\ldots,f_{(m)})$ be a polynomial parameterization onto a (complex) variety. If $\Po = \varphi(S)$ for some set $S \subset \mb{R}^p$ with a non-empty interior then $I_{\mb{R}^{k-1}}(\Po) =  \sqrt[\mb{C}]{\langle f_{(1)},\ldots,f_{(m)} \rangle} \cap \mb{R}[\bs{\pi}]$.
\end{customlem}
\begin{proof}
    First assume that $f \in I_{\mb{R}^{k-1}}(\Po)$. It follows that $f(\varphi(\bs{x})) = 0$ for all $\bs{x} \in S$. As $S$ has a non-empty interior $f \circ \varphi$ is the zero polynomial. It follows that $f \in I_{\mb{C}^{k-1}}( V_{\mb{C}^{k-1}}(f_{(1)},\ldots,f_{(m)}))$ because the image of $\varphi$ is $V_{\mb{C}^{k-1}}(f_{(1)},\ldots,f_{(m)})$ by assumption.  The Nullstellensatz \cite[Sec 4.2]{cox2013ideals} implies $ f \in \sqrt[\mb{C}]{\langle f_{(1)},\ldots,f_{(m)} \rangle} \cap \mb{R}[\bs{\pi}]$. Conversely, if  $ f \in \sqrt[\mb{C}]{\langle f_{(1)},\ldots,f_{(m)} \rangle} \cap \mb{R}[\bs{\pi}]$
 then by the Nullstellensatz $f \in I_{\mb{C}^{k-1}}( V_{\mb{C}^{k-1}}(f_{(1)},\ldots,f_{(m)}))$, implying that $f$ vanishes on $\Po \subset V_{\mb{C}^{k-1}}(f_{(1)},\ldots,f_{(m)})$. 
 \end{proof}
We remark that in Lemma \ref{custlem:RealRadicalParam} the polynomial parameterization $\varphi$ doesn't have to be onto $V_{\mb{C}^{k-1}}(f_{(1)},\ldots,f_{(m)})$, it is enough that the Zariski closure of the image of $\varphi$ be equal to $V_{\mb{C}^{k-1}}(f_{(1)},\ldots,f_{(m)})$. Also note that we only need $\varphi(S) \subseteq \Po$ and do not need to explicitly force equality. 

%%%%%%%%%%%%%%%%%%%%%%%%%%%%%%%%%%%%%%%%%%%
%%%%%%%%%%%%%%%%%%%%%%%%%%%%%%%%%%%%%%%%%%%

A few comments are needed to explain the definitions used in the following variant of Lemma \ref{custlem:RealRadicalParam}. A rational map $\varphi: W \rightarrow V$ between two complex varieties $W$ and $V$ is a map $\varphi = (\varphi_1,\ldots,\varphi_{k-1})$ where each $\varphi_i$ is a rational function and there exists a point $\bl{x} \in W$ where the denominator of every $\varphi_i$ does not vanish at $\bl{x}$ \citep[Def 5.5.4]{cox2013ideals}. An irreducible variety 
$V = V_{\mb{C}^{k-1}}(f_{(1)},\ldots,f_{(m)})$ is said to be rational if it is birationally equivalent to $\mb{C}^p$. This means that there exists two rational maps $\varphi: \mb{C}^p \rightarrow V$ and $\psi: V \rightarrow \mb{C}^p$ such that $\varphi \circ \psi$ and $\psi \circ \varphi$ equal the identity map, except at points where this composition is undefined due to vanishing denominators. Lastly, a variety $V$ is said to be irreducible if it cannot be written as a non-trivial union $V = V_1 \cup V_2$ of subvarieties with $\emptyset \subsetneq V_i \subsetneq V$, $i = 1, 2$.

\begin{customlem}{5.1}
\label{custlem:RationalParamIdeal}
Let $\Po = V_{\mb{R}^{k-1}}(f_{(1)},\ldots,f_{(m)}) \cap \Delta_{k-1}$ and assume that $V_{\mb{C}^{k-1}}(f_{(1)},\ldots,f_{(m)})$ is an irreducible, rational variety with a rational parameterization given by $\varphi$. 
    % let $\varphi:\mb{C}^p \rightarrow V_{\mb{C}}(f_{(1)},\ldots,f_{(m)})$ be a rational parameterization of a (complex) variety.
    If $\varphi(S) \subseteq \Po$ for some set $S \subset \mb{R}^p$ with a non-empty interior then $I_{\mb{R}^{k-1}}(\Po) =  \sqrt[\mb{C}]{\langle f_{(1)},\ldots,f_{(m)} \rangle} \cap \mb{R}[\bs{\pi}]$.
\end{customlem}
\begin{proof}
    First assume that $f \in I_{\mb{R}^{k-1}}(\Po)$. It follows that $f(\varphi(\bs{x})) = 0$ for all $\bs{x} \in S$. As $S$ has a non-empty interior $f \circ \varphi$ is the zero rational function. It follows that $f$ vanishes on the image $\varphi(\mb{C}^p)$. As $V_{\mb{C}^{k-1}}(f_{(1)},\ldots,f_{(m)})$ is rational there exists a $\psi$ with $\varphi \circ \psi = \text{Id}$ outside of a subvariety $V' \subsetneq V_{\mb{C}^{k-1}}(f_{(1)},\ldots,f_{(m)})$, implying that for $\bl{x} \notin V'$, $\varphi \circ \psi (\bl{x}) = \bl{x}$. Thus $f(\bl{x}) = 0$ for all $\bl{x} \in  V_{\mb{C}^{k-1}}(f_{(1)},\ldots,f_{(m)}) \backslash V'$. As $V_{\mb{C}^{k-1}}(f_{(1)},\ldots,f_{(m)})$ is irreducible $f$ must vanish on all of $V_{\mb{C}^{k-1}}(f_{(1)},\ldots,f_{(m)})$, implying that $f$ is in $I_{\mb{C}^{k-1}}(V_{\mb{C}^{k-1}}(f_{(1)},\ldots,f_{(m)}))$. The remainder of the proof is identical to Lemma \ref{custlem:RealRadicalParam}. 
 \end{proof}

%%%%%%%%%%%%%%%%%%%%%%%%%%%%%%%%%%%%%%%%%%%%%%%%%%%%%%%%%%%%%%%%%%%
%%%%%%%%%%%%%%%%%%%%%%%%%%%%%%%%%%%%%%%%%%%%%%%%%%%%%%%%%%%%%%%%%%

\begin{customthm}{5}
\label{custthm:PrincipalIdealUBThresholdUMP}
    Suppose that $\Po$ is algebraic, $\overline{\text{int}(\tilde{\Delta}_{k-1}) \cap \Po} = \Po$, and $I_{\mb{R}^{k-1}}(\Po) = \langle f \rangle$ is a principal ideal generated by a single polynomial. If $f$ has a non-zero gradient throughout $\Po$, the NTUB and SUB unbiasedness thresholds are both equal to $2 \,\mathrm{deg}(f)$. For every size $0<\alpha<1$ when the sample size is $n = 2 \,\mathrm{deg}(f)$ there is a unique, uniformly most powerful NTUB test with size-$\alpha$ that is also the unique, uniformly most powerful SUB test. This test has a power function of the form $\beta = c_{\alpha} f^2 + \alpha$, where $c_{\alpha} \in (0,\infty)$ is the largest constant that makes $\beta$ a power polynomial. Every size-$\alpha$ NTUB test of $\Po$ has a power function of the form $cf^2 + \alpha$, $c \in (0,c_\alpha]$ when $n = 2 \,\mathrm{deg}(f)$. More generally, if $\beta$ is the power function of a level-$\alpha$, NTUB test for a sample size of $n$ there exists a polynomial $h$ that is non-negative on $\Delta_{k-1}$ with degree $n - 2\text{deg}(f)$ such that $\beta = f^2h + \alpha$. For a SUB test $h$ must in addition be strictly positive on $\Pa$.
\end{customthm}
\begin{proof}
    Any polynomial $\Tilde{\beta}$ in $I_{\mb{R}^{k-1}}(\Po)$ associated with an NTUB test has the form $g(\bs{\pi})f(\bs{\pi})$. Set $\Po' \coloneqq \text{int}(\tilde{\Delta}_{k-1}) \cap \Po$.
    % Take $\Tilde{\Delta}_{k-1} = \{\bs{\pi} \in \mb{R}^{k-1}: \sum_{i = 1}^{k-1}\pi_i < 1, \pi_i \geq 0\}$ to be the open probability simplex. All polynomials in this question will be viewed as functions of $(\pi_1,\ldots,\pi_{k-1}) \in \Tilde{\Delta}_{k-1}$ defined by $\Tilde{\beta}(\pi_1,\ldots,\pi_{j-1},1 - \sum_{i = 1}^{k-1} \pi_i)$.
    As $\Tilde{\beta}$ is locally minimized on $\Po'$ the directional derivative at $\bs{\pi}_0 \in \Po'$ along any feasible direction $\bl{v}$ with $\bs{\pi}_0 + \bl{v} \in \Tilde{\Delta}_{k-1}$ is zero. Thus
    \begin{align*}
        \partial_{\bl{v}} (gf)(\bs{\pi_0}) = \langle \nabla g(\bs{\pi}_0),\bl{v}\rangle f(\bs{\pi}_0) +  \langle \nabla f(\bs{\pi}_0),\bl{v} \rangle g(\bs{\pi}_0) =   \langle \nabla f(\bs{\pi}_0),\bl{v}\rangle g(\bs{\pi}_0)
    \end{align*}
    for every $\bl{v}$ that is a feasible direction at $\bs{\pi}_0$. It is claimed that there exists a feasible $\bl{v}$ where $\langle \nabla f(\bs{\pi}_0), \bl{v}\rangle  \neq 0$. Regardless of the point $\bs{\pi}_0$, the space of feasible directions at $\bs{\pi}_0$ spans $\mb{R}^{k-1}$. Hence, $\langle \nabla f(\bs{\pi}_0), \bl{v}\rangle = 0$ for all feasible $\bl{v}$ would imply the contradiction that $ \nabla f(\bs{\pi}_0) = 0$. As there must exist a $\bl{v}$ with $\langle \nabla f(\bs{\pi}_0), \bl{v}\rangle  \neq 0$ we conclude that $g(\bs{\pi}_0) = 0$. As $g$ vanishes on $\Po'$, by continuity, $g$ also vanishes on $\overline{\Po'} = \Po$ and so $g = fh \in I_{\mb{R}^{k-1}}(\Po)$. This gives the general representation $\tilde{\beta} = f^2 h$, or equivalently $\beta = f^2 h + \alpha$. The polynomial $h$ must be non-negative on $\Pa$ for if $h(\bs{\pi}_A) < 0$ then $\beta(\bs{\pi}_A) = f^2(\bs{\pi}_A) h(\bs{\pi}_A) + \alpha < \alpha$, contradicting the unbiasedness of $\beta$. As $\Pa$ is dense in $\Delta_{k-1}$, by continuity $h$ is non-negative on $\Delta_{k-1}$. In addition, for a SUB test $h$ cannot vanish on $\Pa$ as then $\beta$ would have power $\alpha$ at a point in $\Pa$. The NTUB unbiasedness threshold is the smallest degree of a non-constant polynomial that  has the form $f^2 h + \alpha$. Every such polynomial has degree $2\,\text{deg}(f)  + \text{deg}(h)$, which is minimized when $h$ is a non-negative constant $c$.
    % As $f^2$ has degree $2\, \mathrm{deg}(f)$ and is a separating polynomial for a SUB test, the SUB unbiasedness threshold is also $2\, \mathrm{deg}(f)$.    
    
    % Every non-zero polynomial in $\langle f \rangle$ has degree no smaller than $\mathrm{deg}(f)$, from which we conclude that the unbiasedness threshold is at least $2\, \mathrm{deg}(f)$. As $f^2$ has degree $2\, \mathrm{deg}(f)$ and is a separating polynomial for a SUB test, the SUB unbiasedness threshold is $2\, \mathrm{deg}(f)$.   

This proof shows that the only non-negative and non-zero polynomials on the simplex that are in $I_{\mb{R}^{k-1}}(\Po)$ and have degree $2\, \text{deg}(f)$ have the form $cf^2$ for $c \in (0,\infty)$. If there exists a SUB power polynomial $\beta$, by recalling that a separating polynomial is obtained as $\Tilde{\beta} = \beta - \alpha$, there must exist some $c$ such that $\beta = cf^2 + \alpha$. Conversely, if $cf^2 + \alpha$ is to be a power polynomial the coefficient of  $\bs{\pi}^{\bl{x}}$ in $cf^2 + \alpha(\sum_{i = 1}^k \pi_i)^{2 \text{deg}(f)}$ must be in the interval $[0, {n \choose \bl{x}}]$. As the coefficient of $\bs{\pi}^{\bl{x}}$ in $\alpha(\sum_{i = 1}^k \pi_i)^{2 \text{deg}(f)}$ is $\alpha {n \choose \bl{x}}$, there exists a small enough positive $c$ that yields a power polynomial. Moreover, as $f^2 \geq 0$ if $c_1 \geq c_2$ then $c_1f^2 + \alpha \geq c_2 f^2 + \alpha$ on $\Tilde{\Delta}_{k-1}$. The uniformly most powerful SUB and NTUB test therefore has a power polynomial $\beta = c_\alpha f^2 + \alpha$ where $c_\alpha$ is the largest $c > 0$ such that $c_\alpha f^2 + \alpha$ is a valid power polynomial.  
\end{proof}

%%%%%%%%%%%%%%%%%%%%%%%%%%%%%%%%%%%%%%%%%%%%%
%%%%%%%%%%%%%%%%%%%%%%%%%%%%%%%%%%%%%%%%%%%%%
%%%%%%%%%%%%%%%%%%%%%%%%%%%%%%%%%%%%%%%%%%%%%

\begin{customthm}{6}
    Let $\Po = \{\bs{\pi} \in \Delta_{k-1}: f(\bs{\pi}) \leq 0\}$ and assume that $I_{\mb{R}^{k-1}}( \partial \Po) = \langle f \rangle$ where $\partial \Po$ is the relative boundary of $\Po$ in $\Delta_{k-1}$.  The NTUB and SUB unbiasedness thresholds are both equal to $\mathrm{deg}(f)$. For every size $0<\alpha<1$, when the sample size is $n =  \,\mathrm{deg}(f)$ there is a unique, uniformly most powerful NTUB test that is also the unique, uniformly most powerful SUB test. This test has a power function of the form $\beta = c_{\alpha} f + \alpha$, where $c_{\alpha} \in (0,\infty)$ is the largest constant that makes $\beta$ a power polynomial. Every size-$\alpha$ NTUB test of $\Po$ has a power function of the form $cf + \alpha$, $c \in (0,c_\alpha]$ when $n = \mathrm{deg}(f)$.  More generally, if $\beta$ is the power function of a level-$\alpha$, NTUB test for a sample size of $n$ there exists a polynomial $h$ that is non-negative on $\Delta_{k-1}$ with degree $n - \text{deg}(f)$ such that $\beta = fh + \alpha$. For a SUB test $h$ must in addition be strictly positive on $\Pa$.
\end{customthm}
\begin{proof}
    The proof is nearly identical to that of Theorem \ref{custthm:PrincipalIdealUBThresholdUMP}. The assumptions that the gradient of $f$ is non-zero on $\Po$ and that the closure of $\Po \cap \text{int}(\tilde{\Delta}_{k-1})$ is equal to $\Po$ are no longer needed. Any separating polynomial must vanish on $\partial \Po$ and so it must lie in the ideal $\langle f \rangle$. Thus, $\beta = fh + \alpha$ for some $h$ when $\beta$ is an NTUB power polynomial. To show that $h$ is non-negative assume by contradiction that $h(\bs{\pi}) < 0$ for some $\bs{\pi} \in \Delta_{k-1}\backslash \partial \Po$. If $\bs{\pi} \in \Po$ then $\beta(\bs{\pi}) = f(\bs{\pi})h(\bs{\pi}) + \alpha > \alpha$, while if $\bs{\pi} \in \Pa$ then $\beta(\bs{\pi}) = f(\bs{\pi})h(\bs{\pi}) + \alpha < \alpha$. Either of these inequalities contradict the unbiasedness of $\beta$. By continuity $h \geq 0$ on $\Delta_{k-1}\backslash \partial \Po$ implies that $h \geq 0$ on $\Delta_{k-1}$. Again, $h$ cannot vanish on $\Pa$ for the associated test to be SUB. 

The non-constant polynomial of the form $fh + \alpha$ with the smallest degree is $cf + \alpha$ for some $c > 0$. This power polynomial, when appropriately normalized, is also the power polynomial for a SUB test as $f$ is strictly positive on $\Pa$. This strict positivity implies that the power function is pointwise increasing on $\Pa$ as $c$ increases. It follows that the unique UMPU power polynomial is obtained when $c$ is the largest possible value such that $cf + \alpha$ is a power polynomial. 
% It is easily seen that $cf$ is a separating polynomial only when $c > 0$.
% Renormalizing this polynomial, any NTUB power polynomial must equal $\beta = cf + \alpha$ for an appropriate choice of $c > 0$. 
% This power polynomial is also the power polynomial for a SUB test as $f$ is strictly positive on $\Pa$. This strict positivity implies that the power function is pointwise increasing on $\Pa$ as $c$ increases. It follows that the unique UMPU power polynomial is obtained when $c$ is the largest possible value such that $cf + \alpha$ is a power polynomial. 
\end{proof}

%%%%%%%%%%%%%%%%%%%%%%%%%%%%%%%%%%%%%%%%%%%%%%%%
%%%%%%%%%%%%%%%%%%%%%%%%%%%%%%%%%%%%%%%%%%%%%%%
%%%%%%%%%%%%%%%%%%%%%%%%%%%%%%%%%%%%%%%%%%%%%%%

\begin{customcor}{3}
     Assume that $f_{(1)},\ldots,f_{(l)}$ are distinct, irreducible polynomials, $I_{\mb{R}^{k-1}}(\mathcal{P}_i) = \langle f_{(i)} \rangle$, each $\mathcal{P}_i$ is algebraic, $\overline{\text{int}(\tilde{\Delta}_{k-1}) \cap \mathcal{P}_i} = \mathcal{P}_i$,  and $\nabla f_{(i)}$ is non-zero on $\mathcal{P}_i$. Then $I_{\mb{R}^{k-1}}(\cup_{i = 1}^l \mathcal{P}_i) = \langle \prod_{i = 1}^l f_{(i)} \rangle$ and the SUB and NTUB thresholds of $\Po = \cup_{i = 1}^l \mathcal{P}_i$ are $2 \sum_{i = 1}^l\text{deg}(f_{(i)})$. 
\end{customcor}
\begin{proof}
The claim that $I_{\mb{R}^{k-1}}(\cup_{i = 1}^l \mathcal{P}_i) = \langle \prod_{i = 1}^l f_{(i)} \rangle$ is a standard result above principal ideals \citep[Prop 4.3.13]{cox2013ideals} where we note that $\text{LCM}(f_{(1)},\ldots,f_{(l)}) = \prod_{i = 1}^l f_{(i)}$ . It remains to show that the proof of Theorem \ref{custthm:PrincipalIdealUBThresholdUMP} carries through. Note that
\begin{align*}
    \nabla \bigg( \prod_{i = 1}^l f_{(i)} \bigg) = \sum_{i = 1}^l \bigg(\prod_{j \neq i} f_{(j)} \bigg) \nabla f_{(i)},
\end{align*}
and this gradient is non-zero at all points $\bs{\pi} \in \mathcal{P}_i \backslash \big(\cup_{j \neq i} \mathcal{P}_j\big)$ for every $i = 1,\ldots,l$. By the proof of Theorem \ref{custthm:PrincipalIdealUBThresholdUMP} the polynomial $h$ will vanish on $\bs{\pi} \in \mathcal{P}_i \backslash \big(\cup_{j \neq i} \mathcal{P}_j\big)$. The polynomial $h \prod_{j \neq i} f_{(j)}$ must then vanish on $\mc{P}_i$, implying that $f_{(i)}$ divides $h \prod_{j \neq i} f_{(j)}$. As $f_{(i)}$ does not divide any of the $f_{(j)}$, by irreducibility it must divide $h$, showing that $h \in \langle f_{(i)} \rangle$. This holds for every $i$ so $h \in \langle \prod_{i = 1}^l f_{(i)} \rangle$ as needed. 
\end{proof}

%%%%%%%%%%%%%%%%%%%%%%%%%%%%%%%%%%%%%%%%%%%%%
%%%%%%%%%%%%%%%%%%%%%%%%%%%%%%%%%%%%%%%%%%%%%
%%%%%%%%%%%%%%%%%%%%%%%%%%%%%%%%%%%%%%%%%%%%%%

\begin{customlem}{6}
    Let $\mathcal{P}_i$ have the NTUB and SUB thresholds of $d_{1i}, d_{2i}$ respectively with corresponding separating polynomials $f_{(1i)}$ and $f_{(2i)}$ for $i = 1,\ldots,l$. The polynomials $\prod_{i = 1}^l f_{(1i)}$ and $\prod_{i = 1}^l f_{(2i)}$ are NTUB and SUB separating polynomials with degrees $\sum_{i = 1}^l d_{1i}$,  $\sum_{i = 1}^l d_{2i}$ for the union null hypothesis $\Po = \cup_{i = 1}^l \mathcal{P}_i$.
\end{customlem}
\begin{proof}
    Both $\prod_{i = 1}^l f_{(1i)}$ and $\prod_{i = 1}^l f_{(2i)}$ vanish on $\Po$ since on $\mathcal{P}_i$ the factors $f_{(1i)}$ and $f_{(2i)}$ vanish. It is clear that $\prod_{i = 1}^l f_{(2i)} > 0$ on $\Pa$. It remains to show that $\prod_{i = 1}^l f_{(1i)}$ is non-zero on $\Delta_{k-1}$, or equivalently that $1 -\sum_{i = 1}^k \pi_i$ does not divide this polynomial. As $1 -\sum_{i = 1}^k \pi_i$ is irreducible if it divides $\prod_{i = 1}^l f_{(1i)}$ it divides one of the factors $f_{(1j)}$. This implies that $f_{(1j)}$ is zero on $\Delta_{k-1}$ contradicting the assumption that $f_{(1j)}$ is an NTUB separating polynomial. 
\end{proof}

%%%%%%%%%%%%%%%%%%%%%%%%%%%%%%%%%%%%%%
%%%%%%%%%%%%%%%%%%%%%%%%%%%%%%%%%%%%%
%%%%%%%%%%%%%%%%%%%%%%%%%%%%%%%%%%%%%%

\begin{customthm}{7}
    The NTUB and SUB unbiasedness thresholds for the null hypothesis $\Po = \{ \bs{\pi} \in \Delta_{pq-1}: \bs{\pi} \; \mathrm{has \; rank \; less \; than } \; r \}$ on the space of $p \times q$ dimensional contingency tables is equal to $2r$. When $n = 2r$ one such SUB power polynomial has the form $\beta(\bs{\pi}) = c\sum_{I,J} \det(\bs{\pi}_{I,J})^2 + \alpha$ for some $c > 0$, where $I,J$ are index sets of size $r$.  
\end{customthm}
\begin{proof}  The polynomial $\sum_{I,J} \det(\bs{\pi}_{I,J})^2$ is a SUB separating polynomial with total degree $2r$. Let $f$ be a NTUB separating polynomial of minimal total degree, and assume for a contradiction that the degree of $f$ is less than $2r$.  By Lemma \ref{custlem:RadicalofDeterminantalIdeal} $f$ has the form $f(\bs{\pi}) = \sum_{I,J} h_{I,J}(\bs{\pi}) \det(\bs{\pi}_{I,J})$ for some polynomials $h_{I,J}$. As the $\det(\bs{\pi}_{I,J})$ are a Gr\"obner basis for $I_{\mb{R}^{k-1}}(\Po)$ the $h_{I,J}$ are determined by polynomial division of $f$ by the $\det(\bs{\pi}_{I,J})$ and can be assumed to all have degrees less than $r$. The polynomial $f$ can without loss of generality be assumed to be  symmetric with respect to row and column permutations of the matrix $\bs{\pi}$ since the following symmetrized version of $f$ ($\bl{U},\bl{V}$ below run over the set of all permutation matrices):
\begin{align*}
   \frac{1}{p! q!} \sum_{ \bl{U} \in \mc{S}_p, \bl{V} \in \mc{S}_q } f(\bl{U} \bs{\pi} \bl{V}) =   \frac{1}{p! q!} \sum_{ \bl{U} \in \mc{S}_p, \bl{V} \in \mc{S}_q } h_{I,J}(\bl{U} \bs{\pi} \bl{V}) \det([\bl{U} \bs{\pi} \bl{V}]_{I,J}) 
\end{align*}
is a non-zero polynomial if $f$ is. This is because the alternative hypothesis set is invariant under row and column permutations so if $f(\bs{\pi}_A) > 0$ for some $\bs{\pi}_A \in \Pa$ then $ \frac{1}{p! q!} \sum_{ \bl{U} \in \mc{S}_p, \bl{V} \in \mc{S}_q } f(\bl{U} \bs{\pi}_A \bl{V}) > 0$. Furthermore, each term $ \det([\bl{U} \bs{\pi} \bl{V}]_{I,J})$ in the symmetrized polynomial is equal to $\pm \det(\bs{\pi}_{I',J'})$ for some $I',J'$. Without loss of generality it can be assumed, by grouping all of the terms $\pm \det(\bs{\pi}_{I,J})$ of the symmetrized polynomial together, that each $h_{I,J}(\bs{\pi})$ has the property that $h_{I,J}$ is antisymmetric with respect to any row and column permutations that permute row indices in $I$ and columns indices in $J$ while leaving the indices in $I^c = [p] \backslash I$ and $J^c = [q] \backslash J$ fixed. The goal will be to show that one of the $h_{I,J}$ terms must be divisible by $\det(\bs{\pi}_{I,J})$.     

We examine the polynomial $f(\bs{\pi}_{I,J},\bl{0}) = h_{I,J}(\bs{\pi}_{I,J},\bl{0}) \det(\bs{\pi}_{I,J})$ restricted to the subset of matrices where the $\bs{\pi}_{I,J}$ block is non-zero but all other entries of $\bs{\pi}$ are zero. There are two cases. The first case is that $f(\bs{\pi}_{I,J},\bl{0})$ is non-zero and consequently is a separating polynomial for the null hypothesis $\Po' = \{ \bs{\pi}_{I,J} \in \Delta_{r^2 - 1}: \det(\bs{\pi}_{I,J}) = 0\}$. This null hypothesis set has a principal ideal $I_{\mb{R}^{k-1}}(\Po') = \langle \det(\bs{\pi}_{I,J}) \rangle$ by Lemma \ref{custlem:RadicalofDeterminantalIdeal} and so by Theorem \ref{custthm:PrincipalIdealUBThresholdUMP} we conclude that $\det(\bs{\pi}_{I,J})^2$ divides $f(\bs{\pi}_{I,J},\bl{0})$, and so $\det(\bs{\pi}_{I,J})$ divides $h_{I,J}(\bs{\pi}_{I,J}, \bl{0})$, providing a contradiction to $\text{deg}(h_{I,J}) < r$.  The second case has $h_{I,J}(\bs{\pi}_{I,J},\bl{0}) = 0$ for all index sets $I,J$. It follows that every monomial $\pi^L = \prod_{j = 1}^{r-1} \pi_{s_j t_j}$ appearing in $h_{I,J}(\bs{\pi})$ must contain a factor with a row or column index in $I^c$ or $J^c$ respectively. Assume without loss of generality that there is a row index $s_1,\ldots,s_{r-1}$ contained in $I^c$. As $h_{I,J}$ has degree at most $r-1$ this leaves at most $r-2$ factors in $\pi^L$ that have row indices contained in $I$. Let $i,i' \in I$ be two indices that do not appear in $s_1,\ldots,s_{r-1}$. Define 
\begin{align*}
    A_{I,J}(g(\bs{\pi})) = \frac{1}{(r!)^2} \sum_{\bl{U} \in \mc{S}_r(I), \bl{V} \in \mc{S}_r(J)} \text{sign}(\bl{U}) \text{sign}(\bl{V}) g(\bl{U}\bs{\pi}\bl{V}),  
\end{align*}
to be the linear, antisymmetrization operator that is a sum over all the permutation matrices that permute the indices appearing in $I,J$ while leaving all the indices in $I^c$ and $J^c$ fixed. As $h_{I,J}$ was assumed to already be antisymmetric $A_{I,J}(h_{I,J}) = h_{I,J}$. If $\tau$ is the transposition that swaps $i$ with $i'$ then $\pi^{\tau(L)} \coloneqq \prod_{j = 1}^{k-1} \pi_{\tau(s_j)t_j} =  \prod_{j = 1}^{k-1} \pi_{\tau(s_j)t_j} = \pi^L$. That is, the polynomial $\pi^L$ is fixed by this transposition. However,
\begin{align*}
   A_{I,J}(\pi^L) =  A_{I,J}( \pi^{\tau(L)}) =   \text{sign}(\tau)A(\pi^L) = -A_{I,J}(\pi^L),
\end{align*}
implying that $A_{I,J}(\pi^L) = 0$. By linearity we conclude that $h_{I,J} = A_{I,J}(h_{I,J}) = A_{I,J}(\sum_L c_L \pi^L) = 0$. This provides the needed contradiction that $h_{I,J} = 0$ for every $I,J$. 
\end{proof}

%%%%%%%%%%%%%%%%%%%%%%%%%%%%%%%%%%%%%%%%%%%%%%%
%%%%%%%%%%%%%%%%%%%%%%%%%%%%%%%%%%%%%%%%%%%%%%%%%
%%%%%%%%%%%%%%%%%%%%%%%%%%%%%%%%%%%%%%%%%%%%%%%%%%

The parameterization $\varphi(\bl{a},\bl{b}) = (\bl{a}^\intercal, 1 - \bl{1}_{p-1}^\intercal \bl{a})^\intercal(\bl{b}^\intercal, 1 - \bl{1}_{q-1}^\intercal \bl{b}) \in \mb{R}^{p \times q}$ can be used in Lemma \ref{custlem:RealRadicalParam} to find the ideal of the set of set of matrices that have rank one and are contained in affine hull of the probability simplex. More care has to be taken when constructing a parameterization of matrices of probabilities that have rank less than $r$ when $r > 2$. Specifically, if $r = 3$ it is not sufficient to set $\pi_{ij}(\bl{a}_1,\bl{a}_2,\bl{b}_1,\bl{b}_2) = \bl{a}_{1i} \bl{b}_{1j}^\intercal + \bl{a}_{2i} \bl{b}_{2j}^\intercal$ for every $i$ and $j$ except for $\pi_{pq}$ which is set to $\pi_{pq} = 1- \sum_{(i,j) \neq (p,q)} \pi_{ij}(\bl{a}_1,\bl{a}_2,\bl{b}_1,\bl{b}_2)$. By enforcing the sum-to-one constraint on $\pi_{pq}$ the resulting matrix may have a rank that is larger than $2$ and the parameterization has an image that is not contained in the model.

\begin{customlem}{7.1}
\label{custlem:RadicalofDeterminantalIdeal}
 If $\Po = \{ \bs{\pi} \in \Delta_{pq-1}: \bs{\pi} \; \mathrm{has \; rank \; less \; than } \; r \}$ then $I_{\mb{R}^{k-1}}(\Po) = \langle \det(\bs{\pi}_{I,J}): I \in [p]^r, J \in [q]^r \rangle$. 
\end{customlem}
\begin{proof}
To ease notation we define $I = \langle 1 - \sum_{i,j}\pi_{ij}, \det(\bs{\pi}_{I,J}): I \in [p]^r, J \in [q]^r \rangle$
The set of matrices with rank less than $r \leq \min(p,q)$ is parameterized by $\mb{R}^{p \times r-1} \times \mb{R}^{q \times r-1}$ via the map $\psi(\bl{A},\bl{B}) = \bl{A}\bl{B}^\intercal$. The only issue with this parameterization is that $\Po$ is not parameterized by a set with non-empty interior in $\mb{R}^{p \times r-1} \times \mb{R}^{q \times r-1}$ since $\Po$ has an additional sum-to-one constraint that is not respected by $\varphi$. To account for this, note that if $\bl{A},\bl{B}$ have columns $\bl{a}_{\cdot i}, \bl{b}_{\cdot i}$ the sum-to-one constraint is equivalent to $\sum_{i = 1}^{r-1} (\bl{a}_{\cdot i}^\intercal \bl{1})(\bl{b}_{\cdot i}^\intercal \bl{1}) = 1$ ,which allows $A_{11}$ to be expressed in terms of the other parameters as long as $(\bl{b}_{\cdot 1}^\intercal \bl{1})$ is non-zero. Defining
\begin{align*}
    A_{11} = \frac{1}{(\bl{b}_{\cdot 1}^\intercal \bl{1})}\bigg(1 - \sum_{i = 2}^{r-1} (\bl{a}_{\cdot i}^\intercal \bl{1})(\bl{b}_{\cdot i}^\intercal \bl{1}) \bigg) - \sum_{i = 2}^p \bl{A}_{i1},
\end{align*}
we obtain a parameterization $\varphi$ of $\Po$ in terms of rational functions. This parameterization extends to the complex numbers where $\varphi(\mb{C}^{p(r-1)-1} \times \mb{C}^{q \times r-1})$ is contained in the set of matrices that have rank less than $r$ and have entries that sum to one. It is claimed that the image of this parameterization is onto this set.  Let $\bl{M} \in \mb{C}^{p\times q}$ be a matrix with rank less than $r$ with entries that sum to one. By the rank assumption there exists $\bl{A},\bl{B}$ with  $\bl{M} = \bl{A} \bl{B}^\intercal$. If $\bl{b}_{\cdot 1}^\intercal \bl{1} \neq 0$ the sum-to-one constraint that $\bl{M}$ satisfies implies that $\bl{M}$ is in the image of $\varphi$. There must exist one column $i^*$ of $\bl{B}$ with $\bl{b}_{\cdot i^*}^\intercal \bl{1} \neq 0$ as otherwise $\bl{M} \bl{1} = \bl{0}$, contradicting the sum to one assumption. Let $\bl{V}$ be a permutation matrix that swaps $1$ with $i^*$. Defining $\Tilde{\bl{A}} = \bl{A}\bs{V}$ and $\Tilde{\bl{B}} = \bl{B}\bl{V}$ we have that $\bl{M} = \Tilde{\bl{A}} \tilde{\bl{B}}^\intercal$ where the first column of $\Tilde{\bl{B}}$ has $\Tilde{\bl{b}}_{\cdot 1}^\intercal \bl{1} \neq 0$. We conclude that the image of $\varphi$ is onto the set of matrices with rank less than $r$ that have entries that sum to one. As this parameterization is onto, the proof of Lemma \ref{custlem:RationalParamIdeal} can be applied to show that $I_{\mb{R}^{k-1}}(\Po) = \sqrt[\mb{C}]{I} \cap \mb{R}[\bs{\pi}]$. 
% \begin{align*}
% I_{\mb{R}^k}(\Po) = \sqrt[\mb{C}]{I} \cap \mb{R}[\bs{\pi}].
% \end{align*}

We now show that $\sqrt[\mb{C}]{I} = I$ is radical. It is known \citep[Cor 16.29]{MillerSturmfelsCombinatorial} that the ideal $\langle  \det(\bs{\pi}_{I,J}): I \in [p]^r, J \in [q]^r \rangle$ is radical and is the vanishing ideal of all matrices with rank less than $k$. We make a change of coordinates $\bs{\pi} \mapsto \bl{X} \coloneqq \bl{C}\bs{\pi} \bl{D}$. The invertible matrices $\bl{C},\bl{D}$ are chosen so that the first row of $\bl{C}$ is $\bl{1}^\intercal$ and the first column of $\bl{D}$ is $\bl{1}$. In particular, $X_{11} = \sum_{i,j} \pi_{ij}$. This change of coordinates maps the ideal $I$ to the ideal $J =  \langle 1 - X_{11}, \det(\bl{X}_{I,J}): I \in [p]^r, J \in [q]^r \rangle$. To see that $\langle  \det(\bs{\pi}_{I,J}): I \in [p]^r, J \in [q]^r \rangle$ gets mapped to $\langle  \det(\bl{X}_{I,J}): I \in [p]^r, J \in [q]^r \rangle$ note that the rank of $\bs{\pi}$ is the same as the rank of $\bs{X}$ and so the vanishing ideal of all matrices $\bs{\pi}$ with rank less than $k$ will get mapped to the vanishing ideal of all matrices $\bl{X}$ with rank less than $k$. As the coordinate change is an isomorphism of the rings $\mb{C}[\bs{\pi}]$ and $\mb{C}[\bl{X}]$ the ideal $I$ is radical if and only if $J$ is.  We show that $\{1 - X_{11}, \det(\bl{X}_{I,J}): I \in [p]^r, J \in [q]^r \}$ is a Gr\"obner basis of $J$ by showing that the remainder of every S-polynomial under polynomial long-division by this set is zero \citep[Thm 3.6.6]{cox2013ideals}. In \citep[Thm 16.28]{MillerSturmfelsCombinatorial} it is shown that $\{ \det(\bl{X}_{I,J}): I \in [p]^r, J \in [q]^r\}$ is a Gr\"obner basis for the ideal $\langle \det(\bl{X}_{I,J}): I \in [p]^r, J \in [q]^r \rangle$ under any antidiagonal term ordering. This implies that any S-polynomial constructed from $\det(\bl{X}_{I,J})$ and $\det(\bl{X}_{I',J'})$ will have a remainder of zero upon dividing by $\{ \det(\bl{X}_{I,J}): I \in [p]^r, J \in [q]^r\}$. The other kind of S-polynomial is one constructed from $\det(\bl{X}_{I,J})$ and $X_{11} - 1$, which has the form $P(\bl{X}) = X_{11}\big( \det(\bl{X}_{I,J}) - \text{LT}(\det(\bl{X}_{I,J}))\big) + \text{LT}(\det(\bl{X}_{I,J}))$. We choose a graded, lexicographic, antidiagonal, term ordering with $X_{11} \leq X_{ij}$ for all $i,j$. The leading term of $P(\bl{X})$ under this term ordering is
\begin{align*}
    Q(\bl{X}) \coloneqq X_{11}\text{LT}\big(\det(\bl{X}_{I,J}) - \text{LT}(\det(\bl{X}_{I,J}))\big),
\end{align*}
 as the ordering is graded. None of the leading terms of $\det(\bl{X}_{I',J'})$ divide $Q(\bl{X})$ since no such leading term includes $X_{11}$ by the lexicographic assumption and no such leading term divides $\text{LT}(\det(\bl{X}_{I,J}) - \text{LT}(\det(\bl{X}_{I,J})))$ either, since $\text{LT}(\det(\bl{X}_{I',J'}))$ and $\text{LT}(\det(\bl{X}_{I,J}) - \text{LT}(\det(\bl{X}_{I,J})))$ are distinct degree $k$ monomials. The division algorithm thus passes through all of the $\det(\bl{X}_{I',J'})$ terms until $X_{11} - 1$ is reached. The leading term $Q(\bl{X})$ is divided by $X_{11}-1$. The whole process then repeats, where again only $X_{11}-1$ divides the new leading term. In aggregate, we repeatedly divide the terms of $X_{11}\big( \det(\bl{X}_{I,J}) - \text{LT}(\det(\bl{X}_{I,J}))\big)$ in $P(\bl{X})$ by $X_{11} - 1$ until $\big( \det(\bl{X}_{I,J}) - \text{LT}(\det(\bl{X}_{I,J}))\big)$ is obtained. Adding this to the term that has no yet been touched, $\text{LT}(\det(\bl{X}_{I,J}))$, we arrive at $\det(\bl{X}_{I,J})$, which clearly divides to zero.  The conclusion is that the remainder of every S-polynomial upon division by the basis is zero, proving that we have a Gr\"obner basis.  

 To complete the proof we note that the initial ideal of $J$ is generated by the leading terms of the Gr\"obner basis: $\langle X_{11}, \text{LT}(\det(\bl{X}_{I,J})) : I \in [p]^r, J \in [q]^r \rangle$. The initial ideal is squarefree \citep[Def 1.3]{MillerSturmfelsCombinatorial}, implying that it is radical. This further implies \citep[Prop 4.2.11]{SullivantAlgStats} that $J$ is radical, as required. Finally we note that the ideal $I_{\mb{R}^{k-1}}(\Po)$ is obtained from $I$ by the ring isomorphism that substitutes $1 - \sum_{(i,j) \neq (p,q)} \pi_{ij}$ for any occurrence of $\pi_{pq}$ appearing in a polynomial in $I$. 
\end{proof}

%%%%%%%%%%%%%%%%%%%%%%%%%%%%%%%%%%%%%%%%%%%%%%%%%%%%
%%%%%%%%%%%%%%%%%%%%%%%%%%%%%%%%%%%%%%%%%%%%%%%%%%%%

\begin{customthm}{8}
  Suppose that $\Po$ is algebraic, $\overline{\mathrm{int}(\tilde{\Delta}_{k-1}) \cap \Po} = \Po$, $I_{\mb{R}^{k-1}}(\Po) = \langle f_{(1)},\ldots,f_{(m)} \rangle$, and the matrix $[\nabla f_{(1)} \cdots \nabla f_{(m)}] \in \mb{R}^{k-1 \times m}$ has rank $m$ on $\Po$, then every NTUB test of $\Po$ for a sample size of $n$ has a power polynomial of the form $\beta(\bs{\pi}) = \bl{f}^\intercal(\bs{\pi}) \bl{H}(\bs{\pi}) \bl{f}(\bs{\pi}) + \alpha$. Here $\bl{f}^\intercal = (f_{(1)},\ldots,f_{(m)})$, the entries of $\bl{H}$ are polynomials $h_{(ij)}$ with $\bl{f}^\intercal \bl{H} \bl{f} \geq 0$ on $\Delta_{k-1}$. If the test is SUB we require in addition that $\bl{f}^\intercal \bl{H} \bl{f} > 0$ on $\Pa$. If the $f_{(1)},\ldots,f_{(m)}$ are a Gr\"obner basis with respect to a graded monomial ordering then the representation of $\tilde{\beta}$ can be constructed so that $\mathrm{deg}(\tilde{\beta}) = \max_{i,j} \mathrm{deg}(h_{(ij)} f_{(i)}f_{(j)})$. The NTUB threshold is equal to the SOS threshold described in Theorem \ref{thm:UBThresholdSOSUpperBounds}. When the degrees of $f_{(1)},\ldots,f_{(m)}$ are all the same the NTUB threshold is also equal to the SUB threshold.
\end{customthm}
\begin{proof}
We assume that $f_{(1)},\ldots,f_{(m)}$ form a Gr\"obner basis with respect to a graded monomial order. This only comes into play when comparing degrees, and computing the UB thresholds. Following the proof of Theorem \ref{custthm:PrincipalIdealUBThresholdUMP}, we have $\tilde{\beta} = \sum_{i = 1}^m g_{(i)} f_{(i)}$ and $\nabla \tilde{\beta}(\bs{\pi}) = 0$ for $\bs{\pi} \in \text{int}(\tilde{\Delta}_{k-1}) \cap \Po$. When $f_{(1)},\ldots,f_{(m)}$ form a Gr\"obner basis with respect to a graded monomial order $\text{deg}(\tilde{\beta}) = \max_i \text{deg}(g_{(i)}f_{(i)})$ since this representation is obtained by polynomial long division.

    We compute
    \begin{align*}
       0 = \nabla \tilde{\beta}(\bs{\pi})  = \sum_{i = 1}^m \big(\nabla g_{(i)}(\bs{\pi}) f_{(i)}(\bs{\pi}) + \nabla f_{(i)}(\bs{\pi}) g_{(i)}(\bs{\pi}) \big) =  \begin{bmatrix}
            \nabla f_{(1)}(\bs{\pi}) &\cdots &\nabla f_{(m)}(\bs{\pi})
        \end{bmatrix} \begin{bmatrix}
            g_{(1)}(\bs{\pi})
            \\
            \vdots
            \\
            g_{(m)}(\bs{\pi})
        \end{bmatrix}.
    \end{align*}
    By the full-rank assumption this implies that $g_{(1)}(\bs{\pi}) = \cdots = g_{(m)}(\bs{\pi}) = 0$ on $\text{int}(\tilde{\Delta}_{k-1}) \cap \Po$. Thus, $g_{(i)} = \sum_{j = 1}^m h_{(ij)} f_{(j)} \in I_{\mb{R}^{k-1}}(\Po)$ for every $i$, where again we have $\text{deg}(g_{(i)}) = \max_j \text{deg}(h_{(ij)}f_{(j)})$. This yields the desired representation $\tilde{\beta} = \sum_{i = 1}^m\sum_{j = 1}^m h_{(ij)} f_{(i)} f_{(j)}$ with 
    \begin{align*}
        \text{deg}(\tilde{\beta}) = \max_i\text{deg}(g_{(i)})\text{deg}(f_{(i)}) = \max_{i,j}\text{deg}(h_{(ij)}f_{(j)})\text{deg}(f_{(i)}) = \max_{i,j} \text{deg}(h_{(ij)}f_{(i)}f_{(j)})
    \end{align*}
    under the Gr\"obner basis assumption.
    The non-negativity and positivity requirements on $\bl{f}^\intercal \bl{H} \bl{f}$ are immediate. 

The sample size required for any UB test with separating polynomial $\tilde{\beta}$ is $\max_{i,j} \text{deg}(h_{(ij)} f_{(i)}f_{(j)})$, which is bounded below by $2 \min_i \text{deg}(f_{(i)})$ for all $\bl{H} \neq \bl{0}$. This lower bound is attained by the NTUB test with the separating polynomial $f_{(i^*)}^2$, where $f_{(i^*)}$ is the smallest degree polynomial in the Gr\"obner basis. As the SUB threshold is bounded below by the NTUB threshold, when all of the $f_{(i)}$s have the same degree $\sum_{i = 1}^m f_{(i)}^2$ is a SUB separating polynomial with minimum total degree.

\end{proof}

%%%%%%%%%%%%%%%%%%%%%%%%%%%%%%%%%%%%%%%%%%%%%
%%%%%%%%%%%%%%%%%%%%%%%%%%%%%%%%%%%%%%%%%%%%%

\begin{customthm}{9}
\label{custthm:PolytopePeelingUMPU}
A necessary condition for a level-$\alpha$  UMPU test for $\Po$ to exist is that
    \begin{enumerate}
        \item For $V^{(0)} = \{I_{01},\ldots,I_{0l_0}\}$ the set $\text{Proj}_{(I_{01},\ldots,I_{0l_0})}\big( \mc{C}_{n,\alpha}(\Po)\big)$ has a componentwise maximum element $(h_{I_{01}}^*,\ldots,h_{I_{0l_0}}^*)$. The vertex set $V^{(0)}$ is described in Definition \ref{def:ConvexPeeling}.
        \item Inductively, for every $j$ with $V^{(j)} = \{I_{j1},\ldots,I_{jl_j}\}$ the set
        \begin{align*}
\text{Proj}_{(I_{j1},\ldots,I_{jl_j})}\big( \mc{C}_{n,\alpha}(\Po) \cap \{ (h_I) : h_{I_{j'a}} = h_{I_{j'a}}^*, \forall j' < j, a \leq l_{j'}\}\big)
        \end{align*}set 
has a componentwise maximum $(h^*_{I_{j1}},\ldots,h^*_{I_{jl_j}})$. 
    \end{enumerate}
    If a UMPU test exists it has the power polynomial $f^2h^* + \alpha$, where $h^*$ is the polynomial  corresponding to the above polynomial coefficients $(h^*_I)$.
\end{customthm}
\begin{proof}
    Assume that there exists a UMPU test with power polynomial $\beta = \Tilde{f}^2 h + \alpha$ and $h$ homogeneous. Take a different power polynomial $\gamma = \Tilde{f}^2 g + \alpha$ of the same form. The coefficients of $h$ and $g$ must lie in the coefficient polytope. We have $\beta - \gamma = \Tilde{f}^2(h - g) \geq 0$ on $\Pa$ by assumption. By continuity this implies that $h - g \geq 0$ on $\Delta_{k-1}$. Evaluate $h(\bl{e}_w) - g(\bl{e}_w) = h_{I_{0w}} - g_{I_{0w}} \geq 0$ at the standard basis vector $\bl{e}_w$. This shows that the coefficients of $h$ satisfy $(h_{I_{01}},\ldots,h_{I_{0k}}) = (h^*_{I_{01}},\ldots,h_{I_{0k}}^*)$, verifying condition 1.

 Inductively, assume that conditions 1. and 2. hold up to $V^{(j-1)}$. It follows, by the assumption that $(h_I)$ corresponds to a UMPU test, that
 \begin{align*}
    \sum_{I \in \mc{S}^{(j)}}(h_I - g_I) \pi^I 
 =  \sum_I (h_I - g_I) \bs{\pi}^I  \geq 0 
 \end{align*}
 on $\Delta_{k-1}$, and hence by homogeneity on $\mb{R}_{\geq 0}^k$,
 for any set of coefficients $(g_I)$ satisfying the constraint
 \begin{align}
\label{eqn:ProjectedPolytopeConstraint}
(g_I)_{I \in \mc{S}^{(j)}} \in     \text{Proj}_{I \in \mc{S}^{(j)}}\big( \mc{C}_{n,\alpha}(\Po) \cap \{ (h_I) : h_{I_{j'a}} = h_{I_{j'a}}^*, \forall j' < j, a \leq l_{j'}\}\big).
\end{align}
 Let $A \in V^{(j)}$ be a vertex of $\Delta^{(j)}_{k-1}$, meaning that there exists a vector $\bl{x}$ with $ (a_1,\ldots,a_k) \bl{x} > I^\intercal \bl{x}$ for all $I \neq (a_1,\ldots,a_k)$ in  $\Delta^{(j)}_{k-1}$. Defining the sequence $\bl{v}_t \coloneqq (t^{x_1},\ldots,t^{x_k}) \in \mb{R}_{\geq 0}^k$ we compute
 \begin{align*}
      \sum_{I \in \mc{S}^{(j)}}(h_I - g_I) \bl{v}_t^I =   \sum_{I \in \mc{S}^{(j)}}(h_I - g_I) t^{\bl{x}^\intercal I} \geq 0. 
 \end{align*}
 As $t \rightarrow \infty$ it cannot be the case that $h_A < g_A$ as otherwise the sign of the above expression would be negative for large $t$. We conclude that for every $A \in V^{(i)}$, $h_A$ must equal the maximum value of $g_A$, namely $h_A^*$, as needed.
\end{proof}

%%%%%%%%%%%%%%%%%%%%%%%%%%%%%%%%%%%%%%%%%%%%%%
%%%%%%%%%%%%%%%%%%%%%%%%%%%%%%%%%%%%%%%%%%%%%%

\begin{customlem}{7}
\label{custlem:UMPUSufficientCondition}
    If $\bl{v}_1,\ldots,\bl{v}_m$ are the vertices of $\mc{C}_{n,\alpha}(\Po)$ and there exists a vertex that is the maximum in the componentwise ordering, then there exists a UMPU test. When $n' = 1$ this is a also a necessary condition for the existence of a UMPU test.
\end{customlem}
\begin{proof}
If $(h_I)$ is a componentwise maximum in the coefficient polytope then $h(\bs{\pi}) \geq g(\bs{\pi})$ on $\Delta_{k-1}$ for any $(g_I)$ in the coefficient polytope. This implies $\tilde{f}^2 h + \alpha \geq \tilde{f}^2 g + \alpha$, showing that $(h_I)$ corresponds to a UMPU test. It remains to check that $h$ is non-negative on $\Delta_{k-1}$, so that the corresponding test is unbiased. As $\bl{0}$ is in the coefficient polytope, $h_I \geq 0$ for every $I$ and $h(\bs{\pi}) \geq 0$ on $\Delta_{k-1}$. When $n' = 1$, $(h_I)$ being a componentwise maximum is equivalent to the necessary condition in Theorem \ref{custthm:PolytopePeelingUMPU}.  
\end{proof}

%%%%%%%%%%%%%%%%%%%%%%%%%%%%%%%%%%%%%%
%%%%%%%%%%%%%%%%%%%%%%%%%%%%%%%%%%%%%%

\begin{customlem}{8}
    The power polynomial $\beta$ of every unbiased, invariant test of $\Po = \{\bs{\pi} \in \Delta_{p^2-1}: \pi_{ij} = \pi_{ji}\}$ can be written in the form
    \begin{align*}
        \beta(\bs{\pi}) = \sum_{i,j = 1}^p (\pi_{ij} - \pi_{ji})^2 h_{(ij)} + \alpha,
    \end{align*}
    where the polynomial $h_{(ij)}$ is invariant under the transposition $\pi_{kl} \mapsto \pi_{lk}$ for all $l,k \in [p]$.
\end{customlem}
\begin{proof}
Using Theorem \ref{thm:NonPrincipalPowerPolyRepresentation}, $\beta$ has a representation of the form $\beta(\bs{\pi)} = \sum_{i < j}\sum_{k < l}(\pi_{ij} - \pi_{ji})(\pi_{kl} - \pi_{lk}) h_{(ijkl)}(\bs{\pi}) + \alpha$. As $\beta$ is assumed to be invariant under the transposition $\tau$ that sends $\pi_{ab} \mapsto \pi_{ba}$ we have $S_{ab}\big(\beta(\bs{\pi})\big) \coloneqq \tfrac{1}{2}\big(\beta(\bs{\pi}) + \beta(\tau(\bs{\pi}))\big) = \beta(\bs{\pi})$. We can apply the linear symmetrization operator $S_{ab}$ term by term. For any term $(\pi_{ij} - \pi_{ji})(\pi_{kl} - \pi_{lk})h_{(ijkl)}$ where $\pi_{ab}$ does not appear in the first two factors we get $S_{ab}((\pi_{ij} - \pi_{ji})(\pi_{kl} - \pi_{lk})h_{(ijkl)}) = (\pi_{ij} - \pi_{ji})(\pi_{kl} - \pi_{lk})S_{ab}(h_{(ijkl)})$. Similarly, $S_{ab}( (\pi_{ab} - \pi_{ba})^2 h_{(abab)})  = (\pi_{ab} - \pi_{ba})^2 S_{ab}(h_{(abab)})$. The last remaining case is if $\pi_{ab}$ appears in one of the factors, but not the other. Assume without loss of generality that $\pi_{ab}$ appears in the first factor. Then we obtain
\begin{align*}
    S_{ab}\big((\pi_{ab} - \pi_{ba})(\pi_{kl} - \pi_{lk}) h_{(abkl)}(\bs{\pi})\big) = (\pi_{ab} - \pi_{ba})(\pi_{kl} - \pi_{lk}) \big((h_{(abkl)}(\bs{\pi}) - h_{(abkl)}(\bs{\tau}(\bs{\pi}))\big).
\end{align*}
The above computations show that we can write $\beta(\bs{\pi)} = \sum_{i < j}\sum_{k < l}(\pi_{ij} - \pi_{ji})(\pi_{kl} - \pi_{lk}) \tilde{h}_{(ijkl)}(\bs{\pi}) + \alpha$ where $\tilde{h}_{(ijkl)} = S_{ab}(\tilde{h}_{(ijkl)})$ in the first two cases and $\tilde{h}_{(abkl)}(\bs{\pi}) = h_{(abkl)}(\bs{\pi}) - h_{(abkl)}(\bs{\tau}(\bs{\pi}))$ in the last case. In this last case $\tilde{h}_{(abkl)}$ is antisymmetric with respect to $\tau$ and thus this polynomial vanishes along the subspace $\{\bs{\pi}: \pi_{ab} = \pi_{ba}\}$. We conclude that $(\pi_{ab} - \pi_{ba})$ divides $\tilde{h}_{(abkl)}$ and we have a representation of the form $(\pi_{ab} - \pi_{ba})(\pi_{kl} - \pi_{lk}) \tilde{h}_{(abkl)} = (\pi_{ab} - \pi_{ba})^2(\pi_{kl} - \pi_{lk}) \hat{h}_{(abkl)}$ where $\tilde{h}_{(abkl)} =(\pi_{ab} - \pi_{ba})\hat{h}_{(abkl)}$. Next, symmetrize $\beta$ with respect to $S_{ab}$ for every possible pair $a < b$. After symmetrization, any term $(\pi_{ij} - \pi_{ji})(\pi_{kl} - \pi_{lk}) h_{(ijkl)}$ with $(i,j) \neq (k,l)$ can be written in the form $(\pi_{ij} - \pi_{ji})^2\big((\pi_{kl} - \pi_{lk}) \hat{h}_{(ijkl)}(\bs{\pi})\big) = (\pi_{ij} - \pi_{ji})^2 \Bar{h}_{(ijkl)}$ where $\Bar{h}_{(ijkl)}= (\pi_{kl} - \pi_{lk}) \hat{h}_{(ijkl)}$. Grouping the terms that involve a square together, we conclude that $\beta$ has a representation of the form 
\begin{align*}
    \beta(\bs{\pi}) = \sum_{i < j} (\pi_{ij} - \pi_{ji})^2 h_{(ij)}(\bs{\pi}) + \alpha.
\end{align*}
for some choice of polynomials $h_{(ij)}$. Each $h_{(ij)}$ can be chosen to be symmetric with respect to the transposition $\tau$ by another round of symmetrization with respect to $S_{ab}$ if needed.
\end{proof}

\section*{Derivations for Examples 3, 5, and 7}

\subsection*{Example 3 (Square Null Hypothesis)}
We compute the asymptotic properties of the test based on $\max(x_1,x_2)$. Fixing a $c_0$, the power function $\text{Pr}_{\bs{\pi}}( \max(x_1,x_2) > c_0)$ is a monotone increasing function of $(\pi_1,\pi_2)$. Thus, the power function of such a test is maximized over the null hypothesis set $\mc{P}_{0,t}$ when $\bs{\pi} = (t,t,1-2t)$. The central limit theorem applied to $\hat{\pi}_i = x_i/n$ for $i = 1,2,3$ will be used to find the correct $c_0$ to control the level of the test. By the central limit theorem
\begin{align*}
    \sqrt{n}(\hat{\pi}_1 - \pi_1,\hat{\pi}_2- \pi_2) \sim \mc{N}\bigg(\bl{0}, \begin{bmatrix}
        \pi_1(1-\pi_1) & -\pi_1 \pi_2
        \\
        \pi_2\pi_1 & \pi_2(1-\pi_2)
    \end{bmatrix}\bigg).
\end{align*}
Define $(z_1,z_2)$ to be normal random variables distributed according to this asymptotic distribution and take $\bs{\pi}_0 = (t,t,1-2t)$. Then 

\begin{align*}
    \text{Pr}_{\bs{\pi}_0}( \max(x_1,x_2) > nt + \sqrt{n}c) & = \text{Pr}_{\bs{\pi}_0}( \sqrt{n}\max(\hat{\pi}_1 - t,\hat{\pi}_2 - t) > c)
    \\
    & = 1 - \text{Pr}( (z_1,z_2) \in (-\infty,c]^2 ).
\end{align*}
Numerically solving $\text{Pr}( (z_1,z_2) \in (-\infty,c]^2 ) = 1-\alpha$ determines $c$. To show that this test is consistent, assume without loss of generality that $\pi_1 > t$. Then 
\begin{align*}
        \text{Pr}_{\bs{\pi}}( \max(x_1,x_2) > nt + \sqrt{n}c) & \geq \text{Pr}_{\bs{\pi}}( \sqrt{n}(\hat{\pi}_1 - t) > c)
        \\
        & = \text{Pr}_{\bs{\pi}}( \sqrt{n}(\hat{\pi}_1 - \pi_1) > c + \sqrt{n}(\pi_1 - t) ),
\end{align*}
with the last term converging to $1$ as $n \rightarrow \infty$ by the central limit theorem. 

The asymptotic power at $\bs{\pi}_* = (t,0,1-t)$ can be computed as
\begin{align*}
      \text{Pr}_{\bs{\pi}_*}( \max(x_1,x_2) > nt + \sqrt{n}c) & =   \text{Pr}_{\bs{\pi}_*}( \sqrt{n}(\hat{\pi}_1 - t) > c) 
      \\
      & = \text{Pr}(z_1 > c) = 1 - \text{Pr}(z_1 \in (-\infty, c]).
\end{align*}
As
\begin{align*}
    \text{Pr}(z_1 \in (-\infty, c]) =   \text{Pr}((z_1,z_2) \in (-\infty, c] \times (-\infty,\infty) )  > \text{Pr}( (z_1,z_2) \in (-\infty,c]^2 )
\end{align*}
the asymptotic power of the test at $\bs{\pi}_*$ is strictly less than the power at $\bs{\pi}_0$. 

%%%%%%%%%%%%%%%%%%%%%%%%%%%%%%%%%%%%%%%%%%%%%
%%%%%%%%%%%%%%%%%%%%%%%%%%%%%%%%%%%%%%%%%%%%%

\subsection*{Example 5 (Independence Null Hypothesis with an Extra Equality Constraint)}

Macaulay 2 code for computing the Gr\"obner basis in Example \ref{ex:IndepModelpq}:
\begin{verbatim}
M2
R = QQ[p_(1,1)..p_(2,3)]

-- Independence ideal for a 2x3 table where we assume that the
-- first row marginal probability is equal to twice marginal
-- probability of the first column.

-- Substitute this expression for p_(2,3)
p23 = 1- p_(1,1) - p_(1,2) - p_(1,3) - p_(2,1) - p_(2,2)

I = ideal( p_(1,1)* p_(2,2) - p_(1,2)*p_(2,1),  p_(1,1)*p23 - p_(1,3)*p_(2,1),
    2*p_(1,1) + 2*p_(2,1) - p_(1,1) - p_(1,2) - p_(1,3)) 

Irad  = radical I

netList(flatten(entries(gens gb Irad)))
\end{verbatim}

%%%%%%%%%%%%%%%%%%%%%%%%%%%%%%%%%%%%%%%%%%%
%%%%%%%%%%%%%%%%%%%%%%%%%%%%%%%%%%%%%%%%%%%

\subsection*{Example 7 (Spherical Null Hypothesis)}
Viewed as a subset of $\mb{R}^{k-1}$, the null hypothesis set is a $(k-2)$-dimensional sphere as it is the intersection of the hyperplane $\text{aff}(\Delta_{k-1})$ with a $(k-1)$-dimensional sphere in $\mb{R}^k$ centered at $\tfrac{1}{k}\bl{1}_k$. The stereographic projection described in the main text also provides a rational parameterization of the $(k-2)$-dimensional sphere. Lemma \ref{custlem:RationalParamIdeal} applies with $\varphi(\mb{R}^{k-2}) \subseteq \Po$. The polynomial $\sum_{i = 1}^{k-1} (\pi_i-\tfrac{1}{k})^2 + (1 - \sum_{i = 1}^{k-1}\pi_i - \tfrac{1}{k})^2 - \delta^2$ is irreducible since it has no linear factors as no linear forms vanish on the sphere. Thus, $\sqrt[\mb{C}]{\langle \sum_{i = 1}^{k-1} (\pi_i-\tfrac{1}{k})^2 + (1 - \sum_{i = 1}^{k-1}\pi_i - \tfrac{1}{k})^2 - \delta^2 \rangle} = \langle \sum_{i = 1}^{k-1} (\pi_i-\tfrac{1}{k})^2 + (1 - \sum_{i = 1}^{k-1}\pi_i - \tfrac{1}{k})^2 - \delta^2 \rangle$, proving the claim in the main text. 

We now compute $c_\alpha$ and provide an expression for the corresponding UMPU test. Define $\gamma = (\delta + \tfrac{1}{k})$ to ease notation. The power polynomial has the form 
\begin{align*}
    \beta(\bs{\pi}) &=  \bigg(c_\alpha\sum_{i,j = 1}^k \pi_i^2 \pi_j^2\bigg) -2c_\alpha \gamma \bigg(\sum_{i = 1}^k \pi_i^2\bigg) \bigg(\sum_{j = 1}^k \pi_j\bigg)^2+ (c_\alpha \gamma^2 + \alpha) \bigg(\sum_{i = 1}^k \pi_i\bigg)^4.
\end{align*}
We assume for simplicity that $k \geq 4$. The coefficients of each monomial $\bs{\pi}^I$ in $\beta$ are:
\begin{align*}
    \begin{cases}
        c_\alpha -2c_\alpha \gamma + (c_\alpha \gamma^2 + \alpha) \;\; &\text{ if } I = (i,i,i,i)
        \\
       -4c_\alpha \gamma + 4(c_\alpha\gamma^2 + \alpha)  \;\; &\text{ if } I = (i,i,i,j)
         \\
                 2c_\alpha - 4c_\alpha \gamma + 6(c_\alpha \gamma^2 + \alpha)  \;\; &\text{ if } I = (i,i,j,j)
         \\
      -4\gamma c_\alpha + 12(c_\alpha\gamma^2 + \alpha)  \;\; &\text{ if } I = (i,i,j,k)
          \\
       24(c_\alpha\gamma^2 + \alpha)   \;\; &\text{ if } I = (i,j,k,l)
    \end{cases},
\end{align*}
where $i,j,k,l \in [k]$ are distinct indices.  
If $\beta$ is to be a power polynomial it must satisfy the following inequalities
\begin{align*}
   &0 \leq (1 -2\gamma + \gamma^2) c_\alpha + \alpha  \leq 1
   \\
   &0 \leq (4\gamma^2 - 4\gamma)c_\alpha + 4\alpha \leq 4
   \\
   &0 \leq (2 - 4\gamma + 6\gamma^2)c_\alpha + 6\alpha \leq 6
   \\
   &0 \leq (12\gamma^2 - 4 \gamma)c_\alpha + 12\alpha \leq 12
   \\
   &0 \leq 24\gamma^2 c_\alpha + 24\alpha \leq 24.
\end{align*}
Note that $\delta$ must be less than the distance between the center point $\tfrac{1}{k}\bl{1}_k$ and a vertex say $\bl{e}_1 = (1,0,\ldots,0)$, which is $1 - \tfrac{1}{k}$. This implies that $0 < \gamma < 1$. Each inequality will enforce only a single constraint on $c_\alpha$, with the constraint determined by the sign of the coefficient of $c_\alpha$ in the above inequalities. For instance, the last inequality adds the constraint that $c_\alpha \leq \frac{1-\alpha}{\gamma^2}$. However, the sign of some of the coefficients of $c_\alpha$ will depend on the choice of $\delta$ and $k$. Let's assume that $\delta = \tfrac{1}{2} - \tfrac{1}{k}$ so that $\gamma = \tfrac{1}{2}$. Other choices of $\delta$ can be handled in an analogous manner. In this case the inequalities imply that  $c_\alpha = \min\{4(1-\alpha),4\alpha\} = 4\alpha$, assuming that $\alpha \leq \tfrac{1}{2}$. The hypothesis test with power polynomial $\beta$ is equal to
\begin{align*}
    \phi(\bs{x}) =     \begin{cases}
  2\alpha \;\; &\text{ if } I(\bl{x}) = \{i,i,i,i\}
        \\
     0  \;\; &\text{ if } I(\bl{x}) = \{i,i,i,j\}
         \\
     2\alpha  \;\; &\text{ if } I(\bl{x}) = \{i,i,j,j\}
         \\
    \tfrac{4}{3}\alpha \;\; &\text{ if } I(\bl{x}) = \{i,i,j,k\}
          \\
       2\alpha  \;\; &\text{ if } I(\bl{x}) = \{i,j,k,l\}
    \end{cases},
\end{align*}
where the notation $I(\bl{x})$ denotes the multiset consisting of the categories of the counts with repetition. For example, if the count vector has a count of four in the single category $i$, the corresponding multiset is $\{i,i,i,i\}$. If the $\bl{x}$ has a count of $2$ in category $i$ and counts of $1$ in categories $j$ and $k$ then $I(\bl{x}) = \{i,i,j,k\}$.

The coefficient polytope vertices for $n = 5$ are 
\begin{align*}
\begin{bmatrix}
    -0.2
    \\
    0.2
    \\
    -0.2
\end{bmatrix},
\begin{bmatrix}
    0.2
    \\
    0.2
    \\
    -0.2
\end{bmatrix},
\begin{bmatrix}
    -0.2
    \\
    -0.2
    \\
    -0.2
\end{bmatrix},
\begin{bmatrix}
    0.2
    \\
    -0.2
    \\
    -0.2
\end{bmatrix},
\begin{bmatrix}
    0.2
    \\
    -0.2
    \\
    0.2
\end{bmatrix},
\begin{bmatrix}
    -0.2
    \\
    0.2
    \\
    0.2
\end{bmatrix},
\begin{bmatrix}
    -0.2
    \\
    -0.2
    \\
    0.2
\end{bmatrix},
\begin{bmatrix}
    1/3
    \\
    1/3
    \\
    1/3
\end{bmatrix},
%     -0.2        0.2       -0.2
 % 0.2        0.2       -0.2
% -0.2       -0.2       -0.2
%  0.2       -0.2       -0.2
 % 0.2       -0.2        0.2
% -0.2        0.2        0.2
% -0.2       -0.2        0.2
%  0.333333   0.333333   0.333333
\end{align*}
with the last vertex clearly being the componentwise maximum. When $n = 6$ the coefficient polytope has $188$ vertices. The vertices that maximize each coordinate are displayed below
\begin{align*}
\begin{bmatrix}
    -0.2
    \\
    0.56
    \\
    0.44
    \\
    0.56
    \\
    1.52
    \\
    0.44
\end{bmatrix},
\begin{bmatrix}
    0.44
    \\
    0.56
    \\
    0.44
    \\
    0.56
    \\
    1.52
    \\
    0.44
\end{bmatrix},
\begin{bmatrix}
    0.44
    \\
    1.52
    \\
    0.44
    \\
    0.56
    \\
    0.56
    \\
    -0.2
\end{bmatrix},
\begin{bmatrix}
    0.44
    \\
    1.52
    \\
    0.44
    \\
    0.56
    \\
    0.56
    \\
    0.44
\end{bmatrix},
\begin{bmatrix}
    0.44
    \\
    0.56
    \\
    0.44
    \\
    1.52
    \\
    0.56
    \\
    0.44
\end{bmatrix},
\begin{bmatrix}
    0.44
    \\
    0.56
    \\
    -0.2
    \\
    1.52
    \\
    0.56
    \\
    0.44
\end{bmatrix},
\begin{bmatrix}
    0.6
    \\
    1.2
    \\
    0.6
    \\
    1.2
    \\
    1.2
    \\
    0.6
\end{bmatrix}.
%     -0.2   0.56   0.44  0.56  1.52   0.44
%  0.44  0.56   0.44  0.56  1.52   0.44
%  0.44  1.52   0.44  0.56  0.56  -0.2
%  0.44  1.52   0.44  0.56  0.56   0.44
%  0.44  0.56   0.44  1.52  0.56   0.44
 % 0.44  0.56  -0.2   1.52  0.56   0.44
 % 0.6   1.2    0.6   1.2   1.2    0.6
\end{align*}
The coordinates of these vertices are ordered so that $h_i$ corresponds to the $i$th coordinate, where $h(\bs{\pi}) = h_1\pi_1^2 + h_2\pi_1\pi_2 + h_3\pi_2^2 + h_4\pi_1\pi_3 + h_5\pi_2\pi_3 + h_6\pi_3^2$. The first, third, and sixth coordinates are uniquely maximized at the last vertex. Thus, if a UMPU test exists it must have an $h$ of the form $h(\bs{\pi}) = 0.6(\pi_1^2 + \pi_2^2 + \pi_3^2) + 1.2(\pi_1\pi_2 + \pi_1\pi_3 + \pi_2\pi_3)$.  

%%%%%%%%%%%%%%%%%%%%%%%%%%%%%%%%%%%%%%%%%%%%%%%%%%%%%%%%%
%%%%%%%%%%%%%%%%%%%%%%%%%%%%%%%%%%%%%%%%%%%%%%%%%%%%%%%%

\end{document}